\documentclass[a4paper,11pt]{article}

\usepackage{amsmath}
\usepackage{amsfonts}
\usepackage{amssymb}
\usepackage{latexsym}
\usepackage{epsfig}
\usepackage{graphicx}
\usepackage{oldgerm}
\usepackage{theorem}
\usepackage{bm}
\usepackage{mathrsfs}

\def\la{\lambda}
\def\La{\Lambda}
\def\C{\mathbb{C}}
\def\R{\mathbb{R}}
\def\N{\mathbb{N}}
\def\Z{\mathbb{Z}}
\def\H{\mathbb{H}}
\def\S{\mathbb{S}}
\def\B{\mathbb{B}}

\def\I{\mathbb{I}}
\def\W{\mathbb{W}}
\def\D{\mathbb{D}}
\def\O{\mathbb{O}}
\def\P{\mathbb{P}}
\def\E{\mathbb{E}}

\def\J{\mathbb{J}}
\def\ker{{\rm ker} \,}
\def\Ai{{\rm Ai}}
\def\gO{{\rm O}}
\def\U{{\rm U}}
\def\SO{{\rm SO}}
\def\Sp{{\rm Sp}}
\def\rH{{\rm H}}
\def\rS{{\rm S}}
\def\Conf{{\rm Conf}}
\def\supp{{\rm supp}\,}
\def\hcap{{\rm hcap}}
\def\rR{{\rm R}}
\def\rI{{\rm I}}
\def\rN{{\rm N}}
\def\rP{{\rm P}}
\def\rO{{\rm O}}
\def\rE{{\rm E}}
\def\re{{\rm e}}
\def\bE{{\bf E}}
\def\bP{{\bf P}}
\def\bp{{\bf p}}
\def\bK{{\bf K}}
\def\x{\bm{x}}
\def\y{\bm{y}}

\def\u{\bm{u}}
\def\t{\bm{t}}
\def\f{\bm{f}}
\def\bv{\bm{v}}
\def\bLambda{\bm{\Lambda}}
\def\blambda{\bm{\lambda}}
\def\bmeta{\bm{\eta}}

\def\X{\bm{X}}
\def\Y{\bm{Y}}
\def\1{{\bf 1}}
\def\0{{\bf 0}}
\def\V{\bm{V}}

\def\RN{{\it R}_N}
\def\AN{{\it A}_{N-1}}
\def\BN{{\it B}_N}

\def\CN{{\it C}_N}

\def\DN{{\it D}_N}

\def\cC{{\cal C}}
\def\cS{{\cal S}}
\def\cP{{\cal P}}
\def\cH{{\cal H}}

\def\cK{{\cal K}}
\def\cN{{\cal N}}
\def\cW{{\cal W}}
\def\cF{{\cal F}}
\def\cI{{\cal I}} 
\def\cA{{\cal A}}
\def\cB{{\cal B}}
\def\cU{{\cal U}}
\def\cM{{\cal M}}
\def\cV{{\cal V}}
\def\cD{{\cal D}}
\def\cE{{\cal E}}
\def\cN{{\cal N}}
\def\cG{{\cal G}}

\def\sS{{\sf S}}

\def\sD{{\sf D}}
\def\sm{{\sf m}}
\def\sK{{\sf K}}
\def\law={\stackrel{\rm (law)}{=}}
\def\weak{\stackrel{N \to \infty}{\, \Longrightarrow \,}}

\def\weakp{\stackrel{p \to \infty}{\, \Longrightarrow \,}}
\def\Det{\mathop{\mathrm{Det}}}
\def\tr{{\rm Tr} \,}
\def\locally{\mathcal{I}_{1,\mathrm{loc}}}

\def\etabar{\overline{\eta}}
\def\zetabar{\overline{\zeta}}
\def\Re{{\rm Re} \,}
\def\Im{{\rm Im} \,}
\def\arg{{\rm arg}\,}
\def\mS{\mathfrak{S}}
\def\bra{\langle}
\def\ket{\rangle}

\theorembodyfont{\itshape}

\newtheorem{thm}{Theorem}[section]
\newtheorem{lem}[thm]{Lemma}
\newtheorem{cor}[thm]{Corollary}
\newtheorem{prop}[thm]{Proposition}

\newtheorem{df}[thm]{Definition}

\newtheorem{example}[thm]{Example}
\newtheorem{rem}[thm]{Remark}

\newcommand{\SSC}[1]{\section{#1}\setcounter{equation}{0}}
\newcommand{\qed}{\hbox{\rule[-2pt]{3pt}{6pt}}}

\begin{document}

\title{\bf
Determinantal Point Processes, \\
Stochastic Log-Gases, 
and Beyond
\footnote{
This manuscript was prepared for the 
mini course given at 
`Workshop on Probability and Stochastic Processes'
held at Orange County, Coorg, India, 
from 23rd to 26th February, 2020,
which was organized by 
the Indian Academy of Sciences, Bangalore.
}
}
\author{
Makoto Katori
\footnote{
Department of Physics,
Faculty of Science and Engineering,
Chuo University,
Kasuga, Bunkyo-ku, Tokyo 112-8551, Japan;
e-mail: katori@phys.chuo-u.ac.jp}
}
\date{7 March 2020}
\pagestyle{plain}
\maketitle

\begin{abstract}

A determinantal point process (DPP) is an ensemble
of random nonnegative-integer-valued Radon measures,
whose correlation functions are all given by determinants
specified by an integral kernel called the correlation kernel.
First we show our new scheme of DPPs in which 
a notion of partial isometies between a pair of Hilbert spaces
plays an important role.
Many examples of DPPs in one-, two-, and higher-dimensional
spaces are demonstrated, where several types
of weak convergence from finite DPPs to
infinite DPPs are given.
Dynamical extensions of DPP are realized in one-dimensional
systems of diffusive particles conditioned
never to collide with each other.
They are regarded as one-dimensional stochastic log-gases,
or the two-dimensional Coulomb gases 
confined in one-dimensional spaces.
In the second section, we consider such 
interacting particle systems in one dimension. 
We introduce a notion of 
determinantal martingale and prove that,
if the system has determinantal martingale representation (DMR),
then it is a determinantal stochastic process (DSP) in the sense that
all spatio-temporal correlation function are expressed by
a determinant.
In the last section, we construct processes of
Gaussian free fields (GFFs) on simply connected 
proper subdomains of $\C$ coupled with 
interacting particle systems defined 
on boundaries of the domains.
There we use multiple Schramm--Loewner evolutions (SLEs)
driven by the interacting particle systems.
We prove that, if the driving processes
are time-changes of the log-gases studied in the
second section, then the obtained GFF with multiple SLEs
are stationary.
The stationarity defines an equivalence relation of
GFFs, which will be regarded as a generalization of
the imaginary surface studied by Miller and Sheffield.

\vskip 0.2cm


\end{abstract}
\footnotesize
\tableofcontents
\vspace{3mm}
\normalsize

\SSC
{Determinantal Point Processes (DPPs)}
\label{sec:DPP}

A determinantal point process (DPP) is 
an ensemble of random nonnegative-integer-valued 
Radon measures $\Xi$ on 
a space $S$ with measure $\lambda$,
whose correlation functions are all
given by determinants specified by an
integral kernel $K$ called the correlation kernel.
We consider a pair of Hilbert spaces, 
$H_{\ell}, \ell=1,2$, which are assumed to
be realized as $L^2$-spaces,
$L^2(S_{\ell}, \lambda_{\ell})$, $\ell=1,2$, 
and introduce 
a bounded linear operator ${\cal W} : H_1 \to H_2$
and its adjoint ${\cal W}^{\ast} : H_2 \to H_1$.
We show that if ${\cal W}$ is a partial isometry of
locally Hilbert-Schmidt class, then we have a unique DPP on
$(\Xi_1, K_1, \lambda_1)$ associated with $\cW^* \cW$. 
In addition, if $\cW^*$ is also 
of locally Hilbert--Schmidt class, then we have a unique pair
of DPPs, $(\Xi_{\ell}, K_{\ell}, \lambda_{\ell})$, $\ell=1,2$.

We also give a practical framework 
which makes ${\cal W}$ and ${\cal W}^{\ast}$ satisfy the above conditions.
Our framework to construct pairs of DPPs implies
useful duality relations between DPPs making pairs.
For a correlation kernel of a given DPP
our formula can provide plural different expressions, 
which reveal different aspects of the DPP.

In order to demonstrate these advantages of 
our framework as well as to show that
the class of DPPs obtained by this method is large enough to
study universal structures in a variety of DPPs, 
we report plenty of examples of DPPs 
in one-, two-, and higher-dimensional spaces $S$,
where several types of weak convergence from finite DPPs
to infinite DPPs are given.

This section is based on the collaborations with
Tomoyuki Shirai (Kyushu University)
\cite{KS2}.

\subsection{Definition and existence theorem of DPP}
\label{thm:def_DPP}

Let $S$ be a base space, which is locally compact Hausdorff space
with countable base, 
and $\lambda$ be a Radon measure on $S$.
The configuration space over $S$ is given by
the set of nonnegative-integer-valued Radon measures; 
\[
{\rm Conf}(S)
=\left\{ \xi = \sum_j \delta_{x_j} : \mbox{$x_j \in S$,
$\xi(\Lambda) < \infty$ for all bounded set $\Lambda \subset S$} \right\}.
\]
${\rm Conf}(S)$ is equipped with the topological Borel $\sigma$-fields
with respect to the vague topology; 
we say $\xi_n, n \in \N :=\{1, 2, \dots\}$ converges
to $\xi$ in the vague topology, if
$\int_{S} f(x) \xi_n(dx) \to
\int_{S} f(x) \xi(dx)$,
$\forall f \in \cC_{\rm c}(S)$, 
where $\cC_{\rm c}(S)$ is the set of
all continuous real-valued functions 
with compact support.
A {\it point process} on $S$ is
a ${\rm Conf}(S)$-valued random variable
$\Xi=\Xi(\cdot, \omega)$ on a probability space
$(\Omega, \cF, \bP)$.
If $\Xi(\{x\}) \in \{0, 1\}$ for any point $x \in S$,
then the point process is said to be {\it simple}.

Assume that $\Lambda_j, j =1, \dots, m$, $m \in \N$ are
disjoint bounded sets in $S$ and
$k_j \in \N_0 :=\{0,1, \dots\}, j=1, \dots, m$ satisfy
$\sum_{j=1}^m k_j = n \in \N_0$.
A symmetric measure $\lambda^n$ on $S^n$
is called the $n$-th {\it correlation measure},
if it satisfies
\[
\bE \left[
\prod_{j=1}^m \frac{\Xi(\Lambda_j)!}
{(\Xi(\Lambda_j)-k_j)!} \right]
=\lambda^n(\Lambda_1^{k_1} \times \cdots
\times \Lambda_m^{k_m}),
\]
where if $\Xi(\Lambda_j)-k_j \leq 0$,
we interpret $\Xi(\Lambda_j)!/(\Xi(\Lambda_j)-k_j)!=0$.
If $\lambda^n$ is absolutely continuous 
with respect to the $n$-product measure $\lambda^{\otimes n}$,
the Radon--Nikodym derivative
$\rho^n(x_1, \dots, x_n)$ is called
the {\it $n$-point correlation function}
with respect to the background measure $\lambda$;
\[
\lambda^n(dx_1 \cdots dx_n)
=\rho^n(x_1, \dots, x_n) 
\lambda^{\otimes n}(dx_1 \cdots dx_n).
\]
Determinantal point process (DPP) 
is defined as follows \cite{Macchi75,Sos00,ST03a,ST03b,HKPV09}. 
\begin{df}
\label{thm:determinantal0}
A simple point process $\Xi$ on $(S, \lambda)$ is said to be
a determinantal point process (DPP) with 
correlation kernel $K : S \times S \to \C$ if it has correlation functions 
$\{\rho^n \}_{n \in \N}$, 
and they are given by
\begin{equation}
\rho^n(x_1, \dots, x_n) = \det_{1 \leq j, k \leq n}
[ K(x_j, x_k) ]
\quad \mbox{for every $n \in \N$, 
and $x_1, \dots, x_n \in S$}.
\label{eqn:DPP}
\end{equation}
The triplet $(\Xi, K, \lambda(d x))$ 
denotes the DPP, $\Xi \in {\rm Conf}(S)$,  
specified by the correlation kernel $K$
with respect to the measure $\lambda(d x)$.
\end{df}

If the integral projection operator $\cK$ on $L^2(S,\lambda)$ 
with a kernel $K$ is of rank $N \in \N$,
then the number of points is $N$ a.s.
If $N < \infty$ (resp. $N=\infty$),
we call the system a {\it finite DPP} (resp. an {\it infinite DPP}).
The density of points with respect to the background measure $\lambda(dx)$
is given by
\[
\rho(x) :=\rho^1(x) = K(x,x).
\]
The DPP is negatively correlated as shown by
\begin{align}
\rho^2(x, x') &= \det \left[
\begin{array}{ll}
K(x, x) & K(x, x') \cr
K(x', x) & K(x', x') 
\end{array}
\right]
\nonumber\\
&= K(x, x) K(x', x') - |K(x,x')|^2 \leq \rho(x) \rho(x'),
\quad x, x' \in S,
\label{eqn:rho2}
\end{align}
provided that $K$ is Hermitian.

Let $H$ be a separable Hilbert space. 
For operators $\cA, \cB$ on $H$, we say
that $\cA$ is positive definite and write 
$\cA \geq O$ if $\langle \cA f, f \rangle_H \geq 0$
for any $f \in H$, and write
$\cA \geq \cB$ if $\cA-\cB \geq O$. 
For any bounded operator $\cA$, the operator $\cA^* \cA$ is
positive definite. Then it admits a unique positive
definite square root
$\sqrt{\cA^* \cA}$ and is denoted by $|\cA|$. 
Let $\{\phi_n\}_{n \ge 1}$ be an orthonormal basis of $H$. 
For $\cA \ge O$, we define the trace of $\cA$ by 
\[
 \tr \cA := \sum_{n=1}^{\infty} \langle \cA \phi_n, \phi_n
 \rangle_H, 
\]
which does not depend on the choice of an orthonormal basis. An operator
$\cA$ is said to be \textit{of trace class} or a
\textit{trace class operator} if the trace norm $\|\cA\|_1 := \tr |\cA|$ is
finite. The trace $\tr \cA$ is defined whenever $\|\cA\|_1 < \infty$. 

Now, we consider the case $H = L^2(S, \la)$.
For a compact set $\Lambda \subset S$,
the projection from $L^2(S, \lambda)$ to 
the space of all functions vanishing outside $\Lambda$ $\lambda$-a.e.
is denoted by $\cP_{\Lambda}$.
$\cP_{\Lambda}$ is the operation of multiplication of
the indicator function ${\bf 1}_{\Lambda}$ of the set $\Lambda$;
${\bf 1}_{\Lambda}(x)=1$ if $x \in \Lambda$, and
${\bf 1}_{\Lambda}(x)=0$ otherwise.
We say that the bounded self-adjoint operator 
$\cA$ on $L^2(S, \lambda)$ is 
{\it of locally trace class} or a {\it locally trace class operator},
if the restriction of $\cA$ to each compact subset $\Lambda$,
is of trace class; that is, 
\begin{equation}
\mbox{
$\tr \cA_{\Lambda} < \infty$ \,
with \, $\cA_{\Lambda} :=\cP_{\Lambda} \cA \cP_{\Lambda}$ \,
for any compact set $\Lambda \subset S.$
}
\label{eqn:loctrace}
\end{equation}
The totality of locally trace
class operators on $L^2(S, \lambda)$ 
is denoted by $\locally(S, \lambda)$. 

Let $(S, \lambda)$ be a $\sigma$-finite measure space. 
We assume that $\cK \in \locally(S, \lambda)$.
If, in addition, $\cK \ge O$, then it  
admits a Hermitian integral kernel $K(x, x')$ such that (cf. \cite{GY05})
\begin{description}
\item{(i)} \quad
$\displaystyle{\det_{1 \leq j, k \leq n}[K(x_j, x_k)] \ge 0}$ for 
$\lambda^{\otimes n}$-a.e.\,$(x_1,\dots, x_n)$ for every $n \in \N$, 

\item{(ii)} \quad
$K_{x'} := K(\cdot, x') \in L^2(S, \lambda)$ for
$\lambda$-a.e.\,$x'$, 

\item{(iii)} \quad
$\tr \cK_{\Lambda} = \int_{\Lambda} K(x, x) \lambda(d x)$, 
$\Lambda \subset S$ and 
\[
 \tr (\cP_{\Lambda} \cK^n \cP_{\Lambda}) = \int_{\Lambda} \langle K_{x'}, \cK^{n-2} K_{x'}
 \rangle_{L^2(S, \lambda)} \lambda(d x'), 
\quad \forall n \in  \{2, 3, \dots\}. 
\]
\end{description}
The following is the {\it existence theorem of DPP}.
\begin{thm}[\cite{Sos00,ST03a,ST03b}]
\label{thm:exist_DPP}
Assume that $\cK \in \locally(S,\lambda)$ and $O \le \cK \le I$. 
Then there exists a
unique DPP $(\Xi, K, \la)$ on $S$.  
\end{thm}

If $\cK \in \locally(S, \lambda)$ is 
a projection onto a closed subspace 
$H \subset L^2(S, \lambda)$, 
one has {\it the DPP associated with $K$ and
$\lambda$}, or one may say 
{\it the DPP associated with the subspace $H$}. 

For $\cK$ its kernel space is denoted as $\ker \cK$
and the orthogonal complement of $\ker \cK$
is written as $(\ker \cK)^{\perp}$.
In this section, we consider the case that
\begin{align*}
&\cK f = f \quad
\mbox{for all $f \in (\ker \cK)^{\perp} \subset L^2(S, \lambda)$}
\nonumber\\
&\Longleftrightarrow \quad
\mbox{$\cK$ is an {\it orthogonal projection}}.
\end{align*}
By definition, it is obvious that the condition
$O \leq \cK \leq I$ is satisfied.
The purpose of the present section is to
introduce a useful method to provide
orthogonal projections $\cK$ and
DPPs whose correlation kernels are
given by the Hermitian integral kernels 
of $\cK$, $K(x, x'), x, x' \in S$.

\subsection{Partial isometries, 
locally Hilbert--Schmidt operators, and DPPs}
\label{sec:isometry}

Let $H_{\ell}, \ell=1,2$ be separable Hilbert spaces
with inner products $\langle \cdot, \cdot \rangle_{H_{\ell}}$.
For a bounded linear operator
$\cW : H_1 \to H_2$,
the adjoint of $\cW$ is defined as the operator
$\cW^{\ast} : H_2 \to H_1$, 
such that
\begin{equation}
\langle \cW f, g \rangle_{H_2} = \langle f, \cW^{\ast} g \rangle_{H_1}
\quad \mbox{for all $f \in H_1$ and $g \in H_2$}.
\label{eqn:ast1}
\end{equation}
A linear operator $\cW$ is called
an {\it isometry}  if
\[
\|\cW f\|_{H_2} = \|f\|_{H_1} \quad \mbox{for all $f \in H_1$}.
\]The kernel space of $\cW$ is denoted as $\ker \cW$
and its orthogonal complement is written as $(\ker \cW)^{\perp}$.
A linear operator $\cW$ is called 
a {\it partial isometry}, if 
\[
\|\cW f\|_{H_2} = \|f\|_{H_1} \quad \mbox{for all $f \in (\ker \cW)^{\perp}$}.
\]
For the partial isometry $\cW$, 
$(\ker \cW)^{\perp}$ is called the {\it initial space}
and the range of $\cW$, ${\rm ran} \cW$, is called the {\it final space}. 
By the definition (\ref{eqn:ast1}), 
$\|\cW f\|_{H_2}^2=\langle \cW f, \cW f \rangle_{H_2}
=\langle f, \cW^{\ast} \cW f \rangle_{H_1}$. As is suggested 
from this equality, we have the following fact for partial isometries. 
Although this might be known, we give a proof below. 
\begin{lem} \label{lem:A_Astar_partial_isometry2}
Let $H_1$ and $H_2$ 
be separable Hilbert spaces and $\cW : H_1 \to H_2 $ be a bounded
operator. Then, the following are equivalent. 
\begin{description}
\item{\upshape{(i)}} \quad $\cW$ is a partial isometry. 
\item{\upshape{(ii)}} \quad $\cW^* \cW$ is a projection on $H_1$, 
which acts as the identity on $(\ker \cW)^{\perp}$. 
\item{\upshape{(iii)}} \quad $\cW = \cW \cW^* \cW$. 
\end{description}
Moreover, $\cW$ is a partial isometry if and only if so is
 $\cW^*$. 
\end{lem}
\vskip 0.3cm
\noindent{\bf Assumption 1} \,
$\cW$ is a partial isometry.
\vskip 0.3cm

By Lemma~\ref{lem:A_Astar_partial_isometry2}, 
under Assumption 1, $\cW^*$ is also a partial
isometry and hence the operator $\cW^{\ast} \cW$ (resp. $\cW\cW^{\ast}$)
is the projection onto the initial space of $\cW$
(resp. the final space of $\cW$).

Now we assume that $H_1$ and $H_2$ are realized as $L^2$-spaces, 
$L^2(S_1, \lambda_1)$ and $L^2(S_2, \lambda_2)$,
respectively. 

A bounded linear operator $\cA : L^2(S_1, \lambda_1) \to
L^2(S_2, \lambda_2)$ 
is a {\it Hilbert--Schmidt operator} 
if Hilbert--Schmidt norm is finite; 
$\|\cA \|_{\mathrm{HS}}^2 :=\tr(\cA^{\ast} \cA) < \infty$.
We say that $\cA$ is a {\it locally Hilbert--Schmidt operator}
or {\it of locally Hilbert--Schmidt class}, 
if $\cA \cP_{\Lambda}$ is
a Hilbert--Schmidt operator for any compact set
$\Lambda \subset S$. 
It is known as the \textit{kernel theorem} that every
Hilbert--Schmidt operator $\cA : L^2(S_1, \la_1) \to
L^2(S_2, \la_2)$ is defined as an integral operator with kernel $k
\in L^2(S_1 \times S_2, \la_1 \otimes \la_2)$ (cf. Theorem
12.6.2 \cite{Aub00}). 

We put the second assumption.

\vskip 0.3cm
\noindent{\bf Assumption 2} \,
\begin{description}
\item{(i)} \quad $\cW$ is a locally Hilbert--Schmidt operator, 
\item{(ii)} \quad
$\cW^*$ is a locally Hilbert--Schmidt operator. 
\end{description}
\vskip 0.3cm

We note that for any compact set $\La_1 \subset S_1$, 
the operator $\cW \cP_{\La_1}$ is of Hilbert--Schmidt class 
if and only if 
the operator $\cP_{\La_1} \cW^* \cW \cP_{\La_1}$
is of trace class since 
\[
 \|\cW \cP_{\La_1}\|_{\mathrm{HS}}^2 := 
\tr\Big( (\cW \cP_{\La_1})^* \cW \cP_{\La_1} \Big)
=\tr\Big( \cP_{\La_1} \cW^* \cW \cP_{\La_1} \Big) < \infty. 
\]
Therefore, Assumption 2 (i) (resp. Assumption 2 (ii)) 
is equivalent to the following
Assumption 2' (i) (resp. Assumption 2' (ii)), 
which guarantees the existence of DPP
associated with $\cW^{\ast} \cW$ (resp. $\cW \cW^{\ast}$).  
\vskip 0.3cm
\noindent{\bf Assumption 2'} \,
\begin{description}
\item{(i)} \quad $\cW^{\ast} \cW \in \locally(S_1, \lambda_1)$, 
\item{(ii)} \quad $\cW \cW^{\ast} \in \locally(S_2, \lambda_2)$.
\end{description} 
\vskip 0.3cm

Given a measure space $(S, \lambda)$,  
if $f \in L^2(\Lambda, \lambda)$ for all compact
subsets $\Lambda$ of $S$, 
then $f$ is said to be locally $L^2$-integrable.
The set of all such functions is denoted by
$L^2_{\mathrm{loc}}(S, \la)$. 
By this definition if $\cP_{\Lambda} f \in L^2(S, \lambda)$
for any compact set $\Lambda \subset S$, then
$f \in L^2_{\mathrm{loc}}(S, \lambda)$. 
The following proposition is a \textit{local} version of
the kernel theorem for Hilbert--Schmidt operators. 

\begin{prop}
\label{prop:kernel-thm}
Suppose Assumption 2 (i) holds.
Then, $\cW$ is regarded as an integral operator associated with
a kernel $W : S_2 \times S_1 \to \C$; 
\begin{equation}
(\cW f)(y) = \int_{S_1} W(y, x) f(x) \lambda_1(dx),
\quad f \in L^2(S_1, \lambda_1), 
\label{eqn:W1}
\end{equation}
such that 
$\Psi_1 \in L^2_{\mathrm{loc}}(S_1, \la_1)$, where
$\Psi_1(x) := \|W(\cdot, x)\|_{L^2(S_2, \la_2)}, x \in S_1$.
\end{prop}

From Proposition~\ref{prop:kernel-thm}, 
under Assumption 2 (ii), the dual operator
$\cW^{\ast}$ also admits an integral kernel 
$W^* : S_1 \times S_2 \to \C$ such that 
$\Psi_2 \in L^2_{\mathrm{loc}}(S_2, \la_2)$, where
$\Psi_2(y) := \|W^*(\cdot, y)\|_{L^2(S_1, \la_1)}, y \in S_2$.
It is easy to see that 
$W^{\ast}(x,y) = \overline{W(y,x)}$ for $\la_1 \otimes\la_2$-a.e.$(x,y)$. 
Then
\begin{equation}
(\cW^{\ast} g)(x)=\int_{S_2} \overline{W(y, x)} g(y) \lambda_2(dy),
\quad g \in L^2(S_2, \lambda_2).
\label{eqn:W1*}
\end{equation}

Following (\ref{eqn:W1}) and (\ref{eqn:W1*}), 
we have
\begin{align*}
(\cW^{\ast} \cW f)(x) &= \int_{S_1} K_{S_1}(x, x') f(x') \lambda_1(dx'),
\quad f \in L^2(S_1, \lambda_1),
\nonumber\\
(\cW \cW^{\ast} g)(y) &= \int_{S_2} K_{S_2}(y, y') g(y') \lambda_2(dy'),
\quad g \in L^2(S_2, \lambda_2), 
\end{align*}
with the integral kernels, 
\begin{align}
K_{S_1}(x, x') &= \int_{S_2} \overline{W(y, x)} W(y, x') \lambda_2(dy)
=\langle W(\cdot, x'), W(\cdot, x) \rangle_{L^2(S_2, \lambda_2)},
\nonumber\\
K_{S_2}(y,y') &= \int_{S_1} W(y,x) \overline{W(y', x)} \lambda_1(dx)
=\langle W(y, \cdot), W(y', \cdot) \rangle_{L^2(S_1, \lambda_1)}.
\label{eqn:K1}
\end{align}
We see that $\overline{K_{S_1}(x', x)}=K_{S_1}(x, x')$
and $\overline{K_{S_2}(y', y)}=K_{S_2}(y, y')$.

Under Assumptions 1 and 2, we obtain the
following theorem 
as an immediate consequence of 
the well-known existence theorem
of DPP (Theorem \ref{thm:exist_DPP}). 
This is a starting-point for our discussion in the present section. 
\begin{thm}\label{thm:main1}
Under Assumptions 1 and 2, 
associated with $\cW^{\ast} \cW$ and $\cW \cW^{\ast}$,
there exists a unique pair of DPPs; 
$(\Xi_1, K_{S_1}, \lambda_1(d x))$
on $S_1$ and
$(\Xi_2, K_{S_2}, \lambda_2(d y))$
on $S_2$.
The correlation kernels 
$K_{S_{\ell}}, \ell=1,2$ are 
Hermitian and given by (\ref{eqn:K1}). 
\end{thm}
\vskip 0.3cm
Note that the densities of
the DPPs, $(\Xi_1, K_{S_1}, \lambda_1(d x))$ 
and $(\Xi_2, K_{S_2}, \lambda_2(d y))$,
are given by
\begin{align*}
\rho_1(x) &=K_{S_1}(x,x) 
=\int_{S_2} |W(y,x)|^2 \lambda_2(dy)
=\|W(\cdot, x)\|_{L^2(S_2,\lambda_2)}^2, 
\quad x \in S_1, 
\nonumber\\
\rho_2(y) &= K_{S_2}(y, y) 
=\int_{S_1} |W(y,x)|^2 \lambda_1(dx)
=\|W(y, \cdot)\|_{L^2(S_1,\lambda_1)}^2,
\quad y \in S_2,
\end{align*}
with respect to the background measures $\lambda_1(dx)$
and $\lambda_2(dy)$, respectively.

We say that 
a pair of DPPs $(\Xi_1, K_{S_1}, \lambda_1(d x))$
on $S_1$ and $(\Xi_2, K_{S_2}, \lambda_2(d y))$
on $S_2$ is associated with $\cW$.

\subsection{Basic properties of DPPs}
\label{sec:basic_DPP}

For $v=(v^{(1)}, \dots, v^{(d)}) \in \R^d$,
$y=(y^{(1)}, \dots, y^{(d)}) \in \R^d$, $d \in \N$,
the inner product of them is given by
$v \cdot y = y \cdot v
: =\sum_{a=1}^d v^{(a)} y^{(a)}$,
and $|v|^2 := v \cdot v$.
When $S \subset \C^d, d \in \N$,
$x \in S$ has $d$ complex components;
$x=(x^{(1)}, \dots, x^{(d)})$ with
$x^{(a)}=\Re x^{(a)}+\sqrt{-1} \Im x^{(a)}$, $a=1, \dots, d$.
In order to describe clearly such a complex structure, 
we set
$x_{\rR}=(\Re x^{(1)}, \dots, \Re x^{(d)}) \in \R^d$,
$x_{\rI}=(\Im x^{(1)}, \dots, \Im x^{(d)}) \in \R^d$,
and write $x=x_{\rR}+\sqrt{-1} x_{\rI}$ in this manuscript.
The Lebesgue measure is written as
$d x = d x_{\rR} d x_{\rI} := \prod_{a=1}^d d \Re x^{(a)} d \Im x^{(a)}$.
The complex conjugate of $x=x_{\rR}+\sqrt{-1} x_{\rI}$
is defined as $\overline{x}=x_{\rR}-\sqrt{-1} x_{\rI}$.
For $x = x_{\rR}+ \sqrt{-1} x_{\rI}$,
$x'=x'_{\rR}+ \sqrt{-1} x'_{\rI} \in \C^d$,
we use the {\it Hermitian inner product}; 
\[
x \cdot \overline{x'} := (x_{\rR}+\sqrt{-1} x_{\rI}) \cdot (x'_{\rR}-\sqrt{-1} x'_{\rI})
=(x_{\rR} \cdot x'_{\rR}+x_{\rI} \cdot x'_{\rI})
-\sqrt{-1} (x_{\rR} \cdot x'_{\rI}- x_{\rI} \cdot x'_{\rR})
\]
and define
\[
|x|^2 := x \cdot \overline{x} = |x_{\rR}|^2+|x_{\rI}|^2,
\quad x \in \C^d.
\]

For $(\Xi, K, \lambda(d x))$ defined 
on $S=\R^d$, $S=\C^d$, 
or on the space having appropriate periodicities, 
we introduce the 
following operations.
\begin{description}
\item{\bf (Shift)} 
For $u \in S$, 
$\cS_{u} \Xi :=\sum_j \delta_{x_j-u}$,
\[
\cS_{u} K(x, x')=K(x+u, x'+u), 
\]
and $\cS_{u} \lambda(d x)=\lambda(u+d x)$.
We write $(\cS_{u} \Xi, \cS_{u} K, \cS_{u} \lambda(d x))$ 
simply as $\cS_{u} (\Xi, K, \lambda(d x))$.
\item{\bf (Dilatation)} 
For $c>0$, we set
$c \circ \Xi :=\sum_j \delta_{c x_j}$ 
\[
c \circ K(x,x') :=
K \left( \frac{x}{c}, \frac{x'}{c} \right), 
\quad x, x' \in c S :=\{ c x : x \in S\},
\]
and $c \circ \lambda(d x):=\lambda(d x/c)$.
We define $c \circ (\Xi, K, \lambda(d x)) 
:= (c \circ \Xi, c \circ K, c \circ \lambda(d x))$.
\item{\bf (Square root)}
For $(\Xi, K, \lambda(dx))$ on $S=[0, \infty)$, we put
$\Xi^{\langle 1/2 \rangle} :=\sum_j \delta_{\sqrt{x_j}}$, 
$K^{\langle 1/2 \rangle}(x, x') :=K(x^2, {x'}^2)$,
and $\lambda^{\langle 1/2 \rangle}(d x) :=(\lambda \circ v^{-1}) (d x)$,
where $v(x)=\sqrt{x}$.
We define $(\Xi, K, \lambda(d x))^{\langle 1/2 \rangle}$
$:= (\Xi^{\langle 1/2 \rangle}, K^{\langle 1/2 \rangle}, 
\lambda^{\langle 1/2 \rangle}(d x))$
on $[0, \infty)$.
\item{\bf (Gauge transformation)} 
For non-vanishing $u: S \to \C$, 
a gauge transformation of $K$ by $u$ is defined as
\[
K(x, x') \mapsto \widetilde{K}_u(x, x') :=u(x) K(x, x') u(x')^{-1}.
\]
In particular, when $u: S \to \U(1)$,
the $\U(1)$-gauge transformation of $K$ is given by
\[
K(x, x') \mapsto \widetilde{K}_u(x, x') :=u(x) K(x, x') \overline{u(x')}.
\]
\end{description}

We will use the following basic properties of DPP.
\begin{description}
\item{\bf [Gauge invariance]} 
For any $u: S \to \C$, a gauge transformation
does not change the probability law of DPP;
\[
(\Xi, K, \lambda(d x)) \law= (\Xi, \widetilde{K}_u, \lambda(d x)).
\]

\item{\bf [Measure change]}
For a measurable function $g : S \to \R_{\geq 0}$,
\begin{equation}
(\Xi, K(x, x^{\prime}), g(x) \lambda(dx))
\law=
(\Xi, \sqrt{g(x)} K(x, x^{\prime}) \sqrt{g(x^{\prime})}, \lambda(dx)).
\label{eqn:measure_change}
\end{equation}

\item{\bf [Mapping and scaling]} 
For a one-to-one measurable mapping
$h: S \to \widehat{S}$, if we set
\[
\widehat{\Xi}=\sum_j \delta_{h(x_j)},
\quad
\widehat{K}(x, x')=K(h^{-1}(x), h^{-1}(y)),
\quad 
\widehat{\lambda}(dx)= (\lambda \circ h^{-1})(dx),
\]
then 
$(\widehat{\Xi}, \widehat{K}, \widehat{\lambda}(dx))$
is a DPP on $\widehat{S}$.
In particular, 
when $h(x)=x-u, u \in S$,
$(\widehat{\Xi}, \widehat{K}, \widehat{\lambda}(dx))
=\cS_u (\Xi, K, \lambda(dx))$,
when $h(x)=cx, c > 0$,
$(\widehat{\Xi}, \widehat{K}, \widehat{\lambda}(dx))
=c \circ (\Xi, K, \lambda(dx))$,
and when $h(x)=\sqrt{x}$ for $S=[0, \infty)$, 
$(\widehat{\Xi}, \widehat{K}, \widehat{\lambda}(dx))
=(\Xi, K, \lambda(dx))^{\langle 1/2 \rangle}$.
If $c \circ \lambda(d x) = c^{-d} \lambda(d x)$, 
then (\ref{eqn:measure_change}) with $g(x) \equiv c >0$ gives 
\[
c \circ (\Xi, K, \lambda(d x)) \law= (c \circ \Xi, K_c, \lambda(d x)),
\quad c >0, 
\]
with
\[
K_c(x, x')=\frac{1}{c^d} K \left( \frac{x}{c}, \frac{x'}{c} \right),
\]
where the base space is given by $cS$.
\end{description}

We will give some limit theorems for DPPs in this manuscript.
Consider a DPP which depends on a continuous
parameter, or a series of DPPs labeled by
a discrete parameter (e.g., the number of points $N \in \N$),
and describe the system by 
$(\Xi, K_p, \lambda_p(dx))$ with 
the continuous or discrete parameter $p$.
If $(\Xi, K_p, \lambda_p(dx))$ converges
to a DPP, $(\Xi, K, \lambda(dx))$, as $p \to \infty$, 
weakly in the vague topology,
we write this limit theorem as
$(\Xi, K_p, \lambda_p(dx)) \weakp (\Xi, K, \lambda(d x))$.
The weak convergence of DPPs is verified by 
the uniform convergence of the kernel $K_p \to K$
on each compact set $C \subset S \times S$ \cite{ST03a}.

\subsection{Duality relations}
\label{sec:duality}
For $f \in \cC_{\rm c}(S)$, the Laplace transform
of the probability measure $\bP$ for a point 
process $\Xi$ is defined as
\begin{equation}
\Psi[f]=\bE \left[
\exp \left( \int_{S} f(x) \Xi(d x) \right) \right].
\label{eqn:Laplace1}
\end{equation}
For the DPP, $(\Xi, K, \lambda(d x))$,
this is given by the Fredholm determinant
on $L^2(S, \lambda)$ \cite{Sim05},
\begin{align*}
\Det_{L^2(S, \lambda)} [ I - (1-e^f) \cK]
& := 1+\sum_{n \in \N}
\frac{(-1)^n}{n!} \int_{S^n} \det_{1 \leq j, k \leq n} [K(x_j, x_k)]
\prod_{\ell=1}^n (1-e^{f(x_{\ell})}) \lambda^{\otimes n}(d \x).
\end{align*}

\begin{lem}
\label{thm:duality}
Between two DPPs, 
$(\Xi_1, K_{S_1}, \lambda_1(d x))$
on $S_1$ and
$(\Xi_2, K_{S_2}, \lambda_2(d y))$
on $S_2$, given by Theorem \ref{thm:main1},
the following equality holds with an arbitrary parameter $\alpha \in \C$,
\begin{equation}
\Det_{L^2(S_1, \lambda_1)} [I+\alpha \cK_{S_1} ]
=\Det_{L^2(S_2, \lambda_2)} [I+\alpha \cK_{S_2} ].
\label{eqn:duality}
\end{equation}
\end{lem}
\noindent{\it Proof} \quad 
We recall that if $\cA \cB$ and $\cB \cA$ are trace class operators on
a Hilbert space $H$ then \cite{Sim05}
\begin{equation}
 \Det_H [I+ \cB \cA] = \Det_H [I + \cA \cB]. 
\label{eqn:commutative} 
\end{equation}
Now we have $\cA : H_1 \to H_2$ and $\cB : H_2 \to H_1$ between  
two Hilbert spaces $H_1$ and $H_2$. 
Let $\widetilde{\cA}$ and $\widetilde{\cB}$ be two operators on 
$H_1\oplus H_2$ defined by 
\[
 \widetilde{\cA}
= 
\begin{pmatrix}
O & O \\
\cA & O \\
\end{pmatrix}, \quad 
\widetilde{\cB} 
= \begin{pmatrix}
O & \cB \\
O & O \\
\end{pmatrix}
\]
Then, $\widetilde{\cA} \widetilde{\cB}$ and $\widetilde{\cB} \widetilde{\cA}$ are diagonal operators 
$O \oplus \cA \cB$ and $\cB \cA \oplus O$, respectively, and hence 
also they are trace class operators. 
By applying (\ref{eqn:commutative}) to $\widetilde{\cA}$ and
$\widetilde{\cB}$ with $H := H_1 \oplus H_2$, we obtain 
\[
 \Det_{H_1} [I+\cB \cA] = \Det_{H_2} [I+ \cA \cB]. 
\]
Consequently, taking $\cA = \sqrt{\alpha} \cW$, $\cB=\sqrt{\alpha} \cW^*$, 
$H_1=L^2(S_1, \lambda_1)$, and $H_2=L^2(S_2, \lambda_2)$ yields 
(\ref{eqn:duality}). 
\qed
\vskip 0.3cm

For $\Lambda_{\ell} \subset S_{\ell}, \ell=1,2$, let
\begin{equation}
\widetilde{\cW} := \cP_{\Lambda_2} \cW \cP_{\Lambda_1},
\quad
\cK_{S_1}^{(\Lambda_2)}
:= \cW^{\ast} \cP_{\Lambda_2} \cW,
\quad
\cK_{S_2}^{(\Lambda_1)}
:= \cW \cP_{\Lambda_1} \cW^{\ast}.
\label{eqn:defK_A}
\end{equation}
They admit the following integral kernels, 
\begin{align}
\widetilde{W}(y, x)
&={\bf 1}_{\Lambda_2}(y) W(y, x) {\bf 1}_{\Lambda_1}(x),
\nonumber\\
K_{S_1}^{(\Lambda_2)}(x, x^{\prime})
&= \int_{\Lambda_2} \overline{W(y, x)} W(y, x^{\prime}) \lambda_2(d y),
\nonumber\\
K_{S_2}^{(\Lambda_1)}(y, y^{\prime})
&= \int_{\Lambda_1} W(y, x) \overline{W(y^{\prime}, x)} \lambda_1(d x).
\label{eqn:K_A}
\end{align}
Using Lemma \ref{thm:duality}, the following theorem is proved.
\begin{thm}
\label{thm:duality2}
Let $(\Xi_1^{(\Lambda_2)}, K_{S_1}^{(\Lambda_2)}, \lambda_1(d x))$
and $(\Xi_2^{(\Lambda_1)}, K_{S_2}^{(\Lambda_1)}, \lambda_2(d y))$ 
be DPPs associated with the kernels
$K_{S_1}^{(\Lambda_2)}$ and $K_{S_2}^{(\Lambda_1)}$
given by (\ref{eqn:K_A}), respectively.
Then, $\Xi_1^{(\Lambda_2)}(\Lambda_1) \law= 
 \Xi_2^{(\Lambda_1)}(\Lambda_2)$, i.e., 
\[
\bP( \Xi_1^{(\Lambda_2)}(\Lambda_1)=m)
= \bP( \Xi_2^{(\Lambda_1)}(\Lambda_2)=m),
\quad \forall m \in \N_0.
\]
\end{thm}
\noindent{\it Proof} \,
As a special case of (\ref{eqn:Laplace1})
with $f(x)={\bf 1}_{\Lambda_1}(x) \log z$
for $\Xi=\Xi_1^{(\Lambda_2)}$,
$z \in \C$, we have the equality,
\begin{equation}
\bE\left[ z^{\Xi_1^{(\Lambda_2)}(\Lambda_1)} \right]
=\Det_{L^2(S_1, \lambda_1)}
[I-(1-z) \cP_{\Lambda_1} \cK_{S_1}^{(\Lambda_2)} \cP_{\Lambda_1} ],
\label{eqn:generating1}
\end{equation}
where $\cK_{S_1}^{(\Lambda_2)}$ is defined by (\ref{eqn:defK_A}). 
Here LHS is the moment generating function
of $\Xi^{(\Lambda_2)}_1(\Lambda_1)$ 
and RHS gives its Fredholm determinantal expression.
By replacing $\cW$ by $\widetilde{\cW}$ and letting $\alpha=-(1-z)$ in the proof of
Lemma \ref{thm:duality}, we obtain the equality,
\[
\Det_{L^2(S_1, \lambda_1)}
[I-(1-z) \cP_{\Lambda_1} \cK_{S_1}^{(\Lambda_2)} \cP_{\Lambda_1} ]
=
\Det_{L^2(S_2, \lambda_2)}
[I-(1-z) \cP_{\Lambda_2} \cK_{S_2}^{(\Lambda_1)} \cP_{\Lambda_2} ].
\]
Through (\ref{eqn:generating1})
and the similar equality for 
$\bE\left[ z^{\Xi_2^{(\Lambda_1)}(\Lambda_2)} \right]$, 
we obtain the corresponding
equivalence between the moment generating functions
of $\Xi_1^{(\Lambda_2)}(\Lambda_1)$ and $\Xi_2^{(\Lambda_1)}(\Lambda_2)$,
and hence the statement of the proposition is proved. \qed


\subsection{Orthonormal functions and correlation kernels}
\label{sec:orthogonal}

In addition to $L^2(S_{\ell}, \lambda_{\ell})$, $\ell=1, 2$, 
we introduce $L^2(\Gamma, \nu)$ as a parameter space
for functions in $L^2(S_{\ell}, \lambda_{\ell}), \ell=1,2$. 
Assume that there are two families of measurable
functions $\{\psi_1(x, \gamma) : x \in S_1, 
\gamma \in \Gamma \}$ and
$\{\psi_2(y, \gamma) : y \in S_2, 
\gamma \in \Gamma \}$ 
such that 
two bounded operators $\cU_{\ell} : L^2(S_{\ell}, \lambda_{\ell}) \to L^2(\Gamma, \nu)$
given by 
\[
(\cU_{\ell}f)(\gamma) := 
\int_{S_{\ell}} \overline{\psi_{\ell}(x, \gamma)} f(x)
\lambda_{\ell}(dx), \quad \ell=1,2, 
\]
are well-defined. 
Then, their adjoints $\cU_{\ell}^* : L^2(\Gamma, \nu) \to 
L^2(S_{\ell}, \lambda_{\ell}), \ell =1,2$ 
are given by 
\[
(\cU_{\ell}^* F)(\cdot) = 
\int_{\Gamma} \psi_{\ell}(\cdot, \gamma) F(\gamma)
\nu(d\gamma). 
\]
A typical example of $\cU_1$ is the Fourier transform,
i.e., $\psi_1(x,\gamma) = e^{\sqrt{-1} x \gamma}$. In this case, for
any $\gamma$, the function $\psi_1(\cdot,\gamma)$ is \textit{not} in $L^2(\R, dx)$.
Now we define $\cW : L^2(S_1, \lambda_1) \to L^2(S_2, \lambda_2)$ by 
$\cW = \cU_2^* \cU_1$, i.e., 
\begin{equation}
 (\cW f)(y) = \int_{\Gamma} \psi_2(y,\gamma) (\cU_{1} f)(\gamma) 
 \nu(d\gamma). 
\label{eqn:W_of}
\end{equation}
Let $I_{\Gamma}$ be an identity in $L^2(\Gamma, \nu)$.
We can see the following. 
\begin{lem} 
If $\cU_{\ell} \cU_{\ell}^* = I_{\Gamma}$ for $\ell=1,2$, then 
both $\cW$ and $\cW^*$ are partial isometries. 
\end{lem}
\noindent{\it Proof} \,
By the assumption, we see that 
\[
\cW \cW^* \cW = 
(\cU_2^* \cU_1) (\cU_1^* \cU_2) (\cU_2^* \cU_1) 
= \cU_2^* \cU_1 = \cW. 
\]
From Lemma~\ref{lem:A_Astar_partial_isometry2}, $\cW$ is a partial isometry. 
By symmetry, the assertion for $\cW^*$ also follows. 
\qed
\vskip 0.3cm
We note that 
$\cW^* \cW = \cU_1^* \cU_1$ and 
$\cW \cW^* = \cU_2^* \cU_2$. Hence,  
$\cU_{\ell}, \ell=1,2$ are also partial isometries.
In addition, 
$\cW^{\ast} \cW$ is a locally trace
class operator if and only if so is $\cU_1^{\ast} \cU_1$. 
Therefore, $\cW$ is locally Hilbert--Schmidt if and only if 
so is $\cU_1$. 

Now we rewrite the condition for 
$\cU_1$ to be of locally Hilbert--Schmidt class in terms of the function
$\psi_1(x, \gamma), x \in S_1, \gamma \in \Gamma$. 

\begin{lem}\label{lem:HScondition}
Let $\Psi_1(x) := \|\psi_1(x, \cdot)\|_{L^2(\Gamma, \nu)}$, $x \in S_1$
and assume that 
$\Psi_1 \in L^2_{\mathrm{loc}}(S_1, \la_1)$. 
Then, the operator $\cU_1$ is of locally Hilbert--Schmidt class. 
\end{lem}
\noindent{\it Proof} \quad
For a compact set $\La \subset S_1$, we see that 
\begin{align*}
|\cP_{\La} \cU_1^* \cU_1 \cP_{\La} f (x)|
&= \left|\mathbf{1}_{\La}(x) 
\int_{\Gamma} \nu(d\gamma) \psi_1(x,\gamma) 
\int_{S_1} \overline{\psi_1(x',\gamma)}
 \mathbf{1}_{\La}(x') f(x')\la_1(dx')\right|  \\
&\le \mathbf{1}_{\La}(x) \Psi_1(x) \int_{S_1} \mathbf{1}_{\La}(x') \Psi_1(x')
|f(x')| \la_1(dx') \\
&\le \cP_{\La} \Psi_1(x) \|\cP_{\La} \Psi_1\|_{L^2(S_1,\la_1)}
\|\cP_{\La} f\|_{L^2(S_1,\la_1)}. 
\end{align*}
By Fubini's theorem, we have 
\[
\cP_{\La} \cU_1^* \cU_1 \cP_{\La} f (x)
= 
\int_{S_1} \la_1(dx') f(x')
\left(\int_{\Gamma} \mathbf{1}_{\La}(x) \psi_1(x,\gamma) 
\overline{\mathbf{1}_{\La}(x') \psi_1(x',\gamma)}
\nu(d\gamma)\right)
\]
and hence
\[
\|\cU_1 \cP_{\La}\|^2_{\mathrm{HS}} = \tr (\cP_{\La} \cU_1^* \cU_1
\cP_{\La}) = 
\int_{S_1} \la_1(dx) \mathbf{1}_{\La}(x) \left(\int_{\Gamma}
|\psi_1(x,\gamma)|^2 \nu(d\gamma)\right) 
= \|\cP_{\La} \Psi_1\|^2_{L^2(S_1, \la_1)} < \infty. 
\]
This completes the proof. 
\qed
\vskip 0.3cm

Now we put the following.

\vskip 0.3cm
\noindent{\bf Assumption 3} \,
For $\ell=1,2$, 
\begin{description}
\item{(i)} \quad $\cU_{\ell} \cU_{\ell}^* = I_{\Gamma}$,
\item{(ii)} \quad $\Psi_{\ell} \in L^2_{\mathrm{loc}}(S_{\ell}, \la_{\ell})$,
where 
$\Psi_{\ell}(x):=\|\psi_{\ell}(x,\cdot)\|_{L^2(\Gamma, \nu)}$,
$x \in S_{\ell}$. 
\end{description}
\vskip 0.3cm
\noindent
Assumption 3(i) can be rephrased as the following 
orthonormality relations: 
\[
\langle 
\psi_{\ell}(\cdot, \gamma), \psi_{\ell}(\cdot, \gamma') \rangle_{L^2(S_{\ell}, \lambda_{\ell})}
\nu(d \gamma)
=\delta(\gamma-\gamma') d \gamma, 
\quad \gamma, \gamma' \in \Gamma,
\quad \ell=1,2.
\]
We often use these relations below. 
\vskip 0.3cm
The following is immediately obtained as a corollary of
Theorem \ref{thm:main1}.
\begin{cor}
\label{thm:main2}
Let $\cW = \cU_2^* \cU_1$ as in the above. 
We assume Assumption 3. Then, there exist unique pair of DPPs;
$(\Xi_1, K_{S_1}, \lambda_1(d x))$
on $S_1$ and
$(\Xi_2, K_{S_2}, \lambda_2(d y))$
on $S_2$.
Here the correlation kernels 
$K_{S_{\ell}}, \ell=1,2$ are given by 
\begin{align}
K_{S_1}(x, x')
&= \int_{\Gamma} \psi_1(x, \gamma) 
\overline{\psi_1(x', \gamma)} \nu(d \gamma)
=\langle \psi_1(x, \cdot), \psi_1(x', \cdot) \rangle_{L^2(\Gamma, \nu)},
\nonumber\\
K_{S_2}(y, y')
&= \int_{\Gamma} \psi_2(y, \gamma) 
\overline{\psi_2(y', \gamma)} \nu(d \gamma)
=\langle \psi_2(y, \cdot), \psi_2(y', \cdot) \rangle_{L^2(\Gamma, \nu)}.
\label{eqn:K_main2}
\end{align}
In particular, the densities of the DPPs are given by
$\rho_1(x)=K_{S_1}(x,x) = \Psi_1(x)^2, x \in S_1$ 
and $\rho_2(y)=K_{S_2}(y,y)=\Psi_2(y)^2, y \in S_2$ 
with respect to the background measures
$\lambda_1(dx)$ and $\lambda_2(dy)$, respectively.
\end{cor}
\vskip 0.3cm
\begin{rem}
\label{thm:Remark1_1}
Consider the symmetric case such that 
$L^2(S_1, \lambda_1)=L^2(S_2, \lambda_2)=: L^2(S, \lambda)$, 
$\psi_1=\psi_2 =: \psi$, 
$\nu=\lambda|_{\Gamma}$, $\Gamma \subseteq S$. 
In this case, $\cW = \cU^* \cU$ with 
\[
 (\cU f)(\gamma) = \int_S  \overline{\psi(x,\gamma)} f(x)
 \lambda(dx). 
\]
Then 
$K_{S_1}=K_{S_2}=W = : K$
is given by
\begin{equation}
K(x, x')=\int_{\Gamma} \psi(x, \gamma) \overline{\psi(x', \gamma)} \lambda(d \gamma).
\label{eqn:K_simple1}
\end{equation}
This is Hermitian; 
$\overline{K(x', x)}=K(x, x')$, and satisfies
the reproducing property 
\[
K(x, x')=\int_{S} K(x, \zeta) K(\zeta, x') \lambda(d \zeta).
\]
\end{rem}
\vskip 0.3cm

Now we consider a simplified version of 
the preceding setting. 
Let $\Gamma \subseteq S_2$ and $\nu=\lambda_2|_{\Gamma}$.
We define 
$\cU_2 : L^2(S_2,\lambda_2) \to L^2(\Gamma, \nu)$ as the
restriction onto $\Gamma$, and then its adjoint $\cU_2^*$ is given by 
$(\cU_2^* F)(y) = F(y)$ for $y \in \Gamma$, and by $0$ for $y \in
S_2 \setminus \Gamma$. 
We write the extension $\tilde{F} = \cU^*_2 F$ for $F \in L^2(\Gamma,
\nu)$. 
It is obvious that $\cU_2 \cU_2^* =
I_{\Gamma}$ and hence $\cU_2$ is a partial isometry. 

For $\Gamma \subseteq S_2$, we assume that there is a family
of measurable functions 
$\{\psi_{1}(x, y) : x \in S_1, y \in \Gamma \}$ 
such that a bounded operator 
$\cU_1 : L^2(S_1,\lambda_1) \to L^2(\Gamma, \nu)$ given by 
\[
 (\cU_1 f)(\gamma) := \int_{S_1} \overline{\psi_1(x,\gamma)}
 f(x) \lambda_1(dx) \quad (\gamma \in \Gamma)
\]
is well-defined. 

\vskip 0.3cm
\noindent{\bf Assumption 3'} \,
\begin{description}
\item{(i)} \quad $\cU_1 \cU_1^* = I_{\Gamma}$,
\item{(ii)} \quad $\Psi_1 \in L^2_{\mathrm{loc}}(S_1, \la_1)$,
where 
$\Psi_1(x):=\|\psi_1(x,\cdot)\|_{L^2(\Gamma, \nu)}$,
$x \in S_1$. 
\end{description}
\vskip 0.3cm
\noindent
Assumption 3'(i) can be rephrased as the following 
orthonormality relation: 
\[
 \langle 
 \psi_{1}(\cdot, y), \psi_{1}(\cdot, y') \rangle_{L^2(S_{1}, \lambda_{1})}
 \lambda_2(d y)
 =\delta(y-y') d y, 
 \quad y, y' \in \Gamma.
 \]
Now we define $\cW : L^2(S_1, \lambda_1) \to L^2(S_2, \lambda_2)$ by 
$\cW = \cU_2^* \cU_1$ as before. In this case, we have  
\[
 (\cW f)(y) 
= \mathbf{1}_{\Gamma}(y) \int_{S_1}
\overline{\tilde{\psi}_1(x,y)}f(x) \lambda_1(dx),  
\]
and hence 
\begin{equation}
W(y, x)= \overline{\tilde{\psi}_1(x, y)} {\bf 1}_{\Gamma}(y).
\label{eqn:w_set2}
\end{equation}
It follows from Assumption 3' that $\cW$ is a partial isometry.  
Corollary \ref{thm:main2} is reduced to the following.
\begin{cor}
\label{thm:main3}
Let $\cW = \cU_2^* \cU_1$ as in the above. We assume Assumption 3'. 
Then there exists a unique DPP, 
$(\Xi, K, \lambda_{1})$
on $S_{1}$ with
the correlation kernel
\begin{equation}
K_{S_1}(x, x^{\prime})
= \int_{\Gamma} \psi_1(x, y) \overline{\psi_1(x^{\prime}, y)} \lambda_2(d y)
=\langle \tilde{\psi}_1(x, \cdot), \tilde{\psi}_1(x^{\prime}, \cdot) \rangle_{L^2(\Gamma, \lambda_2)}.
\label{eqn:K_simple2}
\end{equation}
In particular, the density of the DPP is given by
$\rho_1(x)=K_{S_1}(x,x) = \Psi_1(x)^2, x \in S_1$ 
with respect to the background measures
$\lambda_1(dx)$.
\end{cor}

\subsection
{Examples in one-dimensional spaces} \label{sec:examples_1d}
\subsubsection{Finite DPPs in $\R$ associated with 
classical orthonormal polynomials}
\label{sec:classical}
Let $S_1=S_2 = \R$. 
Assume that we have
two sets of orthonormal functions 
$\{\varphi_n \}_{n \in \N_0}$ 
and $\{\phi_n \}_{n \in \N_0}$
with respect to the measures $\lambda_1$ and $\lambda_2$,
respectively,
\begin{align}
\langle \varphi_n, \varphi_m \rangle_{L^2(\R, \lambda_1)}
&=
\int_{\R} \varphi_n(x) \overline{\varphi_m(x)} \lambda_1(d x) = \delta_{nm},
\nonumber\\
\langle \phi_n, \phi_m \rangle_{L^2(\R, \lambda_2)}
&=
\int_{\R} \phi_n(y) \overline{\phi_m(y)} \lambda_2(d y) = \delta_{nm},
\quad n, m \in \N_0.
\label{eqn:ON1}
\end{align}
Then for an arbitrary but fixed $N \in \N$, 
we set $\Gamma=\{0,1, \dots, N-1\} \subsetneq \N_0$,
$\psi_1(\cdot, \gamma)=\varphi_{\gamma}(\cdot)$,
$\psi_2(\cdot, \gamma)=\phi_{\gamma}(\cdot)$, $\gamma \in \Gamma$, 
and consider $\ell^2(\Gamma)$ as $L^2(\Gamma, \nu)$
in the setting of Section \ref{sec:orthogonal}. 

\vskip 0.3cm
\begin{rem}
\label{thm:Remark1_2}
If $\Gamma$ is a finite set, $|\Gamma|=:N \in \N$,
and the parameter space is given by $\ell^2(\Gamma)$,
Assumption 3(ii) (resp. Assumption 3'(ii))
is concluded from 3(i) (resp. 3'(i)) as shown below.
Since
\[
\Psi(x)^2 :=\|\varphi_{\cdot}(x) \|^2_{\ell^2(\Gamma)}
=\sum_{n \in \Gamma} |\varphi_n(x)|^2, \quad x \in S, 
\]
we have
\[
\int_{S} \Psi(x)^2 \lambda(dx)
= \sum_{n \in \Gamma} \|\varphi_n \|^2_{L^2(S, \lambda)}.
\]
Then, if $\{\varphi_n\}_{n \in \Lambda}$ are normalized,
the above integral is equal to $|\Gamma|=N < \infty$.
This implies $\Psi \in L^2(S, \lambda)
\subset L^2_{\mathrm{loc}}(S, \lambda)$.
\end{rem}
\vskip 0.3cm

Hence Assumption 3 is satisfied for any $N \in \N$.
Then the integral kernel for $\cW$ defined by (\ref{eqn:W_of}) 
is given by
\[
W(y, x)= \sum_{n=0}^{N-1} \overline{\varphi_n(x)} \phi_n(y).
\]
By Corollary \ref{thm:main2}, we have a pair of DPPs on $\R$,
$(\Xi_1, K_{\varphi}^{(N)}, \lambda_1(d x))$
and $(\Xi_2, K_{\phi}^{(N)}, \lambda_2(d y))$, 
where the correlation kernels are given by
\begin{equation}
K_{\varphi}^{(N)}(x, x')=\sum_{n=0}^{N-1} \varphi_n(x) \overline{\varphi_n(x')},
\quad
K_{\phi}^{(N)}(y, y')=\sum_{n=0}^{N-1} \phi_n(y) \overline{\phi_n(y')},
\label{eqn:K_classical}
\end{equation}
respectively.
Here $N$ gives the number of points for each DPPs.
If we can use the three-term relations 
in $\{\varphi_n\}_{n \in \N_0}$ or $\{\phi_n\}_{n \in \N_0}$, 
(\ref{eqn:K_classical}) can be written 
in the Christoffel--Darboux form
(see, for instance, Proposition 5.1.3 in \cite{For10}). 
As a matter of course, if we have three or more than three,
say $M$ distinct sets of orthonormal functions
satisfying Assumption 3 with a common $\Gamma$, then
by applying Corollary \ref{thm:main2} to every pair of them,
we will obtain $M$ distinct finite DPPs.

Even if we have only one set of orthonormal functions,
for example, only the first one $\{\varphi_n\}_{n \in \N_0}$ in (\ref{eqn:ON1}),
we can obtain a DPP (labeled by the number of particles $N \in \N$) 
following Corollary \ref{thm:main3}.
In such a case, we set
\begin{equation}
W(n, x)=\overline{\varphi_n(x)} {\bf 1}_{\Gamma}(n)
\label{eqn:W_n_x1}
\end{equation}
with $\Gamma=\{0,1, \dots, N-1\}$
for (\ref{eqn:w_set2}). Then we have a DPP,
$(\Xi, K_{\varphi}^{(N)}, \lambda_1(d x))$.

Now we give classical examples of DPPs
associated with real orthonormal polynomials.
Let $\lambda_{\rN(m, \sigma^2)}(dx)$ denote the
{\it normal distribution}, 
\[
\lambda_{\rN(m, \sigma^2)}(dx)=\frac{1}{\sqrt{2 \pi} \sigma} 
e^{-(x-m)^2/(2\sigma^2)} d x,
\quad m \in \R, \quad \sigma >0, 
\]
and 
$\lambda_{\Gamma(a, b)}(d y)$ do the
{\it Gamma distribution}, 
\[
\lambda_{\Gamma(a, b)}(d y)
=\frac{b^a}{\Gamma(a)} y^{a-1} e^{-b y} {\bf 1}_{\R_{\geq 0}}(y) dy,
\quad a >0, \quad b >0, 
\]
with the Gamma function
$\Gamma(z) := \int_0^{\infty} u^{z-1} e^{-u} du, \Re z > 0$.
We set 
\begin{align}
\lambda_1(dx) &= \lambda_{\rN(0,1/2)}(d x)
= \frac{1}{\sqrt{\pi}} e^{-x^2} dx, 
\nonumber\\
\varphi_n(x) &= \frac{1}{\sqrt{2^n n!}} H_n(x), \quad n \in \N_0, 
\label{eqn:HermiteZ1}
\end{align}
and
\begin{align}
\lambda_2(d y) &=\lambda_{\Gamma(\nu+1, 1)}(d y)
= \frac{1}{\Gamma(\nu+1)} y^{\nu} e^{-y} {\bf 1}_{\R_{\geq 0}}(y)dy, 
\nonumber\\
\phi_n(y) &= \phi_n^{(\nu)}(y)=
\sqrt{ \frac{n! \Gamma(\nu+1)}{\Gamma(n+\nu+1)}} \, L^{(\nu)}_n(y),
\quad n \in \N_0, 
\label{eqn:LaguerreZ1}
\end{align}
with parameter $\nu \in (-1, \infty)$.
Here 
$\{H_n(x)\}_{n \in \N_0}$ are the {\it Hermite polynomials}, 
\begin{align}
H_n(x) & : =(-1)^n e^{x^2} \frac{d^n}{dx^n} e^{-x^2}
\nonumber\\
&= n! \sum_{k=0}^{[n/2]} 
\frac{(-1)^k (2x)^{n-2k}}{k! (n-2k)!}, \quad n \in \N_0,
\label{eqn:Hermite1}
\end{align}
where $[a]$ denotes the largest integer not greater than 
$a \in \R$, and
$\{L^{(\nu)}_n(x)\}_{n \in \N_0}$ are
the {\it Laguerre polynomials}, 
\begin{align}
L^{(\nu)}_n(x) & : = \frac{1}{n!} x^{-\nu} e^{x} \frac{d^n}{dx^n} \Big(x^{n+\nu} e^{-x} \Big)
\nonumber\\
&= \sum_{k=0}^{n}
\frac{(\nu+k+1)_{n-k}}{(n-k)! k!} (-x)^k, 
\quad n \in \N_0, \quad \nu \in (-1, \infty), 
\label{eqn:Laguerre1}
\end{align}
where $(\alpha)_n := \alpha(\alpha+1) \cdots (\alpha+n-1)
= \Gamma(\alpha+n)/\Gamma(\alpha)$, $n \in \N$, 
$(\alpha)_0 := 1$. 
The correlation kernels (\ref{eqn:K_classical}) are written 
in the {\it Christoffel--Darboux form} as,
\begin{align}
K^{(N)}_{\varphi}(x, x^{\prime}) &=K_{\rm Hermite}^{(N)}(x, x')
= \sum_{n=0}^{N-1} \varphi_n(x) \varphi_n(x^{\prime})
\nonumber\\
&=\sqrt{\frac{N}{2}}
\frac{\varphi_N(x) \varphi_{N-1}(x')-\varphi_N(x') \varphi_{N-1}(x)}
{x-x^{\prime}},
\quad x, x^{\prime} \in \R,
\label{eqn:K_Hermite}
\end{align}
and
\begin{align}
K^{(N)}_{\phi}(y, y') 
&=K_{\rm Laguerre}^{(\nu, N)}(y, y')
=\sum_{n=0}^{N-1} \phi^{(\nu)}_n(y) \phi^{(\nu)}_n(y^{\prime})
\nonumber\\
&=- \sqrt{N(N+\nu)}
\frac{\phi^{(\nu)}_N(y) \phi^{(\nu)}_{N-1}(y')
-\phi^{(\nu)}_N(y') \phi^{(\nu)}_{N-1}(y)}
{y-y'},
\quad y, y' \R_{\geq 0}.
\label{eqn:K_Laguerre}
\end{align}
When $x=x'$ or $y=y'$, we make sense of the above formulas
by using L'H\^opital's rule.
The former is called the {\it Hermite kernel}
and the latter is the {\it Laguerre kernel}. 

By definition, 
for a finite DPP $(\Xi, K, \lambda(dx))$ with $N$ points in $S$, 
the probability density with respect to
$\lambda^{\otimes N}(dx_1 \cdots dx_N)$ is given by
$\rho^{N}(x_1, \dots, x_N)
=\det_{1 \leq j, k \leq N}[K(x_j, x_k)], 
\x=(x_1, \dots, x_N) \in S^N$.
Using the {\it Vndermonde determinantal formula}, 
$\det_{1 \leq j, k \leq N}(z_k^{j-1})
=\prod_{1 \leq j < k \leq N}(z_k-z_j)$,
which will be given also as the type $\AN$ of
Weyl denominator formula (\ref{eqn:Weyl_denominator}) below, 
we can verify that 
the probability densities
of the DPPs
$(\Xi, K^{(N)}_{\rm Hermite}, \lambda_{\rN(0,1/2)}(dx))$
and 
$(\Xi, K^{(N)}_{\rm Laguerre}, \lambda_{\Gamma(\nu+1, 1)}(dy))$
with respect to 
the Lebesgue measures
$d \x =\prod_{j=1}^N d x_j$ and 
$d \y =\prod_{j=1}^N d y_j$ are given
as
\begin{align}
\bp^{(N)}_{\rm Hermite}(\x) 
&= \frac{1}{Z^{(N)}_{\rm Hermite}}
\prod_{1 \leq j < k \leq N} (x_k-x_j)^2
\prod_{\ell=1}^N e^{-x_{\ell}^2},
\quad \x=(x_1, \dots, x_N) \in \R^N,
\nonumber\\
\bp^{(\nu, N)}_{\rm Laguerre}(\y) 
&= \frac{1}{Z^{(\nu, N)}_{\rm Laguerre}}
\prod_{1 \leq j < k \leq N} (y_k-y_j)^2
\prod_{\ell=1}^N y_{\ell}^{\nu} e^{-y_{\ell}},
\quad \nu \in (-1, \infty), \quad \y \in \R_{\geq 0}^N,
\label{eqn:P_Laguerre}
\end{align}
with the normalization constants 
$Z^{(N)}_{\rm Hermite}$ and
$Z^{(\nu, N)}_{\rm Laguerre}$. 

The DPP $(\Xi, K^{(N)}_{\rm Hermite}, \lambda_{\rN(0,1/2)}(dy))$
describes
the eigenvalue distribution of 
$N \times N$ Hermitian random 
matrices in the {\it Gaussian unitary ensemble} (GUE).
When $\nu \in \N_0$, 
the DPP 
$(\Xi, K^{(N)}_{\rm Laguerre}$, $\lambda_{\Gamma(\nu+1, 1)}(dx))$
describes the distribution of the nonnegative square roots
of eigenvalues of $M^{\dagger} M$,
where $M$ is $(N+\nu) \times N$  complex
random matrix in the 
{\it chiral Gaussian ensemble} (chGUE)
and $M^{\dagger}$ is its Hermitian conjugate.
The probability density (\ref{eqn:P_Laguerre}) 
can be extended to any $\nu \in (-1, \infty)$ and it
is called the {\it complex Laguerre ensemble}
or the {\it complex Wishart ensemble}.
Many other examples of one-dimensional DPPs 
are given as eigenvalue ensembles
of Hermitian random matrices in the literatures of
random matrix theory (see, for instance, \cite{Meh04,For10,KT04}).

\begin{rem}
\label{thm:log1}
It should be noted that we can regard (\ref{eqn:P_Laguerre})
as the Gibbs measures
\begin{align*}
\bp^{(N)}_{\rm Hermite}(\x) 
&= \frac{1}{Z^{(N)}_{\rm Hermite}}
e^{-\beta V^{(N)}_{\rm Hermite}}(\x), 
\quad \x \in \R^N,
\nonumber\\
\bp^{(\nu, N)}_{\rm Laguerre}(\y) 
&= \frac{1}{Z^{(\nu, N)}_{\rm Laguerre}}
e^{-\beta V^{(N)}_{\rm Laguerre}}(\y),
\quad \nu \in (-1, \infty), \quad \y \in \R_{\geq 0}^N,
\end{align*}
with the {\it inverse temperature} $\beta=2$ associated with
the potentials,
\begin{align*}
V^{(N)}_{\rm Hermite}(\x)
&= - \sum_{1 \leq j < k \leq N} \log |x_k-x_j|+ \frac{1}{2} \sum_{j=1}^N x_j^2,
\quad \x \in \R^N, 
\nonumber\\
V^{(N)}_{\rm Laguerre}(\y)
&= - \sum_{1 \leq j < k \leq N} \log|y_k-y_j|
+ \frac{1}{2} \sum_{j=1}^N (-\nu \log y_j +y_j),
\quad \nu \in (-1, \infty), \quad \y \in \R_{\geq 0}^N.
\end{align*}
Here the interactions are given by logarithmic two-body potentials.
\end{rem}

\subsubsection{Duality relations between 
DPPs in continuous and discrete spaces}
\label{sec:duality_application1}
We consider the simplified setting (\ref{eqn:W_n_x1}) of $W$
with $\Gamma=\N_0$.
If we set
$\Lambda_1=[r, \infty) \subset S_1=\R, r \in \R$ and
$\Lambda_2=\{0,1, \dots, N-1\} \subset S_2=\Gamma=\N_0, N \in \N$
in (\ref{eqn:K_A}), we obtain
\begin{align}
K_{\R}^{\{0,1, \dots, N-1 \}}(x, x')
&= \sum_{n=0}^{N-1} \varphi_n(x) \overline{\varphi_n(x')},
\quad x, x' \in \R,
\nonumber\\
K^{[r, \infty)}_{\N_0}(n, n')
&= \int_r^{\infty} \overline{\varphi_n(x)} \varphi_{n'}(x) \lambda_1(dx),
\quad n, n' \in \N_0.
\label{eqn:Kr_inf}
\end{align}
When $\lambda_1(dx)$ and $\{\varphi_n\}_{n \in \N_0}$
are given by (\ref{eqn:HermiteZ1}) or by (\ref{eqn:LaguerreZ1}),
the kernels (\ref{eqn:Kr_inf}) are given by
\begin{align*}
& K_{{\rm DHermite}^+(r)}(n, n')
= ( \pi 2^{n+n'} n! n'!)^{-1/2}
\int_r^{\infty} H_n(x) H_{n'}(x) e^{-x^2} dx
\nonumber\\
& \qquad = - ( \pi n! n'! 2^{n+n'+2} )^{-1/2} 
e^{-r^2} 
\frac{H_{n+1}(r) H_{n'}(r) - H_n(r) H_{n'+1}(r)}{n-n'},
\end{align*}
and, provided $r > 0$, 
\begin{align*}
& K_{{\rm DLaguerre}^+(r, \nu+1)}(n, n')
= \left( \frac{n ! n' !}{\Gamma(n+\nu+1) \Gamma(n'+\nu+1)} \right)^{1/2}
\int_r^{\infty} L^{(\nu)}_n(x) L^{(\nu)}_{n'}(x) x^{\nu} e^{-x} dx
\nonumber\\
& \qquad = \left( \frac{n ! n' !}{\Gamma(n+\nu+1) \Gamma(n'+\nu+1)} \right)^{1/2}
r^{\nu+1} e^{-r}
\frac{L^{(\nu+1)}_{n-1}(r) L^{(\nu)}_{n'}(r)
- L^{(\nu)}_{n}(r) L^{(\nu+1)}_{n'-1}(r)}{n-n'},
\end{align*}
with the convention that $L^{(\nu)}_{-1}(r)=0$, 
respectively (see Propositions 3.3 and 3.4 in \cite{BO17}).
Borodin and Olshanski called
the correlation kernels $K_{{\rm DHermite}^+(r)}$ and
$K_{{\rm DLaguerre}^+(r, \nu+1)}$ the {\it discrete Hermite kernel}
and the {\it discrete Laguerre kernel}, respectively \cite{BO17}.
Theorem~\ref{thm:duality2} gives 
\begin{equation}
\bP(\Xi_1^{\{0,1, \dots, N-1\}}([r, \infty))=m)
=\bP(\Xi_2^{[r, \infty)}(\{0,1, \dots, N-1\})=m),
\quad \forall m \in \N_0,
\label{eqn:duality_eq3}
\end{equation}
where LHS denotes the probability that
the number of points in the interval $[r, \infty)$
is $m$ for the $N$-point continuous DPP on $\R$ such as
$(\Xi_1, K^{(N)}_{\rm Hermite}, \lambda_{N(0, 1/2)}(dx))$
or $(\Xi_1, K^{(N)}_{\rm Laguerre}, \lambda_{\Gamma(\nu+1, 1)}(dx))$, 
$\nu \in (-1, \infty)$, 
while RHS does the probability that
the number of points in $\{0,1,\dots, N-1\}$ is $m$ for
the discrete DPP on $\N_0$ such as
$(\Xi_2, K_{{\rm DHermite}^{+}(r)})$
or $(\Xi_2, K_{{\rm DLaguerre}^{+}(r, \nu+1)})$, $\nu \in (-1, \infty)$. 
The {\it duality between continuous and discrete ensembles}
of Borodin and Olshanski (Theorem 3.7 in \cite{BO17}) is
a special case with $m=0$ of the equality (\ref{eqn:duality_eq3}).

\subsubsection{Finite DPPs in intervals
related with classical root systems}
\label{sec:Lie_group}

Let $N \in \N$ and consider the four types of
{\it classical root systems} denoted by $\AN, \BN, \CN$, and $\DN$.
We set
$S^{\AN}=\S^1=[0, 2 \pi)$, the unit circle, 
with a uniform measure 
$\lambda^{\AN}(dx)=\lambda_{[0, 2 \pi)}(dx) :=dx/(2 \pi)$, and
$S^{\RN}=[0, \pi]$, the upper half-circle, with 
$\lambda^{\RN}(dx)=\lambda_{[0, \pi]}(dx) :=dx/\pi$ for
$\RN=\BN, \CN, \DN$. 

For a fixed $N \in \N$, we introduce the four sets of
functions $\{\varphi^{\RN}_n\}_{n=1}^N$ on $S^{\RN}$ defined as
\begin{align*}
\varphi^{\RN}_n(x) = \begin{cases}
\displaystyle{
e^{-i(\cN^{\AN}-2 J^{\AN}(n))x/2}
}, 
& \RN=\AN,
\cr
\displaystyle{
\sin \big[(\cN^{\RN}-2J^{\RN}(n)) x/2 \big]
}, 
& \RN=\BN, \CN, 
\cr
\displaystyle{
\cos \big[(\cN^{\DN}-2J^{\DN}(n)) x/2 \big]
}, 
& \RN=\DN,
\end{cases}
\end{align*}
where
\begin{equation}
\cN^{\RN} = \begin{cases}
N, \quad & \RN = \AN, \\
2N-1, \quad & \RN = \BN, \\
2(N+1), \quad & \RN =\CN, \\
2(N-1), \quad & \RN = \DN.
\end{cases}
\label{eqn:N_R_classic}
\end{equation}
and
\begin{equation}
J^{\RN}(n) = \begin{cases}
n-1/2, \quad & \RN = \AN, 
\\
n-1, \quad & \RN=\BN, \DN,
\\
n, \quad & \RN=\CN.
\end{cases}
\label{eqn:J_R_classic}
\end{equation}
It is easy to verify that they satisfy the following
orthonormality relations,
\begin{align*}
\langle \varphi^{\AN}_n, \varphi^{\AN}_m \rangle_{L^2(\S^1, \lambda_{[0, 2 \pi)})}
&=\delta_{n m}, 
\nonumber\\
\langle \varphi^{\RN}_n, \varphi^{\RN}_m \rangle_{L^2([0, \pi], \lambda_{[0, \pi]})}
&=\delta_{n m}, 
\quad \RN=\BN, \CN, \DN, \quad
\mbox{if $n, m \in \{1, \dots, N\}$}.
\end{align*}
We put $\Gamma = \{1, \dots, N\}, N \in \N$ and
$L^2(\Gamma, \nu) = \ell^2(\Gamma)$.
By the argument given in Remark \ref{thm:Remark1_2}, 
Assumption 3 is verified, and hence
Corollary \ref{thm:main2} gives the four types of DPPs; 
$(\Xi, K^{\AN}, \lambda_{[0, 2 \pi)}(dx))$ on $\S^1$,
and 
$(\Xi, K^{\RN}, \lambda_{[0, \pi]}(dx))$ on $[0, \pi]$,
$\RN=\BN, \CN, \DN$, with the correlation kernels,
\begin{align*}
K^{\RN}(x, x^{\prime})
&= \sum_{n=1}^N \varphi^{\RN}_n(x) \overline{\varphi^{\RN}_n(x^{\prime})}
\nonumber\\
&= \begin{cases}
\displaystyle{
\frac{\sin\{N(x-x^{\prime})/2\}}{\sin\{(x-x^{\prime})/2\}}
},
& \RN=\AN,
\cr
\displaystyle{ \frac{1}{2} \left[
\frac{\sin\{N(x-x^{\prime})\}}{\sin\{(x-x^{\prime})/2\}}
- \frac{\sin\{N(x+x^{\prime})\}}{\sin\{(x+x^{\prime})/2\}}
\right]},
& \RN=\BN, 
\cr
& \cr
\displaystyle{\frac{1}{2} \left[
\frac{\sin\{(2N+1)(x-x^{\prime})/2\}}{\sin\{(x-x^{\prime})/2\}}
- \frac{\sin\{(2N+1)(x+x^{\prime})/2\}}{\sin\{(x+x^{\prime})/2\}}
\right]},
& \RN=\CN, 
\cr
& \cr
\displaystyle{\frac{1}{2} \left[
\frac{\sin\{(2N-1)(x-x^{\prime})/2\}}{\sin\{(x-x^{\prime})/2\}}
+ \frac{\sin\{(2N-1)(x+x^{\prime})/2\}}{\sin\{(x+x^{\prime})/2\}}
\right]},
& \RN=\DN.
\cr
\end{cases}
\end{align*}

The {\it Weyl denominator formulas} for classical root systems
play a fundamental role in Lie theory and related area.
For a reduced root systems they are given in the form,
\[
\sum_{w \in W} \det(w) e^{w(\rho)-\rho}
=\prod_{\alpha \in R_{+}}(1-e^{-\alpha}),
\]
where $W$ is the Weyl group, $R_+$ the set of positive roots
and $\rho=\frac{1}{2} \sum_{\alpha \in R_+} \alpha$.
For classical root systems $\AN, \BN, \CN$ and $\DN$, 
$N \in \N$, 
the explicit forms are given as follows,
\begin{align}
\mbox{ (type $\AN$)}
\quad & \det_{1 \leq j, k \leq N} \Big(z_k^{j-1} \Big)
= \prod_{1 \leq j < k \leq N} (z_k-z_j),
\nonumber\\
\mbox{ (type $\BN$)}
\quad &
\det_{1 \leq, j, k \leq N} 
\Big( z_k^{j-N}-z_k^{N+1-j} \Big)
= \prod_{\ell=1}^{N} z_{\ell}^{1-N} (1-z_{\ell}) 
\prod_{1 \leq j < k \leq N} (z_k-z_j)(1-z_j z_k),
\nonumber\\
\mbox{ (type $\CN$)}
\quad &
\det_{1 \leq, j, k \leq N} 
\Big( z_k^{j-N-1}-z_k^{N+1-j} \Big)
= \prod_{\ell=1}^{N} z_{\ell}^{-N} (1-z_{\ell}^2) 
\prod_{1 \leq j < k \leq N} (z_k-z_j)(1-z_j z_k),
\nonumber\\
\mbox{ (type $\DN$)}
\quad &
\det_{1 \leq, j, k \leq N} 
\Big( z_k^{j-N}+z_k^{N-j} \Big)
= 2 \prod_{\ell=1}^{N} z_{\ell}^{1-N}
\prod_{1 \leq j < k \leq N} (z_k-z_j)(1-z_j z_k), 
\label{eqn:Weyl_denominator}
\end{align}
respectively. 
See, for instance, \cite{RS06}.
If we change the variables as
\begin{equation}
z_k=e^{-2 \sqrt{-1} \zeta_k}, \quad \zeta_k \in \C, \quad k=1, \dots, N,
\label{eqn:x_to_z}
\end{equation}
then, the following equalities are derived from
the above.

\begin{lem}
\label{thm:Weyl_trigonometric}
For $\zeta_k \in \C, k=1, \dots,N$, 
the following equalities are established.
\begin{align*}
{\rm (type} \, \AN {\rm )}
\quad &
\det_{1 \leq j, k \leq N} 
\Big[ e^{- \sqrt{-1} (\cN^{\AN}-2J^{\AN}(j)) \zeta_k} \Big]
= (2i) ^{N(N-1)/2}
\prod_{1 \leq j < k \leq N} 
\sin(\zeta_k-\zeta_j).
\nonumber\\
{\rm (type} \, \BN {\rm )}
\quad &
\det_{1 \leq j, k \leq N} 
\Big[ \sin\{(\cN^{\BN}-2 J^{\BN}(j)) \zeta_k\} \Big]
\nonumber\\
& \qquad
= 2^{N(N-1)}
\prod_{\ell=1}^N \sin \zeta_{\ell}
\prod_{1 \leq j < k \leq N} 
\sin(\zeta_k-\zeta_j) \sin(\zeta_k+\zeta_j),
\nonumber\\
{\rm (type} \, \CN {\rm )}
\quad &
\det_{1 \leq j, k \leq N} 
\Big[ \sin \{ (\cN^{\CN}-2 J^{\CN}(j)) \zeta_k\} \Big]
\nonumber\\
& \qquad
= 2^{N(N-1)}
\prod_{\ell=1}^N \sin(2 z_{\ell}) 
\prod_{1 \leq j < k \leq N} 
\sin (\zeta_k-\zeta_j) \sin (\zeta_k+\zeta_j),
\nonumber\\
{\rm (type} \, \DN {\rm )}
\quad &
\det_{1 \leq j, k \leq N} 
\Big[ \cos \{(\cN^{\DN}-2 J^{\DN}(j)) \zeta_k \} \Big]
\nonumber\\
& \qquad
= 2^{(N-1)^2}
\prod_{1 \leq j < k \leq N} 
\sin(\zeta_k-\zeta_j) \sin (\zeta_k+\zeta_j),
\end{align*}
where $\cN^{\RN}$ and $J^{\RN}(j)$, 
$\RN=\AN, \BN, \CN, \DN$, are given by
(\ref{eqn:N_R_classic}) and (\ref{eqn:J_R_classic}). 
\end{lem}

By Lemma \ref{thm:Weyl_trigonometric}, 
the probability densities for these finite DPPs
with respect to 
the Lebesgue measures,
$d \x =\prod_{j=1}^N d x_j$ are given
as
\begin{align*}
\bp^{\AN}(\x) 
&= \frac{1}{Z^{\AN}} \prod_{1 \leq j < k \leq N}
\sin^2 \frac{x_k-x_j}{2}, 
\quad \x \in [0, 2\pi)^N,
\nonumber\\
\bp^{\BN}(\x) 
&= \frac{1}{Z^{\BN}} 
\prod_{1 \leq j < k \leq N} 
\left( \sin^2 \frac{x_k-x_j}{2} \sin^2 \frac{x_k+x_j}{2} \right)
\prod_{\ell=1}^N \sin^2 \frac{x_{\ell}}{2},
\quad \x \in [0, \pi]^N,
\nonumber\\
\bp^{\CN}(\x) 
&= \frac{1}{Z^{\CN}} 
\prod_{1 \leq j < k \leq N} 
\left( \sin^2 \frac{x_k-x_j}{2} \sin^2 \frac{x_k+x_j}{2} \right)
\prod_{\ell=1}^N \sin^2 x_{\ell},
\quad \x \in [0, \pi]^N,
\nonumber\\
\bp^{\DN}(\x) 
&= \frac{1}{Z^{\DN}} 
\prod_{1 \leq j < k \leq N} 
\left( \sin^2 \frac{x_k-x_j}{2} \sin^2 \frac{x_k+x_j}{2} \right),
\quad \x \in [0, \pi]^N,
\end{align*}
with the normalization constants $Z^{\RN}$.

The DPP, $(\Xi, K^{\AN}, \lambda_{[0, 2 \pi)}(dx))$ is known as
the {\it circular unitary ensemble} (CUE)
in random matrix theory (see Section 11.8 in \cite{Meh04}).
These four types of DPPs, 
$(\Xi, K^{\AN}, \lambda_{[0, 2 \pi)}(d x))$,
$(\Xi, K^{\RN}, \lambda_{[0, \pi]}(d x))$, 
$\RN=\BN, \CN, \DN$
are realized as the eigenvalue distributions of
random matrices in the {\it classical groups}, 
$\U(N)$, $\SO(2N+1)$, $\Sp(N)$, and $\SO(2N)$, respectively.
(See Section 2.3 c) in \cite{Sos00}
and Section 5.5 in \cite{For10}.)

As mentioned in Remark \ref{thm:log1}, 
$\bp^{R_N}(\x)$ can be regarded as Gibbs measures
$e^{-\beta V^{R_N}(\x)}/Z^{R_N}$
with $\beta=2$, $R_N=\AN, \BN, \CN, \DN$.
For example, for type $\AN$
the potential is given as
\[
V^{\AN}(\x)=-\sum_{1 \leq j < k \leq N}
\log \left|\sin \frac{x_k-x_j}{2} \right|, \quad \x \in (0, 2 \pi]^N.
\]

\subsubsection{Infinite DPPs in $\R$ 
associated with classical orthonormal functions} 
\label{sec:classical_functions}

Here we give examples of infinite DPPs obtained
by Corollary \ref{thm:main3}.

\begin{description}
\item{(i)} \quad DPP with the {\it sinc kernel} : 
We set $S_1=\R$, $\lambda_1(dx)=dx$, 
$\Gamma=(-1, 1)$, $\nu(dy)=\lambda_2(d y)=dy$, and put
\[
\psi_1(x, y)=\frac{1}{\sqrt{2 \pi}} e^{ \sqrt{-1} x y} .
\]
We see that 
\[
\Psi_1(x)^2 :=\| \psi_1(x, \cdot) \|^2= \frac{1}{\pi}, 
\quad x \in \R,
\] 
and thus Assumption 3'(ii) is satisfied.
The correlation kernel $K_{S_1}$ is given by
\[
K_{\rm sinc}(x, x^{\prime})
= \frac{1}{2 \pi}\int_{-1}^1 e^{\sqrt{-1} y(x-x^{\prime})} dy
= \frac{\sin(x-x^{\prime}) }{ \pi (x-x^{\prime}) },
\quad x, x^{\prime} \in \R.
\]

\item{(ii)} \quad DPP with the {\it Airy kernel} :
We set $S_1=\R$, $\lambda_1(dx)=dx$, 
$\Gamma=\R_{\geq 0}$, $\nu(dy)=\lambda_2(d y)=dy$, 
and put
\[
\psi_1(x,y) =\Ai(x+y),
\]
where $\Ai(x)$ denotes the {\it Airy function} \cite{NIST10}
\[
\Ai(x)=\frac{1}{\pi} \int_0^{\infty} \cos \left( \frac{k^3}{3}+kx \right) dk.
\]
We see that 
\[
\Psi_1(x)^2 :=\| \psi_1(x, \cdot) \|^2
=-x \Ai(x)^2 + \Ai'(x)^2, \quad x \in \R,
\]
and thus Assumption 3(ii) is satisfied. 
The correlation kernel $K_{S_1}$ is given by
\[
K_{\rm Airy}(x, x^{\prime})
= \int_0^{\infty} \Ai(x+y) \Ai(x^{\prime}+y) dy
= \frac{\Ai(x) \Ai^{\prime}(x^{\prime}) - \Ai(x^{\prime}) \Ai^{\prime}(x)}
{x-x^{\prime}},
\quad x, x^{\prime} \in \R,
\]
where $\Ai^{\prime}(x)=d \Ai(x)/dx$.

\item{(iii)} \quad DPP with the {\it Bessel kernel} :
We set
$S_1=[0,\infty)$, $\lambda_1(dx)= dx$,
$\Gamma=[0, 1]$,
$\nu(dy)=\lambda_2(d y)=dy$.
With parameter $\nu \in (-1, \infty)$ we put
\[
\psi_1(x,y) =\sqrt{x y} J_{\nu}(x y),
\]
where $J_{\nu}$ is the 
{\it Bessel function of the first kind} defined by
\begin{equation}
J_{\nu}(x)= \sum_{n=0}^{\infty} 
\frac{(-1)^n}{n! \Gamma(\nu+n+1)}
\left(\frac{x}{2} \right)^{2n+\nu},
\quad x \in \C \setminus (-\infty, 0).
\label{eqn:Bessel_function}
\end{equation}
We see that 
\[
\Psi_1(x)^2 := \|\psi_1(x, \cdot) \|^2
= x \{J_{\nu}(x)^2 - J_{\nu-1}(x)J_{\nu+1}(x)\}/2,
\quad x \in \R_{\geq 0}, 
\]
and thus Assumption 3'(ii) is satisfied.
The correlation kernel $K_{S_1}$ is given by
\begin{align}
K_{\rm Bessel}^{(\nu)}(x, x')
&= \int_0^1 \sqrt{x y} J_{\nu}(x y) \sqrt{x' y} J_{\nu}(x' y) dy 
\nonumber\\
&= \frac{\sqrt{x x'}}{x^2-(x')^2}
 \Big\{ J_{\nu}(x) x' J_{\nu}^{\prime}(x')
 - x J_{\nu}^{\prime}(x) J_{\nu}(x') \Big\},
 \quad x, x' \in \R_{\geq 0},
 \label{eqn:Bessel_kernel}
\end{align}
where $J_{\nu}^{\prime}(x)=d J_{\nu}(x)/dx$.
\end{description}

These three kinds of infinite DPPs,
$(\Xi, K_{\rm sinc}, d x)$, $(\Xi, K_{\rm Airy}, d x)$,
and $(\Xi, K_{\rm Bessel}^{(\nu)}, {\bf 1}_{\R_{\geq 0}}(x) d x)$,
are obtained as the scaling limits
of the finite DPPs,
$(\Xi, K^{(N)}_{\rm Hermite}, \lambda_{\rN(0, 1/2)}(d x))$ 
and $(\Xi, K^{(\nu, N)}_{\rm Laguerre}$, $\lambda_{\Gamma(\nu+1, 1)}(d x))$, 
given in Section \ref{sec:classical}
as follows.

\begin{description}
\item{(i)} {\it Bulk scaling limit},
\[
\sqrt{2N} \circ
(\Xi, K^{(N)}_{\rm Hermite}, \lambda_{\rN(0, 1/2)}(d x))
\weak (\Xi, K_{\rm sinc}, d x). 
\]

\item{(ii)} {\it Soft-edge scaling limit},
\[
\sqrt{2}N^{1/6} \circ \cS_{\sqrt{2N}} \,
(\Xi,  K^{(N)}_{\rm Hermite}, \lambda_{\rN(0, 1/2)}(d x))
\weak (\Xi, K_{\rm Airy}, d x). 
\]

\item{(iii)} {\it Hard-edge scaling limit},
for $\nu \in (-1, \infty)$, 
\[
4N \circ
\Big(
(\Xi, K^{(\nu, N)}_{\rm Laguerre}, \lambda_{\Gamma(\nu+1,1)}(d x))^{\langle 1/2 \rangle}
\Big) 
\weak (\Xi, K^{(\nu)}_{\rm Bessel}, {\bf 1}_{\R_{\geq 0}}(x) d x).
\]
\end{description}

\noindent
See, for instance, \cite{Meh04,For10,AGZ10,Kat15_Springer}, 
for more details.

The DPPs with the sinc kernel and 
the Bessel kernel with the special values of
parameter $\nu$ can be obtained as the bulk
scaling limits of the DPPs,
$(\Xi, K^{\RN}, \lambda^{\RN}(d x))$,
$\RN=\AN, \BN, \CN, \DN$
given in Section \ref{sec:Lie_group} as
\begin{align}
\frac{N}{2} \circ (\Xi, K^{\AN}, \lambda_{[0, 2 \pi)}(d x))
&\weak (\Xi, K_{\rm sinc}, d x),
\nonumber\\
\left.
\begin{array}{l}
N \circ (\Xi, K^{\BN}, \lambda_{[0, \pi]}(d x))
\cr
N \circ (\Xi, K^{\CN}, \lambda_{[0, \pi]}(d x))
\end{array}
\right\}
&\weak (\Xi, K_{\rm Bessel}^{(1/2)}, {\bf 1}_{\R_{\geq 0}}(x) d x),
\nonumber\\
N \circ (\Xi, K^{\DN}, \lambda_{[0, \pi]}(d x))
&\weak (\Xi, K_{\rm Bessel}^{(-1/2)}, {\bf 1}_{\R_{\geq 0}}(x) d x),
\label{eqn:Lie_bulk}
\end{align}
where 
\begin{align*}
K_{\rm Bessel}^{(1/2)}(x, x^{\prime})
&= \frac{\sin (x-x^{\prime})}{\pi (x-x^{\prime})} 
- \frac{\sin(x+x^{\prime})}{\pi (x+x^{\prime})}, 
\quad x, x^{\prime} \in \R_{\geq 0}, 
\nonumber\\
K_{\rm Bessel}^{(-1/2)}(x, x^{\prime})
&= \frac{\sin(x-x^{\prime})}{\pi (x-x^{\prime})} 
+ \frac{\sin(x+x^{\prime})}{\pi (x+x^{\prime})}, 
\quad x, x^{\prime} \in \R_{\geq 0}.
\end{align*}
Since $J_{1/2}(x)=\sqrt{2/(\pi x)} \sin x$ and
$J_{-1/2}(x)=\sqrt{2/(\pi x)} \cos x$, 
the above correlation kernels 
are readily obtained from (\ref{eqn:Bessel_kernel})
by setting $\nu=1/2$ and $-1/2$, respectively.

\subsection
{Examples in two-dimensional spaces} \label{sec:examples_2d}
\subsubsection{Infinite DPPs on $\C$ : Ginibre and Ginibre-type DPPs}
\label{sec:Ginibre_ACD}
Let 
$\lambda_{\rN(m, \sigma^2; \C)}(d x)$ denote
the {\it complex normal distribution}, 
\begin{align*}
\lambda_{\rN(m, \sigma; \C)}(d x)
&:= \frac{1}{\pi \sigma^2} e^{-|x-m|^2/\sigma^2} d x
\nonumber\\
&= \frac{1}{\pi \sigma^2} 
e^{-(x_{\rR}-m_{\rR})^2/\sigma^2 -(x_{\rI}-m_{\rI})^2/\sigma^2} 
d x_{\rR}  d x_{\rI}, 
\end{align*}
$m \in \C, m_{\rR} :=\Re m, m_{\rI} :=\Im m, \sigma>0$.
We set
$S=\C$,
\begin{align*}
\lambda(dx) &=\lambda_{\rN(0,1;\C)}(dx)
= \frac{1}{\pi} e^{-|x|^2} dx
\nonumber\\
&=\lambda_{\rN(0,1/2)}(dx_{\rR}) \lambda_{\rN(0,1/2)}(dx_{\rI}),
\end{align*}
and 
\begin{align*}
\psi^{A}(x, \gamma) &=
e^{-(x_{\rR}^2-x_{\rI}^2)/2+2 x \gamma}, 
\nonumber\\
\psi^{C}(x, \gamma) &=\sqrt{2} 
\sinh(2 x \gamma)
e^{-(x_{\rR}^2-x_{\rI}^2)/2},
\nonumber\\
\psi^{D}(x, \gamma) &=\sqrt{2}
\cosh(2 x \gamma)
e^{-(x_{\rR}^2-x_{\rI}^2)/2}.
\end{align*}
It is easy to confirm that
\begin{align*}
\frac{1}{\pi} 
\int_{\R} \psi^A(x, \gamma) \overline{\psi^A(x, \gamma')} e^{-x_{\rI}^2} d x_{\rI}
&=e^{-(x_{\rR}^2-4 x_{\rR} \gamma)} \delta(\gamma-\gamma'),
\cr
\frac{1}{\pi}
\int_{\R} \psi^R(x, \gamma) \overline{\psi^R(x, \gamma')} e^{- x_{\rI}^2} d x_{\rI}
&= e^{- x_{\rR}^2} 
\cosh(4 x_{\rR} \gamma)
\times 
\begin{cases}
\delta(\gamma-\gamma')-\delta(\gamma+\gamma'), & R=C, \cr
\delta(\gamma-\gamma')+\delta(\gamma+\gamma'), & R=D.
\end{cases}
\end{align*}
Therefore, we have
\begin{align*}
& \langle \psi^A(\cdot, \gamma), \psi^A(\cdot, \gamma') \rangle_{L^2(\C, \lambda_{\rN(0, 1; \C)})}
\nu(d \gamma)
=\delta(\gamma-\gamma') d \gamma, \quad
\gamma, \gamma' \in \Gamma^A :=\R,
\cr
& \langle \psi^R(\cdot, \gamma), \psi^R(\cdot, \gamma') \rangle_{L^2(\C, \lambda_{\rN(0,1; \C)})}
\nu(d \gamma)
=\delta(\gamma-\gamma') d \gamma, \quad
\gamma, \gamma' \in \Gamma^R :=(0, \infty), \quad R=C, D,
\end{align*}
with $\nu(d \gamma)=\lambda_{\rN(0, 1/4)}(d \gamma)$. 
We also see that
\begin{align*}
\Psi^A(x)^2 &:=\| \psi^A(x, \cdot) \|_{L^2(\Gamma^A, \nu)}^2
=e^{|x|^2},
\nonumber\\
\Psi^C(x)^2 &:=\| \psi^C(x, \cdot) \|_{L^2(\Gamma^C, \nu)}^2
=\sinh |x|^2,
\nonumber\\
\Psi^D(x)^2 &:=\| \psi^D(x, \cdot) \|_{L^2(\Gamma^D, \nu)}^2
=\cosh |x|^2, \quad x \in \C.
\end{align*}
Thus Assumption 3 is satisfied and we can
apply Corollary \ref{thm:main2}. 
The kernels (\ref{eqn:K_main2}) of obtained DPPs are given as
\begin{align*}
K^A(x, x^{\prime})
&= \sqrt{\frac{2}{\pi}} e^{-\{(x_{\rR}^2-x_{\rI}^2)+({x'_{\rR}}^2-{x'_{\rI}}^2) \}/2}
\int_{-\infty}^{\infty} 
e^{- 2 \{\gamma^2-(x+\overline{x'}) \gamma\}} d\gamma,
\nonumber\\
K^C(x, x^{\prime})
&= 2 \sqrt{\frac{2}{\pi}} e^{-\{(x_{\rR}^2-x_{\rI}^2)+({x'_{\rR}}^2-{x'_{\rI}}^2) \}/2}
\int_0^{\infty} 
e^{-2 \gamma^2}
\sinh(2 x \gamma) \sinh(2 \overline{x^{\prime}} \gamma) d\gamma, 
\nonumber\\
K^D(x, x^{\prime})
&= 2 \sqrt{\frac{2}{\pi}} e^{-\{(x_{\rR}^2-x_{\rI}^2)+({x'_{\rR}}^2-{x'_{\rI}}^2) \}/2}
\int_0^{\infty} 
e^{-2 \gamma^2}
\cosh(2 x \gamma) \cosh(2 \overline{x^{\prime}} \gamma) d\gamma. 
\end{align*}
The integrals are performed and we obtain 
$
K^R(x, x^{\prime})
= e^{\sqrt{-1} x_{\rR} x_{\rI}} K^{R}_{{\rm Ginibre}}(x, x^{\prime})
e^{-\sqrt{-1} x_{\rR}^{\prime} x_{\rI}^{\prime}}
$,
$R=A, C, D$, 
with
\begin{align}
K^{A}_{{\rm Ginibre}}(x, x^{\prime})
&= e^{x \overline{x^{\prime}}},
\label{eqn:K_Ginibre_A}
\\
K^{C}_{{\rm Ginibre}}(x, x^{\prime})
&= \sinh(x \overline{x^{\prime}}),  
\label{eqn:K_Ginibre_C}
\\
K^{D}_{{\rm Ginibre}}(x, x^{\prime})
&= \cosh(x \overline{x^{\prime}}), \quad x, x^{\prime} \in \C. 
\label{eqn:K_Ginibre_D}
\end{align}
Due to the gauge invariance of DPP mentioned in Section \ref{sec:basic_DPP}, 
the obtained three types of infinite DPPs on $\C$ are written as
$(\Xi, K_{\rm Ginibre}^R, \lambda_{\rN(0,1; \C)}(dx))$, $R=A, C, D$.
The DPP, $(\Xi, K_{\rm Ginibre}^A, \lambda_{\rN(0,1; \C)}(dx))$ with (\ref{eqn:K_Ginibre_A})
describes the eigenvalue distribution of the
Gaussian random complex matrix in the bulk scaling limit,
which is called the {\it complex Ginibre ensemble}
\cite{Gin65,Meh04,HKPV09,For10,Shirai15}.
This density profile is uniform 
with the Lebesgue measure $dx$ on $\C$ as
\[
\rho_{{\rm Ginibre}}(x) dx
=K^{A}_{{\rm Ginibre}}(x, x) \lambda_{\rN(0, 1; \C)}(d x)
=\frac{1}{\pi} dx_{\rR} dx_{\rI},
\quad x \in \C.
\]
On the other hands, the Ginibre DPPs
of types $C$ and $D$ 
with the correlation kernels (\ref{eqn:K_Ginibre_C})
and (\ref{eqn:K_Ginibre_D}) 
are rotationally symmetric around the origin, 
but non-uniform on $\C$.
The density profiles 
with the Lebesgue measure $dx$ on $\C$ are given by
\begin{align*}
\rho^{C}_{{\rm Ginibre}}(x) dx
&=
K^{C}_{{\rm Ginibre}}(x, x) \lambda_{\rN(0, 1; \C)}(d x)
=\frac{1}{2 \pi} (1-e^{-2 |x|^2}) dx_{\rR} dx_{\rI}, \quad x \in \C, 
\nonumber\\
\rho^{D}_{{\rm Ginibre}}(x) dx
&=K^{D}_{{\rm Ginibre}}(x, x) \lambda_{\rN(0, 1; \C)}(d x)
=\frac{1}{2 \pi} (1+e^{-2 |x|^2})dx_{\rR} dx_{\rI}, \quad x \in \C.
\end{align*}
They were first obtained in \cite{Kat19b} by taking
the limit $W \to \infty$ keeping the density of points of the infinite DPPs 
in the strip on $\C$, $\{z \in \C : 0 \leq \Im z \leq W\}$. 

\subsubsection{Representations of Ginibre and Ginibre-type 
kernels in the Bargmann--Fock space and the eigenspaces of Landau levels}
\label{sec:Fock}

One of the advantages of our framework is that 
the we can obtain pairs of DPPs satisfy 
useful duality relations.
Now we concentrate on one of a pair of DPPs
constructed in our framework, $(\Xi_1, K_{S_1}, \lambda_1)$.
The correlation kernel $K_{S_1}$ is given by the first equation
of (\ref{eqn:K1}), that is, 
\[
K_{S_1}(x, x') = \int_{S_2} \overline{W(y, x)} W(y, x') \lambda_2(dy)
=\langle W(\cdot, x'), W(\cdot, x) \rangle_{L^2(S_2, \lambda_2)},
\quad x, x' \in S_1, 
\]
which is an integral kernel for $f \in L^2(S_1, \lambda_1)$.
We can regard this equation as a {\it decomposition formula}
of $K_{S_1}$ by a product of $W$ and $\overline{W}$.
Since $W$ is an integral kernel for an isometry
$L^2(S_1, \lambda_1) \to L^2(S_2, \lambda_2)$, 
as a matter of course, it depends on a choice of
another Hilbert space $L^2(S_2, \lambda_2)$.
We note that a given DPP, $(\Xi_1, K_{S_1}, \lambda_1)$, 
choice of $L^2(S_2, \lambda_2)$ is not unique. 
Such multivalency gives plural different expressions
for one correlation kernel $K_{S_1}$ and they reveal
different aspects of the DPP.

Here we demonstrate this fact using the three kinds of
Ginibre DPPs associated with 
$L^2(\C, \lambda_{\rN(0,1;\C)})$.
In the previous section we have chosen the parameter spaces as
$L^2(\Gamma^R, \lambda_{\N(0,1/4)})$
with $\Gamma^A=\R$ and $\Gamma^C=\Gamma^D=(0, \infty)$.
We will choose another parameter spaces below.

Let
$S_1=\C$
and $S_2=\N_0$ with 
$\lambda_1(dx)=\lambda_{\rN(0, 1; \C)}(d x)$.
We put
\begin{equation}
\varphi_n(x)=\frac{x^n}{\sqrt{n!}},
\quad n \in \N_0.
\label{eqn:BF_phi}
\end{equation}
Note that $\{\varphi_n(x)\}_{n \in \N_0}$ forms 
a complete orthonormal system 
of the {\it Bargmann--Fock space}, 
which is the space of square-integrable analytic functions
on $\C$ with respect to the complex Gaussian
measure; 
\[
\langle \varphi_n, \varphi_m \rangle_{L^2(\C, \lambda_{\rN(0,1;\C)})}
=\delta_{n m}, \quad n, m \in \N_0.
\]
We assume that $\Gamma=S_2=\N_0$.
We can see that 
\[
\|\varphi_{\cdot}(x) \|_{\ell^2(\Gamma)}
=\sum_{n \in \N_0} |x|^{2n}/n!=e^{|x|^2}, 
\quad x \in \C.
\]
Hence Assumption 3' is satisfied.
By Corollary \ref{thm:main3}, 
we obtain the DPP on $\C$ in which
the correlation kernel 
with respect to $\lambda_{\rN(0, 1: \C)}$ is given by
\begin{align*}
K_{\rm BF}(x, x^{\prime}) 
&= \sum_{n \in \N_0} \varphi_n(x) \overline{\varphi_n(x^{\prime})}
=\sum_{n=0}^{\infty} \frac{(x \overline{x^{\prime}})^n}{n!}
\nonumber\\
&= e^{x \overline{x^{\prime}}},
\quad x, x^{\prime} \in \C.
\end{align*}
This is the reproducing kernel in the Bargmann--Fock space
and obtained DPP is identified with 
$(\Xi, K^A_{\rm Ginibre}, \lambda_{\rN(0,1; \C)}(dx))$.
See \cite{Shirai15,BQ17,AGR19}.

If we set $\Gamma=2 \N_0+1 =\{1,3,5,\dots\}$ or 
$\Gamma=2 \N_0=\{0,2,4,\dots\}$, 
we will obtain the DPPs with the following kernels 
\begin{align*}
K_{\rm BF}^{\rm odd}(x, x^{\prime})
&= \sum_{k=0}^{\infty} \frac{(x \overline{x^{\prime}})^{2k+1}}{(2k+1)!}
=\sinh(x \overline{x^{\prime}}), 
\nonumber\\
K_{\rm BF}^{\rm even}(x, x^{\prime})
&= \sum_{k=0}^{\infty} \frac{(x \overline{x^{\prime}})^{2k}}{(2k)!}
=\cosh(x \overline{x^{\prime}}),
\quad x, x^{\prime} \in \C. 
\end{align*}
The obtained DPPs are identified with
$(\Xi$, $K^C_{\rm Ginibre}$, $\lambda_{\rN(0,1; \C)}(dx))$ 
and $(\Xi$, $K^D_{\rm Ginibre}$, $\lambda_{\rN(0,1; \C)}(dx))$,
respectively.

The Ginibre DPP of type A is extended to 
Ginibre-type DPPs indexed by $q \in \N_0$,  
$(\Xi$, $K^{(q)}_{\mathrm{Ginibre\text{-}type}}, 
\lambda_{\rN(0,1; \C)}(dx))$, $q \in \N_0$, 
which are introduced in \cite{Shirai15} and 
also known as the infinite pure {\it polyanalytic ensembles} (cf. \cite{AGR19}). 
Each Ginibre-type DPP with index 
$q \in \N_0$ is associated with the correlation kernel 
\begin{equation}
 K^{(q)}_{\mathrm{Ginibre\text{-}type}}(x,x') := 
L_q^{(0)}(|x-x'|^2) K_{\mathrm{Ginibre}}^A(x, x'), 
\quad x, x' \in \C, 
\label{eqn:K_G_type}
\end{equation}
where $L_q^{(0)}$ is the $q$-th Laguerre polynomial (\ref{eqn:Laguerre1})
with parameter $\nu=0$ and $K_{\mathrm{Ginibre}}^A$ is defined by
(\ref{eqn:K_Ginibre_A}). 
The correlation kernel
(\ref{eqn:K_G_type}) admits the similar representation
in terms of the {\it complex Hermite polynomials} defined by 
\[
 H_{p,q}(\zeta, \zetabar) = (-1)^{p+q} e^{\zeta \zetabar} 
 \frac{\partial^p}{\partial \zetabar^p} \frac{\partial^q}{\partial \zeta^q} 
e^{-\zeta \zetabar}, \quad \zeta \in \C, 
\quad p, q \in \N_0, 
\]
which were introduced by It\^{o} \cite{Ito53}. 
We note that their generating function is given by 
\[
\sum_{p=0}^{\infty} \sum_{q=0}^{\infty} 
H_{p,q}(\zeta,\zetabar) \frac{s^p t^q}{p!q!} 
= \exp(\zeta s + \zetabar t - st)
\]
and the set $\{ H_{p,q}(\zeta, \zetabar)/\sqrt{p!q!} : p, q
\in \N_0\}$ forms a complete orthonormal system of $L^2(\C,
\lambda_{\rN(0, 1; \C)}(d\zeta))$. 
Let $S_1=\C$ and $S_2=\N_0$ with 
$\lambda_1(dx)=\lambda_{\rN(0, 1; \C)}(d x)$, 
and for fixed $q \in \N_0$, define
\[
\varphi^{(q)}_n(x) := \frac{1}{\sqrt{n!q!}} H_{n,q}(x, \overline{x}),
\quad x \in \C, 
\quad n \in \N_0.
\]
Then $\{\varphi^{(q)}_n(x)\}_{n \in \N_0}$ forms 
a complete orthonormal system of the eigenspace corresponding
to the $q$-th {\it Landau level}, which coincides with 
the Bargmann--Fock space when $q=0$. 
Since the following formula is known 
\[
L_q^{(0)}(|\zeta-\eta|^2) e^{\zeta \etabar}
= \sum_{p=0}^{\infty} \frac{1}{p!q!}
H_{p,q}(\zeta,\zetabar)
\overline{H_{p,q}(\eta,\etabar)}, 
\quad \zeta, \eta \in \C, 
\quad q \in \N_0, 
\]
we obtain the following expansion formula for (\ref{eqn:K_G_type}), 
\[
K^{(q)}_{\mathrm{Ginibre\text{-}type}}(x,x') 
= \sum_{n=0}^{\infty}  \varphi^{(q)}_n(x) \overline{\varphi^{(q)}_n(x')},
\quad x, x' \in \C, 
\quad q \in \N_0. 
\]

\subsubsection{Application of duality relations}
\label{sec:duality_application2}
We consider the simplified setting (\ref{eqn:W_n_x1}) of $W$
with (\ref{eqn:BF_phi}) and 
$\Gamma=\N_0$.  
If we set $\lambda_1(dx)=\lambda_{\rN(0,1;\C)}(dx)$, 
$\Lambda_1$ be a disk (i.e., two-dimensional 
ball) $\B^2_r$ with radius $r \in (0, \infty)$
centered at the origin in $S_1=\C \simeq \R^2$ and
$\Lambda_2=S_2=\N_0$ in (\ref{eqn:K_A}), we obtain
\begin{align*}
K_{\C}^{(\N_0)}(x, x')
&= \sum_{n=0}^{\infty} \varphi_n(x) \overline{\varphi_n(x')}
=e^{x \overline{x'}} 
\nonumber\\
&= K_{\rm Ginibre}^{A}(x, x'),
\quad x, x' \in \C,
\end{align*}
where $K_{\rm Ginibre}^A$ denotes the 
correlation kernel of the Ginibre DPP 
of type $A$, and
\begin{align*}
K^{(\B^2_r)}_{\N_0}(n, n')
&= \int_{\B^2_r} \overline{\varphi_n(x)} \varphi_{n'}(x) \lambda_{\rN(0,1;\C)}(dx)
= \frac{1}{\pi \sqrt{n! n'!}} \int_0^r ds \, e^{-s^2} s^{n+n'+1} 
\int_0^{2 \pi} d \theta \, e^{i \theta (n'-n)}
\nonumber\\
&= 2 \delta_{n n'} \int_0^r \frac{s^{2n+1} e^{-s^2} }{n!} ds
=\delta_{n n'} \int_0^{r^2} \lambda_{\Gamma(n+1, 1)} (du),
\quad n, n' \in \N_0.
\end{align*}
Define
\[
\lambda_n(r) :=
\int_0^{r^2} \frac{u^n e^{-u}}{n!} du
= \sum_{k=n+1}^{\infty} \frac{r^{2k} e^{-r^2}}{k!},
\quad n \in \N_0, \quad r \in (0, \infty),
\]
where the second equality is due to Eq.(4.1) in \cite{Shirai15}.
That is, if we write the Gamma distribution with parameters
$(a, b)$ as $\Gamma(a, b)$ (see Section \ref{sec:classical})
and the Poisson distribution with parameter $c$ as
${\rm Po}(c)$,
\[
\lambda_n(r) := \bP(R_n \leq r^2) = \bP(Y_{r^2} \geq n+1),
\]
provided $R_n \sim \Gamma(n+1, 1)$ and
$Y_{r^2} \sim {\rm Po}(r^2)$.
Then DPP 
$(\Xi_2^{(\B^2_r)}, K^{(\B^2_r)}_{\N_0})$ on $\N_0$
is the product measure 
$\bigotimes_{n \in \N_0} \mu_{\lambda_n(r)}^{\rm Bernoulli}$
under the natural identification between 
$\{0,1\}^{\N_0}$ and the power set of $\N_0$, 
where $\mu_p^{\rm Bernoulli}$ denotes the Bernoulli measure of probability $p \in [0, 1]$.
Theorem~\ref{thm:duality2} gives the duality relation
\[
\bP(\Xi_{\rm Ginibre}^A(\B^2_r)=m)
=
\bP(\Xi_2^{(\B^2_r)}(\N_0)=m), 
\quad
\forall m \in \N_0,
\]
where we have identified the DPP, 
$(\Xi_1^{(\N_0)}, K_{\C}^{(\N_0)}, \lambda_1(dx))$
with the Ginibre DPP of type A, 
$(\Xi_{\rm Ginibre}^A$, $K_{\rm Ginibre}^A, \lambda_{\rN(0, 1; \C)})$.
If we introduce a series of random variables
$X_n^{(r)} \in \{0, 1\}, n \in \N_0$,
which are mutually independent and 
$X_n^{(r)} \sim \mu_{\lambda_n(r)}^{\rm Bernoulli}$, $n \in \N_0$, 
then the above implies the equivalence in probability law
\[
\Xi_{\rm Ginibre}^A(\B^2_r) \law= 
\Xi_2^{(\B^2_r)}(\N_0) 
\law= \sum_{n \in \N_0} X_n^{(r)},
\quad r \in (0, \infty).
\]
Similarly, we have the following equalities
by the results in Section \ref{sec:Fock} and 
Theorem \ref{thm:duality2},
\[
\Xi_{\rm Ginibre}^C(\B^2_r) 
\law= \sum_{n \in 2\N_0+1} X_n^{(r)},
\quad
\Xi_{\rm Ginibre}^D(\B^2_r) 
\law= \sum_{n \in 2\N_0} X_n^{(r)},
\quad
r \in (0, \infty).
\]

The argument above is valid for general 
{\it radially symmetric DPPs}
associated with radially symmetric finite measure $\lambda_1(dx) =
p(|x|) dx$ on $\C$. 
Let 
$\varphi_n(x) = a_n x^n, n \in \N_0$ 
be an
orthonormal system in $L^2(\C, \lambda_1)$ where 
$a_n > 0, n \in \N_0$ are the normalization constants, 
and we set 
\begin{align*}
K_{\C}^{(\N_0)}(x, x')
&= \sum_{n=0}^{\infty} \varphi_n(x) \overline{\varphi_n(x')}
= \sum_{n=0}^{\infty} a_n^2 (x \overline{x'})^n 
\quad x, x' \in \C,
\nonumber\\
K^{(\B^2_r)}_{\N_0}(n, n')
&= \int_{\B^2_r} \overline{\varphi_n(x)} \varphi_{n'}(x)
 \lambda_1(dx)
=\delta_{n n'} \lambda_n(r) \quad n, n' \in \N_0, 
\end{align*}
where 
\[
\la_n(r) := \frac{1}{Z_n} \int_0^{r^2} u^n p(\sqrt{u}) du 
\]
with $Z_n = \int_0^{\infty} u^n p(\sqrt{u})du$. 
Then DPP 
$(\Xi_1^{(\N_0)}, K_{\C}^{(\N_0)}, p(|x|)dx)$ on $\C$
is radially symmetric 
and DPP 
$(\Xi_2^{(\B^2_r)}, K^{(\B^2_r)}_{\N_0})$ on $\N_0$
is again identified
with the product measure 
$\bigotimes_{n \in \N_0} \mu_{\lambda_n(r)}^{\rm
Bernoulli}$.  
For example, if $p(s) = \pi^{-1} e^{-s^2}$ and $a_n = 1/\sqrt{n!}$,
then 
$(\Xi_1^{(\N_0)}, K_{\C}^{(\N_0)}, p(|x|) dx)$
is the Ginibre DPP of type A. 
The function $\la_n(r)$ is considered as a probability distribution
function on $[0,\infty)$ and hence there exist independent random
variables $R_n, n \in \N_0$ such that 
\[
 \la_n(r) = \bP(R_n \le r^2). 
\]
If we define $X_n^{(r)} = \mathbf{1}_{\{R_n \le r^2\}}$ for
each $n \in
\N_0$, then Theorem \ref{thm:duality2} gives the duality relation
\[
\Xi_1^{(\N_0)}(\B^2_r) \law= \Xi_2^{(\B^2_r)}(\N_0)
\law= \sum_{n \in \N_0} X_n^{(r)},
\quad r \in (0, \infty).
\]
Indeed, $\{X_n^{(r)}, n \in \N_0\}$ 
are mutually independent $\{0,1\}$-valued random variables
whose laws are given by 
$\{\mu_{\lambda_n(r)}^{\rm Bernoulli}$, $n \in\N_0\}$. 
If we take a set $\La_2 \subset \N_0$, then DPP 
$(\Xi_1^{(\La_2)}, K_{\C}^{(\Lambda_2)}, p(|x|) dx)$
satisfies 
\[
\Xi_1^{(\La_2)}(\B^2_r) \law= \Xi_2^{(\B^2_r)}(\Lambda_2)
\law= \sum_{n \in \Lambda_2} X_n^{(r)},
\quad r \in (0, \infty).
\]
We note that if we write 
$\Xi_1^{(\N_0)} = \sum_{i} \delta_{x_i}$,
then $\sum_i \delta_{|x_i|^2}$ is equal to
$\sum_{n \in \N_0} \delta_{R_n}$ in law, which was 
discussed in Theorem 4.7.1 in \cite{HKPV09} by constructing
$\{R_n \}_{n \in \N_0}$ in
terms of {\it size-biased sampling}.  

\subsection{Examples in spaces with arbitrary dimensions} 
\label{sec:example_d_dim}
\subsubsection{Euclidean family of infinite DPPs on $\R^{d}$}
\label{sec:EuclideanDPP}

For $d \in \N$, let 
$S_1=S_2=\R^d$,
$\lambda_1(d x)=d x$, 
$\lambda_2(d y)=\nu(d y)=d y$,
$\psi_1(x, y)=e^{\sqrt{-1} x \cdot y}/(2\pi)^{d/2}$, and
$\Gamma=\B^d \subsetneq \R^d$,
where
$\B^d$ denotes the unit ball centered at the origin; 
$\B^d := \{ y \in \R^d : |y| \leq 1\}$.
We see 
\[
\Psi_1(x)^2 := \|\psi_1(x, \cdot)\|_{L^2(\Gamma, d \nu)}^2
=|\B^d|/(2 \pi)^d, \quad x \in \R^d, 
\]
where the volume of $\B^{d}$ is denoted by
$|\B^d| =\pi^{d/2}/\Gamma((d+2)/2)$.
Then Assumption 3' is satisfied and 
Corollary \ref{thm:main3} gives the DPP
in $S_1=\R^d$ whose correlation kernel
with respect to $\lambda_1(dx)=dx$ is given by
\begin{equation}
K^{(d)}(x, x')
= \frac{1}{(2\pi)^{d}}
\int_{\B^{d}} e^{\sqrt{-1} (x-x') \cdot y} d y.
\label{eqn:Fourier}
\end{equation}

The kernels $K^{(d)}$ on $\R^d, d \in \N$ have been
studied by Zelditch and others 
(see \cite{Zel00,SZ02,Zel09,CH15} and references therein),
who regarded them as the Szeg\"o kernels
for the reduced Euclidean motion group.
Here we call the DPPs associated with
the correlation kernels in this form 
the {\it Euclidean family of DPPs} on $\R^d, d \in \N$.
We can verify other expressions of $K^{(d)}, d \in \N$ 
using the Bessel function of the first kind
(\ref{eqn:Bessel_function}) as follows \cite{KS2}.

\begin{df}
\label{thm:EuclideanDPP}
The Euclidean family of DPP on $\R^d, d \in \N$ is defined by
$\Big( \Xi, K^{(d)}_{\rm Euclidean}, d x \Big)$
with the correlation kernel
\begin{align*}
K^{(d)}_{\rm Euclid}(x, x')
&= \frac{1}{(2 \pi)^d}
\int_{\R^d} {\bf 1}_{\B^d}(y)
e^{\sqrt{-1} (x-x') \cdot y} d y
= \frac{1}{(2 \pi)^d}
\int_{\B^d} 
e^{\sqrt{-1} (x-x') \cdot y} d y,
\nonumber\\
&=\frac{1}{(2 \pi)^{d/2}}
\frac{1}{|x-x'|^{(d-2)/2}}
\int_0^1 s^{d/2} J_{(d-2)/2} (|x-x'| s) ds
\nonumber\\
&= \frac{1}{(2 \pi)^{d/2}} 
\frac{J_{d/2}(|x-x'|)}{|x-x'|^{d/2}}, 
\quad x, x' \in \R^d.
\end{align*}
\end{df}

We see that
\[
K^{(d)}_{\rm Euclid}(x, x)
=\lim_{r \to 0} \frac{1}{(2\pi)^{d/2}} \frac{J_{d/2}(r)}{r^{d/2}}
=\frac{1}{2^d \pi^{d/2} \Gamma((d+2)/2)}.
\]
Then the Euclidean family of DPP is uniform on $\R^d$
with the density 
\[
\rho^{(d)}_{\rm Euclid}= \frac{1}{2^d \pi^{d/2} \Gamma((d+2)/2)}
\]
with respect to the Lebesgue measure $d x$ of $\R^d$. 

For lower dimensions, the correlation kernels 
and the densities are given as follows,
\begin{align*}
K^{(1)}_{\rm Euclid}(x, x')
&= \frac{\sin(x-x')}{\pi(x-x')}
=K_{\rm sinc}(x, x')
\quad
\mbox{with} \quad \rho^{(1)}_{\rm Euclid} = \frac{1}{\pi},
\nonumber\\
K^{(2)}_{\rm Euclid}(x, x')
&= \frac{J_1(|x-x'|)}{2 \pi |x-x'|}
\quad
\mbox{with} \quad \rho^{(2)}_{\rm Euclid} = \frac{1}{4 \pi},
\nonumber\\
K^{(3)}_{\rm Euclid}(x, x')
&= \frac{1}{2 \pi^2 |x-x'|^2}
\left( \frac{\sin |x-x'|}{|x-x'|}
-\cos |x-x'| \right)
\quad
\mbox{with} \quad \rho^{(3)}_{\rm Euclid} = \frac{1}{6 \pi^2}.
\end{align*}
This family of DPPs includes the DPP with the sinc kernel 
$K_{\rm sinc}$ as the
lowest dimensional case with $d=1$.
Note that, 
if $d$ is odd,
\[
K^{(d)}_{\rm Euclid}(x,x')
=\left(- \frac{1}{2\pi r} \frac{d}{dr} \right)^{(d-1)/2}
\frac{\sin r}{\pi r}
\quad \mbox{with} \quad
r=|x-x'|.
\]
This is proved by Rayleigh's formula
for the spherical Bessel function of the first kind
(Eq. (10.49.14) in \cite{NIST10});
\[
j_m(x) := \sqrt{\frac{\pi}{2 x}}
J_{m+1/2}(x)
=x^m \left( - \frac{1}{x} 
\frac{d}{dx} \right)^m \frac{\sin x}{x},
\quad m \in \N.
\]

\subsubsection{Heisenberg family of infinite DPPs on $\C^{d}$}
\label{sec:HeisenbergDPP}

The Ginibre DPP of type A on $\C$ given in Section \ref{sec:Ginibre_ACD}
can be generalized to the DPPs on $\C^d$ for $d \geq 2$.
This generalization was done by \cite{AGR16,APRT17,AGR19}
as the family of DPP
called the {\it Weyl--Heisenberg ensembles},
but here we derive the DPPs on $\C^d$, $d \in \N$,
following Corollary \ref{thm:main3} given in
Section \ref{sec:orthogonal}. 

Let $S_1=\C^d$, $S_2=\Gamma=\R^d$,  
\begin{align*}
\lambda_1(d x) & = \prod_{a=1}^d \lambda_{\rN(0, 1; \C)}(d x^{(a)})
= \frac{1}{\pi^d} e^{-|x|^2} = \frac{1}{\pi^d} e^{-(|x_{\rR}|^2+|x_{\rI}|^2)}
\nonumber\\
&=: \lambda_{\rN(0, 1; \C^d)}(d x),
\nonumber\\
\lambda_2(d y) &= \nu(dy) 
= \prod_{a=1}^d \lambda_{\rN(0, 1/4)}(d y^{(a)})
=\left( \frac{2}{\pi} \right)^{d/2} e^{-2 |y|^2}, 
\end{align*}
and
\[
\psi_1(x, y)
= e^{-(|x_{\rR}|^2-|x_{\rI}|^2)/2 + 2 (x_{\rR} \cdot y + \sqrt{-1} x_{\rI} \cdot y)},
\quad x=x_{\rR}+\sqrt{-1} x_{\rI} \in \C^d, \quad y \in \R^d.
\]
We see that 
\[
\Psi_1(x)^2:=\|\psi_1(x, \cdot) \|_{L^2(\R^d, \nu)}^2=e^{|x|^2},
\quad x \in \C^d. 
\]
Hence Assumptions 3' is satisfied, 
and then, by Corollary \ref{thm:main3}, 
we obtain the DPP on $\C^d$ with the correlation kernel,
\begin{align*}
K^{(d)}(x, x') 
&= \left( \frac{2}{\pi} \right)^{d/2}
e^{-\{(|x_{\rR}|^2-|x_{\rI}|^2)+(|x^{\prime}_{\rR}|^2-|x^{\prime}_{\rI}|^2)\}/2}
\int_{\R^d} e^{-2 [ |y|^2-\{(x_{\rR}+\sqrt{-1} x_{\rI}) 
+ (x^{\prime}_{\rR}-\sqrt{-1} x^{\prime}_{\rI}) \} \cdot y ]}
d y
\nonumber\\
&= \frac{e^{\sqrt{-1} x_{\rR} \cdot x_{\rI}}}{e^{\sqrt{-1} x^{\prime}_{\rR} \cdot x^{\prime}_{\rI}}}
K_{\rm Heisenberg}^{(d)}(x, x')
\end{align*}
with 
\[
K^{(d)}_{\rm Heisenberg}(x, x^{\prime})
= e^{x \cdot \overline{x^{\prime}}}, \quad x, x^{\prime} \in \C^d.
\]

The kernels in this form on $\C^d, d \in \N$ have been
studied by Zelditch and his coworkers
(see \cite{Zel00,BSZ00} and references therein),
who identified them with the Szeg\"o kernels
for the reduced Heisenberg group.
Here we call the DPPs associated with
the correlation kernels in this form
the {\it Heisenberg family of DPPs} on $\C^d, d \in \N$.
This family includes the Ginibre DPP of type A as the
lowest dimensional case with $d=1$.

\begin{df}
\label{thm:Heisenberg_class}
The Heisenberg family of DPP on $\C^d, d \in \N$
is defined by
$\Big(\Xi, K^{(d)}_{\rm Heisenberg}, \lambda_{\rN(0, 1; \C^d)}(d x) \Big)$
with
\[
K^{(d)}_{\rm Heisenberg}(x, x')
=e^{x \cdot \overline{x'}},
\quad x, x \in \C^d.
\]
\end{df}
\vskip 0.3cm
Since
\[
K^{(d)}_{\rm Heisenberg}(x, x) \lambda_{\rN(0, 1; \C^d)}(dx)
= \frac{1}{\pi^d} d x, \quad x \in \C^d,
\]
every DPP in the Heisenberg family 
is uniform on $\C^d$ and the density with respect to
the Lebesgue measure $dx$ is given by $1/\pi^d$.

\subsection{Open problems}
\label{sec:open_problems}

With $L^2(S, \lambda)$ and $L^2(\Gamma, \nu)$, 
we can consider the system of 
{\it biorthonormal functions},
which consists of a pair of distinct families of measurable functions 
$\{\psi(x, \gamma) : x \in S, \gamma \in \Gamma \}$
and $\{\varphi(x, \gamma) : x \in S, \gamma \in \Gamma \}$
satisfying the biorthonormality relations
\begin{equation}
\langle \psi(\cdot, \gamma), \varphi(\cdot, \gamma') \rangle_{L^2(S, \lambda)}
\nu(d \gamma)=\delta(\gamma-\gamma') d \gamma,
\quad \gamma, \gamma' \in \Gamma.
\label{eqn:bi_ortho}
\end{equation}
If the integral kernel defined by
\begin{equation}
K^{\rm bi}(x, x')
=\int_{\Gamma} \psi(x, \gamma) \overline{\varphi(x', \gamma)} \nu(d \gamma),
\quad x, x' \in S, 
\label{eqn:K_bi}
\end{equation}
is of finite rank, we can construct a finite DPP on $S$
whose correlation kernel is given by (\ref{eqn:K_bi})
following a standard method of random matrix theory
(see, for instance, Appendix C in \cite{Kat19a}).
By the biorthonormality (\ref{eqn:bi_ortho}),
it is easy to verify that 
$K^{\rm bi}$ is a projection kernel,
but it is not necessarily an orthogonal projection.
This observation means that such a DPP
is not constructed by the method reported in this 
manuscript.
Generalization of the present framework in order to cover
such DPPs associated with biorthonormal systems
is required.

Moreover, the dynamical extensions of DPPs
called {\it determinantal processes}
(see, for instance, \cite{Kat15_Springer})
shall be studied in the context of the present section.

\clearpage

\SSC
{One-Dimensional Stochastic Log-Gases}
\label{sec:SLG}
\subsection{Eigenvalue and singular-value processes}
\label{eigenvalue}
For $N \in \N:=\{1,2, \dots\}$, 
let $\rH(N)$ and $\U(N)$ be the space of 
$N \times N$ Hermitian matrices and
the group of $N \times N$ unitary matrices, respectively.
The space of $N \times N$ real symmetric matrices and 
the group of $N \times N$ orthogonal matrices are
denoted by $\rS(N)$ and $\gO(N)$, respectively.
As a matter of course, $\rS(N) \subset \rH(N)$
and $\gO(N) \subset \U(N)$.
In the probability space $(\Omega, \cF, \rP)$, 
we consider complex-valued processes
$M_{ij}(t) \in \C, t \geq 0, 1\leq i, j \leq N$ 
with the condition $\overline{M_{ji}(t)}=M_{ij}(t)$,
where $\overline{z}$ denotes the
complex conjugate of $z \in \C$.
We consider an $\rH(N)$-valued process
$M(t)=(M_{ij}(t))_{1 \leq i, j \leq N}$, $t \geq 0$.
Let $\R_{\geq 0} := \{x \in \R: x \geq 0\}$.
For $S=\R$ and $\R_{\geq 0}$, 
we define the Weyl chambers in $S^N$ as
\begin{equation}
\W_N(S) := \Big\{ \x=(x_1, \dots, x_N) \in S^N: x_1 < \cdots < x_N \Big\},
\label{eqn:Weyl1}
\end{equation}
and their closures as
$\overline{\W_N(S)} := \{ \x \in S^N :
x_1 \leq \cdots \leq x_N\}$.
For each $t \geq 0$, there exists 
$U(t)=(U_{ij}(t))_{1 \leq i, j \leq N} \in \U(N)$ such that
it diagonalizes $M(t)$ as
\[
U(t)^{\ast} M(t) U(t) = \Lambda(t)
:= {\rm diag}(\Lambda_1(t), \dots, \Lambda_N(t)),
\]
with the eigenvalues $\{\Lambda_i(t)\}_{i=1}^N$ 
of $M(t)$ which are assumed to be in the non-decreasing order,
\[
\Lambda_1(t) \leq \cdots \leq \Lambda_N(t)
\quad
\Longleftrightarrow
\quad
\bLambda :=(\Lambda_1(t), \dots, \Lambda_N(t))
\in \overline{\W_N(\R)}.
\]
For $d M(t):=(d M_{ij}(t))_{1 \leq i, j \leq N}$, $t \geq 0$, define
a set of quadratic variations, 
\[
\Gamma_{ij, k \ell}(t):= \langle (U^{\ast} d M  U)_{ij}, 
(U^{\ast} d M  U)_{k \ell} \rangle_t,
\quad 1 \leq i, j, k, \ell \leq N, \quad t \geq 0.
\]
We denote by $\1(\omega)$ the indicator function
of a condition $\omega$;
$\1(\omega)=1$ if $\omega$ is satisfied,
and $\1(\omega)=0$ otherwise.
In particular, given a subspace $\Gamma \subset \R^N$ we define
$\1_{\Gamma}(\x) :=\1(\x \in \Gamma)$ 
for $\x \in \R^N$, and 
$\delta_{ij} := \1(i=j)$.
The following is proved \cite{Bru89,KT04,Kat15_Springer}.
See Section 4.3 of \cite{AGZ10} for details of proof.

\begin{prop}
\label{thm:Bru}
Assume that $(M_{ij}(t))_{t \geq 0}, 1 \leq i, j \leq N$ are
continuous semi-martingales.
The eigenvalue process $(\bLambda(t))_{t \geq 0}$ satisfies
the following system of stochastic differential equations (SDEs),
\[
d \Lambda_i(t)=d \cM_i(t)+ d J_i(t), \quad t \geq 0, 
\quad 1 \leq i \leq N, 
\]
where $(\cM_i(t))_{t \geq 0}, 1 \leq i \leq N$
are martingales with quadratic variations
\[
\langle \cM_i, \cM_j \rangle_t
= \int_0^t \Gamma_{ii, jj}(s) ds, \quad t \geq 0, 
\quad 1 \leq i, j \leq N, 
\]
and $(J_i(t))_{t \geq 0}, 1 \leq i \leq N$ are the
processes with finite variations given by
\[
dJ_i(t)= \sum_{j=1}^N \frac{1}{\Lambda_i(t)-\Lambda_j(t)}
\1(\Lambda_i(t) \not= \Lambda_j(t))
\Gamma_{ij, ji} (t)dt
+ d \Upsilon_i(t), \quad t \geq 0,
\quad 1 \leq i \leq N. 
\]
Here $(d \Upsilon_i(t))_{t \geq 0}, 1 \leq i \leq N$
are the finite-variation parts of $(U(t)^{\ast} d M(t) U(t))_{ii}$.
\end{prop}

This proposition is given in \cite{KT04} as a generalization
of Theorem 1 in Bru \cite{Bru89}. 
A proof is given at Section 3.2 in \cite{Kat15_Springer}.

We will show four basic examples of $M(t) \in \rH(N), t \geq 0$
and applications of Proposition \ref{thm:Bru} \cite{KT04}.
Let $\nu \in \N_0 := \N \cup \{0\}$ and
$(B_{ij}(t))_{t \geq 0}$, $(\widetilde{B}_{ij}(t))_{t \geq 0}$, 
$1 \leq i \leq N+\nu$, $1 \leq j \leq N$
be independent one-dimensional standard Brownian motions.
For $1 \leq i \leq j \leq N$, put
\[
S_{ij}(t)= \begin{cases}
B_{ij}(t)/\sqrt{2}, \quad & (i<j),
\cr
B_{ii}(t), \quad & (i=j),
\end{cases}
\quad 
A_{ij}(t)= \begin{cases}
\widetilde{B}_{ij}(t)/\sqrt{2}, \quad & (i<j),
\cr
0, \quad & (i=j),
\end{cases}
\]
and put
$S_{ij}(t)=S_{ji}(t)$ and $A_{ij}(t)=-A_{ji}(t)$, $t \geq 0$
for $1 \leq j < i \leq N$.

\noindent
\begin{example}
\label{thm:Bru_example}
\begin{description}
\item{\rm (i)} \,
Put
$M_{ij}(t)=S_{ij}(t)+\sqrt{-1} A_{ij}(t), t \geq 0$, 
$1 \leq i, j \leq N$.
By definition
$\langle d M_{ij}, d M_{k \ell} \rangle_t=\delta_{i \ell} \delta_{j k} dt$, $t \geq 0$, 
$1 \leq i, j, k, \ell \leq N$. Hence, by unitarity of $U(t), t \geq 0$, 
we see that $\Gamma_{ij, k \ell}(t) = \delta_{i \ell, j k}$, 
which gives $\langle d \cM_i, d \cM_j \rangle_{t}=\Gamma_{ii, jj}(t) dt =\delta_{ij} dt$
and $\Gamma_{ij, ji}(t) \equiv 1$, $t \geq 0$, $1 \leq i, j \leq N$.
Then Proposition \ref{thm:Bru} proves that
the eigenvalue process $\bLambda(t)=(\Lambda_1(t), \dots, \Lambda(t))$, $t \geq 0$, 
satisfies the following system of SDEs with $\beta=2$,
\begin{equation}
d \Lambda_i(t)= d B_i(t) + \frac{\beta}{2} \sum_{\substack{1 \leq j \leq N, \cr j \not=i}}
\frac{dt}{\Lambda_i(t)-\Lambda_j(t)}, \quad t \geq 0, 
\quad 1 \leq i \leq N.
\label{eqn:Dyson1}
\end{equation}
Here $(B_i(t))_{t \geq 0}, 1 \leq i \leq N$ are
independent one-dimensional standard Brownian motions,
which are different from 
$(B_{ij}(t))_{t \geq 0}$ and $(\widetilde{B}_{ij}(t))_{t \geq 0}$ used above
to define $(S_{ij}(t))_{t \geq 0}$ 
and $(A_{ij}(t))_{t \geq 0}$, $1 \leq i, j \leq N$.
\item{\rm (ii)} \,
Put
$M(t)=(S_{ij}(t))_{1 \leq i, j \leq N} \in \rS(N), t \geq 0$.
In this case
$\langle d M_{ij}, d M_{k \ell} \rangle_t
=(\delta_{i \ell} \delta_{j k}+\delta_{ik} \delta_{j \ell}) dt/2$, 
$t \geq 0$, $1 \leq i, j, k, \ell \leq N$ and $(U(t))_{t \geq 0}$ is $\gO(N)$-valued.
Then we see that
$\Gamma_{ij, k \ell}(t)=(\delta_{i \ell} \delta_{jk}+\delta_{ik} \delta_{j \ell})/2$,
$t \geq 0, 1 \leq i, j \leq N$, 
and hence Proposition \ref{thm:Bru} proves that 
the eigenvalue process 
satisfies (\ref{eqn:Dyson1}) with $\beta=1$.
\item{\rm (iii)} \,
Consider an $(N+\nu) \times N$ rectangular-matrix-valued process 
given by $K(t)=(B_{ij}(t)+\sqrt{-1} \widetilde{B}_{ij}(t))_{1 \leq i \leq N+\nu, 1 \leq j \leq N}$,
$t \geq 0$, 
and define an $\rH(N)$-valued process by
\[
M(t)=K^\ast(t) K(t), \quad t \geq 0,
\]
where $K^\ast(t)$ is the Hermitian conjugate of $K(t)$, that is,
$K^{\ast}_{ij}(t)=\overline{K_{jk}(t)},
1 \leq i \leq N, 1 \leq j \leq N+\nu$.
The matrix $M$ is positive definite and hence the eigenvalues are non-negative;
$\Lambda_i(t) \in \R_{\geq 0}$, 
$t \geq 0, 1 \leq i \leq N$. 
We see that the finite-variation part of $dM_{ij}(t)$ is 
equal to $2(N+\nu) \delta_{ij} dt$, $t \geq 0$, 
and 
$\langle d M_{ij}, d M_{k \ell} \rangle_t=2(M_{i \ell}(t) 
\delta_{jk}+M_{k \ell}(t) \delta_{i \ell}) dt$,
$t \geq 0$, $1 \leq i, j, k, \ell \leq N$,
which implies that
$d \Upsilon_i(t)=2(N+\nu) dt$, 
$\Gamma_{ij, ji}=2(\Lambda_i(t)+\Lambda_j(t))$, 
and $\langle d \cM_i, d \cM_j \rangle_t = \Gamma_{ii, jj}(t) dt 
= 4 \Lambda_i(t) \delta_{ij} dt$,
$t \geq 0$, $1 \leq i, j \leq N$.
Then we have the SDEs for eigenvalue processes,
\begin{align}
d \Lambda_i(t)
&=2 \sqrt{\Lambda_i(t)} d \widetilde{B}_i(t) 
+ \beta \left[ (\nu+N) +
\sum_{\substack{1 \leq j \leq N, \cr j \not=i}}
\frac{\Lambda_i(t)+\Lambda_j(t)}{\Lambda_i(t)-\Lambda_j(t)} \right] dt
\nonumber\\
&= 2 \sqrt{\Lambda_i(t)} d \widetilde{B}_i(t) 
+ \beta \left[ (\nu+1) + 2 \Lambda_i(t)
\sum_{\substack{1 \leq j \leq N, \cr j \not=i}}
\frac{1}{\Lambda_i(t)-\Lambda_j(t)} \right] dt, 
\nonumber\\
& \hskip 5cm t \geq 0, \quad 1 \leq i \leq N,
\label{eqn:Wishart1}
\end{align}
with $\beta=2$,
where $(\widetilde{B}_i(t))_{t \geq 0}, 1 \leq i \leq N$ are
independent one-dimensional standard Brownian motions,
which are different from $(B_{ij}(t))_{t \geq 0}$ 
and $(\widetilde{B}_{ij}(t))_{t \geq 0}$, 
$1 \leq i, j \leq N$, used above
to define the rectangular-matrix-valued process
$(K(t))_{t \geq 0}$.
The positive roots of eigenvalues of $M(t)$ give the
{\it singular values} of the rectangular matrix $K(t)$,
which are denoted by 
$\cS_i(t) =\sqrt{\Lambda_i(t)}, t \geq 0, 1 \leq i \leq N$.
The system of SDEs for them is readily obtained
from (\ref{eqn:Wishart1}) as
\begin{align}
d \cS_i(t) &= d \widetilde{B}_i(t)+ \frac{1}{2} \left[\frac{\beta(\nu+1)-1}{\cS_i(t)}
+ \beta \sum_{\substack{1 \leq j \leq N, \cr j \not=i}}
\left( \frac{1}{\cS_i(t)-\cS_j(t)} + \frac{1}{\cS_i(t)+\cS_j(t)} \right) \right] dt,
\nonumber\\
& \hskip 5cm t \geq 0, \quad 1 \leq i \leq N,
\label{eqn:Wishart2}
\end{align}
with $\beta=2$ and $\nu \in \N_0$. 
\item{\rm (iv)} \,
Put $K(t)=(B_{ij}(t))_{1 \leq i \leq N+\nu, 1 \leq j \leq N}$, $t \geq 0$
and consider the process in $\rS(N)$, 
$M(t)=K(t)^{\ast} K(t)$, $t \geq 0$.
In this case the finite-variation part of $dM_{ij}(t)$ is $(N+\nu) \delta_{ij} dt$,
$t \geq 0$, and 
$\langle d M_{ij}, d M_{k \ell} \rangle_t=(M_{ik}(t) \delta_{j \ell}+M_{i \ell}(t) \delta_{jk}
+M_{jk}(t) \delta_{i \ell}+M_{j \ell}(t) \delta_{ik}) dt$, $t \geq 0$, 
$1 \leq i, j, k, \ell \leq N$,
which imply that
$d \Upsilon_i(t)=(N+\nu) dt$, 
$\Gamma_{ij, ji}=(\Lambda_i(t)+\Lambda_j(t))(1+\delta_{ij})$, 
and $\langle d \cM_i, d \cM_j \rangle_t = \Gamma_{ii, jj}(t)
=4 \Lambda_i(t) \delta_{ij} dt$, $t \geq 0$, 
$1 \leq i, j \leq N$. 
Then we have (\ref{eqn:Wishart1}) with $\beta=1, \nu \in \N_0$
as the SDEs for the eigenvalue process of $(M(t))_{t \geq 0}$,
and (\ref{eqn:Wishart2}) with $\beta=1, \nu \in \N_0$ as
the SDEs for the singular-value process of $(K(t))_{t \geq 0}$.
\end{description}
\end{example}
\vskip 0.3cm

Other examples of $M(t) \in \rH(N), t \geq 0$ are
shown in \cite{KT04}, in which the eigenvalue processes
following the SDEs (\ref{eqn:Dyson1}) with $\beta=4$ are
also shown.

Dyson \cite{Dys62} introduced the system of SDEs similar to (\ref{eqn:Dyson1})
with $\beta=1,2$, and 4, but with drift terms of
the Ornstein--Uhlenbeck type so that the stationary measures
of the processes are given by the eigenvalue distributions
of the {\it Gaussian orthogonal ensemble} (GOE),
the {\it Gaussian unitary ensemble} (GUE), and
the {\it Gaussian symplectic ensemble} (GSE) 
studied in random matrix theory 
\cite{Meh04,For10,AGZ10,ABD11,Tao12}
for $\beta=1,2$, and $\beta=4$, respectively.
For general $\beta \in \R$, the system of SDEs (\ref{eqn:Dyson1})
starting from the configuration 
$\blambda :=\bLambda(0) \in \W_N(\R)$ can be
defined up to the first collision time,
\[
T^{\blambda} := \inf \{ t > 0 : \Lambda_i(t) = \Lambda_j(t) \quad
\mbox{for some $1 \leq i \not= j \leq N$} \}.
\]
It was proved that when $\beta \geq 1$,
$T^{\blambda} = \infty$ a.s. and 
(\ref{eqn:Dyson1}) has a strong and pathwise unique 
non-colliding solution for $\blambda \in \W_N(\R)$
\cite{RS93}. The statement was extended
to general initial configuration $\blambda \in \overline{\W_N(\R)}$
\cite{CL97,Dem09,GM13,GM14}. 
The transition probability density of the process
(\ref{eqn:Dyson1}) with $\beta >0$ was studied
by Baker--Forrester \cite{BF97} and 
R\"{o}sler \cite{Ros98}
using the multivariate special functions (see also \cite{Dem08}).

The $\rS(N)$-valued process given as
Example \ref{thm:Bru_example} (iv) above was first studied by
Bru \cite{Bru91}, in which the SDEs (\ref{eqn:Wishart1}) with
$\beta=1$ was considered.
When the initial configuration $\blambda$ is specially chosen as
the null vector $\0 :=(0, \dots, 0)$, 
the probability density function (PDF) at time $t=1$
realizes the eigenvalue distribution called
the {\it real Wishart ensemble} \cite{Wis28} 
or the {\it chiral GOE} studied in random matrix theory \cite{KT04,ABD11}.
K\"onig and O'Connell \cite{KO01} 
studied the $\rH(N)$-valued process of Example \ref{thm:Bru_example} (iii) 
and called it
the {\it Laguerre process}.
Graczyk and Ma{\l}ecki \cite{GM13,GM14} proved that,  
if $\beta \geq 1$, 
(\ref{eqn:Wishart1}) has a strong and pathwise unique 
non-colliding solution for general initial configurations
$\blambda \in \overline{\W_N(\R_{\geq 0})}$.

When $\beta=2$, the eigenvalue processes 
(\ref{eqn:Dyson1}), (\ref{eqn:Wishart1}), and 
the singular-value process (\ref{eqn:Wishart2})
are realized in $\R$ as the $N$-tuples of the one-dimensional
standard Brownian motions, the $2(\nu+1)$-dimensional
squared Bessel processes, 
and the $2(\nu+1)$-dimensional Bessel processes 
{\it conditioned never to collide with each other} 
\cite{Gra99,KO01,KT04,KT07,KT11b}.
Here $\nu \in (-1, \infty)$, and if $\nu \in (-1, 0]$ 
the reflection boundary condition will be assumed at the origin
for (\ref{eqn:Wishart2}). 
These noncolliding diffusion processes are proved to be 
{\it determinantal stochastic processes} \cite{BR05,KT07b} 
and all spatio-temporal correlation functions
are explicitly expressed by determinants.
The determinants are governed by the functions
called the spatio-temporal correlation kernels, which can be
simply expressed using the Hermite and
Laguerre polynomials if the initial configurations
$\blambda$ are given by $\0$.
We will show these facts in
Sections \ref{sec:DSPs}--\ref{sec:multi}.
For more detail, 
see \cite{KT10,KT11,Kat14,Kat15_Springer}. 

\subsection{Stochastic log-gases and 2D-Coulomb gases 
confined in one-dimensional spaces}
\label{sec:log_gases2}

In the last section, we will consider the
Schramm--Loewner evolution (SLE).
Schramm used a parameter $\kappa >0$ in order to parametrize
time change of the Brownian motion \cite{Sch00}.
Accordingly, we change the parameter $\beta \to \kappa$
by setting
\begin{equation}
\beta=\frac{8}{\kappa},
\label{eqn:kappa1}
\end{equation}
and perform the time change $t \to \kappa t$.
Since 
$(B(\kappa t))_{t \geq 0} \law= (\sqrt{\kappa} B(t))_{t \geq 0}$,
if we put
$Y^{\R}_i(t) := \Lambda_i(\kappa t)$,
$Y^{\R_{\geq 0}}_i(t) := \cS_i(\kappa t)$, $t \geq 0$, $1 \leq i \leq N$,
the systems of SDEs (\ref{eqn:Dyson1}) and (\ref{eqn:Wishart1})
are written as
\begin{align}
d Y^{\R}_i(t) &= \sqrt{\kappa} dB_i(t)
+ 4 \sum_{\substack{1 \leq j \leq N, \cr j \not=i}}
\frac{dt}{Y^{\R}_i(t)-Y^{\R}_j(t)}, \quad
t \geq 0, \quad 1 \leq i \leq N,
\label{eqn:DysonB}
\\
d Y^{\R_{\geq 0}}_i(t) &= \sqrt{\kappa} d \widetilde{B}_i(t)
+ 4 \sum_{\substack{1 \leq j \leq N, \cr j \not=i}}
\left( \frac{1}{Y^{\R_{\geq 0}}_i(t)-Y^{\R_{\geq 0}}_j(t)}
+ \frac{1}{Y^{\R_{\geq 0}}_i(t)+Y^{\R_{\geq 0}}_j(t)} \right) dt
\nonumber\\
& \qquad \qquad \quad 
+\left\{ 4 (1+\nu) - \frac{\kappa}{2} \right\}
\frac{dt}{Y^{\R_{\geq 0}}_i(t)}, 
\quad
t \geq 0, \quad 1 \leq i \leq N.
\label{eqn:Bru_WishartB}
\end{align}
In the last section, we will call
$\Y^{\R}(t) =(Y^{\R}_1(t), \dots, Y^{\R}_N(t)) \in \R^N, t \geq 0$,
the $(8/\kappa)$-{\it Dyson model} and
$\Y^{\R_{\geq 0}}(t) =(Y^{\R_{\geq 0}}_1(t), \dots, Y^{\R_{\geq 0}}_N(t)) 
\in (\R_{\geq 0})^N, t \geq 0$,
the $(8/\kappa, \nu)$-{\it Bru--Wishart process},
respectively.

The above stochastic processes
$\Y^{S}(t) =(Y^S_1(t), \dots, Y^S_N(t)), t \geq 0, S=\R$ or $\R_{\geq 0}$
can be written in the form,
\begin{equation}
dY^S_i(t)=\sqrt{\kappa} dB_i(t)
+ \left. 
\frac{\partial \phi^S(\x)}{\partial x_i} 
\right|_{\x=\Y^S(t)} dt,
\quad t \geq 0, \quad 1 \leq i \leq N,
\label{eqn:log_gas1}
\end{equation}
if we introduce the following potential energies,
\begin{equation}
\phi^{S}(\x) :=
\begin{cases}
\displaystyle{
4 \sum_{1 \leq i < j \leq N} 
\log(x_j-x_i)}, 
\quad \mbox{for $S=\R$},
&
\cr
\displaystyle{
4 \sum_{1 \leq i < j \leq N}
\Big[ \log(x_j-x_i)+\log(x_j+x_i) \Big]
+ \left\{ 4 (\nu+1) - \frac{\kappa}{2} \right\} \sum_{i=1}^N \log x_i, 
} &
\cr
\hskip 4.5cm \mbox{for $S=\R_{\geq 0}$}.
&
\end{cases}
\label{eqn:potential}
\end{equation}
In both cases of (\ref{eqn:potential}),
the potentials
are given by logarithmic functions,
and  the drift terms are gradient forces of these potentials.
(See Remark \ref{thm:log1} in Section \ref{sec:classical}.)
In this sense 
the $(8/\kappa)$-Dyson model
and the $(8/\kappa, \nu)$-Bru--Wishart process 
are regarded as {\it stochastic log-gases} in $\R$ \cite{For10}.
Since the logarithmic potential describes
the two-dimensional Coulomb law in electrostatics,
the present processes are also considered
as stochastic models of two-dimensional
charged-particles ({\it $2D$-Coulomb gas}) confined on 
a line $\R$ or to a half-line $\R_{\geq 0}$.

\subsection{Determinantal martingales and
determinantal stochastic processes (DSPs)}
\label{sec:DSPs}
We consider the same setting as in Section \ref{sec:DPP}.
Let $S$ be a base space, which is locally compact Hausdorff space
with countable base, 
and $\lambda$ be a Radon measure on $S$.
The configuration space over $S$ is given by
the set of nonnegative-integer-valued Radon measures; 
\[
\Conf(S)
=\left\{ \xi = \sum_j \delta_{x_j} : \mbox{$x_j \in S$,
$\xi(\Lambda) < \infty$ for all bounded set $\Lambda \subset S$} \right\}.
\]
${\rm Conf}(S)$ is equipped with the topological Borel $\sigma$-fields
with respect to the vague topology.

We consider {\it interacting particle systems} 
as $\Conf(S)$-valued continuous-time processes and write them as
\begin{equation}
\Xi(t)=\sum_{j=1}^N \delta_{X_j(t)},
\quad t \geq 0,
\label{eqn:Xi1}
\end{equation}
where $\X(t)=(X_1(t), \dots, X_N(t)), t \geq 0$ 
are defined by a solution of
a given SDEs.
We call $\x \in S^N$ a {\it labeled configuration}
and $\xi \in \Conf(S)$ an {\it unlabeled configuration}.
The probability law of $(\Xi(t))_{t \geq 0}$
starting from a fixed configuration $\xi \in \Conf(S)$
is denoted by $\bP_{\xi}$ 
and the process
specified by the initial configuration
is expressed by
$((\Xi(t))_{t \geq 0}, \bP_{\xi})$.
The expectations with respect to $\bP_{\xi}$
is denoted by $\bE_{\xi}$.
We introduce a filtration $(\cF_{\Xi}(t))_{t \geq 0}$
generated by $(\Xi(t))_{t \geq 0}$ satisfying the usual conditions
(see, for instance, p.45 in  \cite{RY05}).
We set 
\[
\Conf_{0}(S)= \{ \xi \in \Conf(S) : 
\xi(\{x\}) \leq 1 \mbox { for any }  x \in S \},
\]
which gives a collection of configurations
of simple point processes (i.e., 
without multiple points).

Let $0 \leq T < \infty$.
We consider the expectation
of an $\cF_{\Xi}(T)$-measurable bounded function $F$,
$\bE_{\xi} [F(\Xi(\cdot))]$.
It is sufficient 
to consider the case that $F$ is given as
\begin{equation}
F\left(\Xi(\cdot)\right)
= \prod_{m=1}^M g_m(\X(t_m))
\label{eqn:F_gm}
\end{equation}
for an arbitrary $M \in \N$, $0 \leq t_1< \cdots <t_M \leq T < \infty$
with bounded measurable functions $g_m$ 
on $S^N$, $1 \leq m \leq M$.
Since the particles are unlabeled 
in the process $((\Xi(t))_{t \geq 0}, \bP_{\xi})$,
$g_m$'s are {\it symmetric functions}.

We consider a continuous-time Markov process $(Y(t))_{t \geq 0}$ 
on $S$
which is a connected open set in $\R$.
It is a diffusion process in $S$ or a process showing
a position on the circumference 
$S=[0, 2 \pi r)$ of a diffusion process moving around
on the circle with a radius $r>0$;
$\S^{1}(r) := \{x \in [0, 2 \pi r) : x+2 \pi r = x\}$.
The probability space is denoted by $(\Omega, \cF, \rP^v)$
with expectation $\rE^v$, 
when the initial state is fixed to be $v \in S$.
When $v$ is the origin, the subscript is omitted.
We introduce a filtration $\{\cF(t) : t \geq 0 \}$
generated by $Y$ so that it satisfies the usual conditions
(see, for instance, p.45 in  \cite{RY05}).
We assume that the process has
a {\it transition density},
$p(t, y|x)$, 
$t \in [0, \infty), x, y \in S$
such that for any measurable bounded 
function $f(t, x), (t, x) \in [0, \infty) \times S$,
\begin{equation}
\rE[f(t, Y(t)) | \cF(s)]
= \int_{S} dy \, f(t, y) p(t-s, y|Y(s)) \quad
\mbox{a.s., \, $0 \leq s \leq t < \infty$}.
\label{eqn:def_p}
\end{equation}

Recall that $\W_N(S), N \in \N$ denotes the Weyl chamber in $S$;
\[
\W_N(S)=\{\x=(x_1, \dots, x_N) \in S^{N} : x_1 < x_2 < \cdots < x_N \}.
\]
Given $\u=(u_1, \dots, u_N) \in \W_N(S)$, 
we have a measure 
$\xi(\cdot)=\sum_{j=1}^N \delta_{u_j}(\cdot) \in \Conf_0(S)$.
Depending on $\xi$, we 
assume that there is a one-parameter
family of continuous functions
\[
\cM_{\xi}^v(\cdot, \cdot) : 
[0, \infty) \times S \to \R
\]
with parameter $v \in S$, 
such that the processes 
$\cM_{\xi}(t, Y(t))=\{\cM_{\xi}^v(t, Y(t)) : v \in \{u_1, \dots, u_N\} \}$, 
$t \geq 0$
satisfy the following conditions.
\begin{description}
\item{\bf (M1)} \quad
$(\cM_{\xi}^{u_k}(t, Y(t)))_{t \geq 0}$, $1 \leq k \leq N$ are continuous martingales; 
\[
\rE[\cM_{\xi}^{u_k}(t, Y(t)) | \cF(s)]=\cM_{\xi}^{u_k}(s, Y(s)) \quad 
\mbox{a.s. for all $0 \leq s \leq t < \infty$.}
\]
\item{\bf (M2)} \quad
For any time $t \geq 0$, $\cM_{\xi}^{u_k}(t, x), 1 \leq k \leq N$
are linearly independent functions of $x$.
\item{\bf (M3)} \quad
For $1 \leq j, k \leq N$,
\[
\lim_{t \to 0} \rE^{u_j}[\cM_{\xi}^{u_k}(t, Y(t))] 
= \delta_{jk}.
\]
\end{description}
We call $\{\cM^v_{\xi}(t, x) : v \in \{u_1, \dots, u_N\} \}$ {\it martingale functions}.

Let $(Y_j(t))_{t \geq 0}, 1 \leq j \leq N$
be a collection of $N$ independent copies 
of $(Y(t))_{t \geq 0}$.
We consider the $N$-component vector-valued Markov process
$\Y(t)=(Y_1(t), \dots, Y_N(t))$, $t \geq 0$,
for which the initial values are fixed to be
$Y_j(0)=u_j \in S, 1 \leq j \leq N$,
provided $\u=(u_1, \dots, u_N) \in \W_N(S)$.
We consider a determinant of the martingales
\begin{equation}
\cD_{\xi}(t, \Y(t)) 
=\det_{1 \leq j,k \leq N} \Big[ 
\cM_{\xi}^{u_k}(t, Y_j(t)) \Big],
\quad t \geq 0. 
\label{eqn:D1}
\end{equation}
The condition {\bf (M2)} is necessary so that 
it is not zero constantly.
This determinant is a continuous martingale 
and we call it a {\it determinantal martingale}.

\begin{df}
\label{thm:DM_rep}
Given $\xi \in \Conf_0(S)$, consider a process $((\Xi(t))_{t \geq 0}, \bP_{\xi})$. 
If there exists a pair $(Y, \cM_{\xi})$
defining $\cD_{\xi}$ by (\ref{eqn:D1}) such that
for any $\cF_{\Xi}(t)$-measurable bounded function $F$,
$0 \leq t \leq T < \infty$, 
the equality 
\[
\bE_{\xi}[F(\Xi(\cdot))]
=\rE^{\u} \left[
F \left( \sum_{j=1}^N \delta_{Y_j(\cdot)} \right) 
\cD_{\xi}(T, \Y(T)) \right]
\]
holds, 
then we say $((\Xi(t))_{t \in \geq 0}, \bP_{\xi})$ has 
a determinantal-martingale
representation (DMR)
associated with $(Y, \cM_{\xi})$.
\end{df}

Consider an arbitrary but fixed $M \in \N$.
Assume that we have 
a sequence of times
$\t=(t_1,\dots,t_M)$ with 
$0 \leq t_1 < \cdots < t_M < \infty$, 
and a sequence of functions
$\f=(f_{t_1},\dots,f_{t_M}) \in \cC_{\rm c}(S)^M$.
Then the {\it multitime Laplace transform of 
probability measure} $\bP_{\xi}$ is defined by
\begin{equation}
\Psi_{\xi}^{\t}[\f]
:= \bE_{\xi} \left[ \exp \left\{ \sum_{m=1}^{M} 
\int_{S} f_{t_m}(x) \Xi(t_m, dx) \right\} \right].
\label{eqn:GF1}
\end{equation}
We assume that 
it is expanded with respect to 
$\{1-e^{f_{t_m}(\cdot)}: 1 \leq m \leq M\}$ as
\begin{align}
\Psi_{\xi}^{\t}[\f]
&=1+ \prod_{m=1}^M \sum_{1 \leq n_m \leq N}
\frac{(-1)^{n_m}}{n_m !}
\int_{S^{n_m}} \lambda_{t_m}^{\otimes n_m} (d \x^{(m)}_{n_m})
\prod_{\ell=1}^{n_m} (1- e^{f_{t_m}(x^{(m)}_{\ell})}) 
\nonumber\\
& \hskip 6cm \times
\rho_{\xi} 
\Big( t_{1}, \x^{(1)}_{n_1}; \dots ; t_{M}, \x^{(M)}_{n_M} \Big),
\label{eqn:LaplaceB}
\end{align}
where $(\lambda_t)_{t \geq 0}$ is a time-dependent
measure on $S$, 
$\x^{(m)}_{n_m}$ denotes
$(x^{(m)}_1, \dots, x^{(m)}_{n_m})$, $1 \leq m \leq M$, 
$1 \leq n_m \leq N$, and
$d \x^{(m)}_{n_m}= \prod_{j=1}^{n_m} dx^{(m)}_j$,
$1 \leq m \leq M$. 
This expansion formula of $\Psi_{\xi}^{\t}[\f]$
defines the
{\it spatio-temporal correlation functions}
$\rho_{\xi}$ for the process $((\Xi(t))_{t \geq 0}, \bP_{\xi})$.

\begin{df}
\label{thm:determinantal}
A process $((\Xi(t))_{t \geq 0}, \bP_{\xi})$ is said to be 
a determinantal stochastic process (DSP)
with spatio-temporal correlation kernel 
$\bK_{\xi} : ([0, \infty) \times S)^2 \to \R$, 
if any multitime Laplace transform
of $\bP_{\xi}$ (\ref{eqn:GF1}) can be expanded 
in the form (\ref{eqn:LaplaceB})
with a time-dependent measure $(\lambda_t)_{t \geq 0}$ on $S$
and 
all spatio-temporal correlation functions
are given by determinants as 
\begin{equation}
\rho_{\xi} \Big(t_1,\x^{(1)}_{n_1}; \dots;t_M,\x^{(M)}_{n_M} \Big) 
=\det_{\substack
{1 \leq j \leq n_{m}, 1 \leq k \leq n_{n}, \\
1 \leq m, n \leq M}
}
\Bigg[
\bK_{\xi}(t_m, x_{j}^{(m)}; t_n, x_{k}^{(n)} )
\Bigg],
\label{eqn:rho1}
\end{equation}
$0 \leq t_1 < \cdots < t_M < \infty$, 
$1 \leq n_m \leq N$,
$\x^{(m)}_{n_m} \in S^{n_m}, 1 \leq m \leq M \in \N$.
The DSP is denoted by
$((\Xi(t))_{t \geq 0}, \bK_{\xi}, (\lambda_t)_{t \geq 0})$.
\end{df}
By Definition \ref{thm:determinantal},
the multitime Laplace transform of $\bP_{\xi}$ 
(\ref{eqn:GF1}) is written as follows.
We regard this as a definition of the 
{\it multitime Fredholm determinant} of an integral 
kernel $\sK_{\xi}$ specified by $\bK_{\xi}$,
\begin{align}
&1+ \prod_{m=1}^M \sum_{1 \leq n_m \leq N}
\frac{(-1)^{n_m}}{n_m !}
\int_{S^{n_m}} \lambda_{t_m}^{\otimes n_m} (d \x^{(m)}_{n_m})
\prod_{\ell=1}^{n_m} (1- e^{f_{t_m}(x^{(m)}_{\ell})}) 
\nonumber\\
& \hskip 5cm \times
\det_{\substack
{1 \leq j \leq n_{m}, 1 \leq k \leq n_{n}, \\
1 \leq m, n \leq M}
}
\Bigg[
\bK_{\xi}(t_m, x_{j}^{(m)}; t_n, x_{k}^{(n)} )
\Bigg]
\nonumber\\
&\qquad =:  \mathop{{\rm Det}}_
{\bigotimes_{m=1}^M L^2(S, \lambda_{t_m})}
\Big[I-(1-e^{f_{\cdot}}) \sK_{\xi} \Big].
\label{eqn:Fredholm1}
\end{align}

The main theorem in this section is the following.
\begin{thm}
\label{thm:DM_det}
If $((\Xi(t))_{t \geq 0}, \bP_{\xi})$ has DMR associated 
with $(Y, \cM_{\xi})$, 
then it is a DSP 
$((\Xi(t))_{t \geq 0}, \bK_{\xi}, dx)$
with the spatio-temporal correlation kernel
\begin{equation}
\bK_{\xi}(s,x;t,y)
=\int_{S} \xi(dv) p(s, x|v) \cM_{\xi}^{v}(t,y)
- \1(s>t) p(s-t,x|y),
\label{eqn:KK1}
\end{equation}
$(s,x), (t,y) \in [0, \infty) \times S$,
where $p$ is the transition density of the process $Y$.
\end{thm}

This type of correlation kernel (\ref{eqn:KK1}) 
was first obtained
by Eynard and Mehta 
for a multi-matrix model \cite{EM98}
and by Nagao and Forrester \cite{NF98}
for the noncolliding Brownian motion started at a special
initial distribution 
$\bp^{(N)}_{\rm Hermite}$ given by (\ref{eqn:P_Laguerre})
(the GUE eigenvalue distribution), 
and has been extensively studied
\cite{FNH99,TW04,Meh04,BR05,For10,KT10,KT13}.
Note that, while all correlation kernels
obtained discussed in Section \ref{sec:DPP} are symmetric, 
the present spatio-temporal correlation kernel 
is asymmetric with respect to
the exchange of two points $(s, x)$ and $(t, y)$
on the spatio-temporal plane $[0, \infty) \times S$ 
and shows {\it causality} in the system.

In the following, first we states a lemma and a proposition, and then
prove Theorem \ref{thm:DM_det} using them.
For $n \in \N$, an index set $\{1,2, \dots, n\}$
is denoted by $\I_{n}$.
Fixing $N \in \N$ with $N' \in \I_N$, 
we write
\[
\J \subset \I_N, \sharp \J=N'
\quad
\Longleftrightarrow 
\quad
\J=\{j_1, \dots, j_{N'}\}, \quad
1 \leq j_1 < \dots < j_{N'} \leq N.
\]
For  $\x=(x_1, \dots, x_N) \in \R^N$, 
put $\x_{\J} : =(x_{j_1}, \dots, x_{j_{N'}})$.
In particular, we write
$\x_{N'} :=\x_{\I_{N'}}, 1 \leq N' \leq N$.
(By definition $\x_N=\x$.)
A collection of all permutations of 
elements in $\J$ is denoted by $\mS(\J)$.
In particular, we write $\mS_{N'} :=\mS(\I_{N'}), 1 \leq N' \leq N$.

The following lemma shows the
{\it reducibility} of the determinantal martingale
in the sense that,
if we observe a symmetric function depending
on $N'$ variables, $N' \leq N$,
then the size of determinantal
martingale can be reduced from $N$ to $N'$.

\begin{lem}
\label{thm:reduce}
Let $\u=(u_1, \dots, u_N) \in \W_N(S)$ and
$\xi=\sum_{j=1}^N \delta_{u_j} \in \Conf_0(S)$.
Assume that there exists a pair $(Y, \cM_{\xi})$
satisfying conditions (M1)--(M3) and $\cD_{\xi}$
is defined by (\ref{eqn:D1}).
Let $1 \leq N' \leq N$.
For $0 < t \leq T < \infty$ and an 
$\cF_{\Xi}(t)$-measurable symmetric function
$F_{N'}$ on $\R^{N'}$,
\begin{align}
& \sum_{\J \subset \I_N, \sharp \J=N'}
\rE^{\x} \left[
F_{N'}(\Y_{\J}(t))
\cD_{\xi}(T, \Y(T)) \right]
\nonumber\\
& \quad
= \int_{\W_{N'}^{\rm A}} \xi^{\otimes N'} (d \bv)
\rE^{\bv} \left[
F_{N'}(\Y_{N'}(t))
\cD_{\xi}(T, \Y_{N'}(T)) \right].
\label{eqn:reducibility}
\end{align}
\end{lem}
\noindent{\it Proof} \quad
By the definition (\ref{eqn:D1}), 
LHS of (\ref{eqn:reducibility}) is equal to
\begin{align}
& \sum_{\J \subset \I_N, \sharp \J=N'}
\rE^{\u} \left[ F_{N'}(\Y_{\J}(t))
\det_{i, j \in \I_N}
[\cM_{\xi}^{x_j}(T, Y_i(T))] \right]
\nonumber\\
& \, = \sum_{\J \subset \I_N, \sharp \J=N'}
\rE^{\u} \left[ F_{N'}(\Y_{\J}(t))
\sum_{\sigma \in \mS_N} {\rm sgn}(\sigma)
\prod_{i=1}^N \cM_{\xi}^{u_{\sigma(i)}}(T, Y_{i}(T)) \right]
\nonumber\\
& \, = \sum_{\J \subset \I_N, \sharp \J=N'}
\sum_{\sigma \in \mS_N} {\rm sgn}(\sigma)
\nonumber\\
& \qquad \times \ 
\rE^{\u} \left[ F_{N'}(\Y_{\J}(t))
\prod_{i \in \J} \cM_{\xi}^{u_{\sigma(i)}}(T, Y_i(T))
\prod_{j \in \I_N \setminus \J}
\cM_{\xi}^{u_{\sigma(j)}}(T, Y_j(T)) \right]
\nonumber
\end{align}
Since $(Y_i(t))_{t \geq 0}, 1 \leq i \leq N$ are independent, it is equal to
\begin{align}
&  \sum_{\J \subset \I_N, \sharp \J=N'}
\sum_{\sigma \in \mS_N} {\rm sgn}(\sigma)
\rE^{u} \left[ F_{N'}(\Y_{\J}(t))
\prod_{i \in \J} \cM_{\xi}^{u_{\sigma(i)}}(T, Y_i(T)) \right]
\nonumber\\
& \qquad \qquad \qquad \qquad \qquad \times
\prod_{j \in \I_N \setminus \J}
\rE^{\u} \left[ \cM_{\xi}^{u_{\sigma(j)}}(T, Y_j(T)) \right].
\label{eqn:A1}
\end{align}
By the condition {\bf (M1)} of $\cM_{\xi}$, 
\[
\prod_{j \in \I_N \setminus \J}
\rE^{\u} \left[ \cM_{\xi}^{u_{\sigma(j)}}(T, Y_j(T)) \right]
=
\prod_{j \in \I_N \setminus \J}
\rE^{u} \left[ \cM_{\xi}^{u_{\sigma(j)}}(t, Y_j(t)) \right],
\quad ^{\forall} t \in [0, T], 
\]
and by the condition {\bf (M3)} of $\cM_{\xi}$, 
this is equal to 
$\prod_{j \in \I_N \setminus \J}
\delta_{j \sigma(j)}$.
Then (\ref{eqn:A1}) becomes
\begin{align*}
&  \sum_{\J \subset \I_N, \sharp \J=N'}
\sum_{\sigma \in \mS(\J)} {\rm sgn}(\sigma)
\rE^{\u} \left[ F_{N'}(\Y_{\J}(t))
\prod_{i \in \J} \cM_{\xi}^{u_{\sigma(i)}}(T, Y_i(T)) \right]
\nonumber\\
& \quad =
\sum_{\J \subset \I_N, \sharp \J=N'}
\rE^{\u} \left[ F_{N'}(\Y_{\J}(t))
\det_{i, j \in \J} [ \cM_{\xi}^{u_j}(T, Y_i(T))] \right]
\nonumber\\
& \quad =
\int_{\W_{N'}^{\rm A}} \xi^{\otimes N'}(d \bv)
\rE^{\bv} \left[ F_{N'}(\Y_{N'}(t))
\det_{i, j \in \I_{N'}} [ \cM_{\xi}^{v_j}(T, Y_i(T))] \right],
\end{align*}
where equivalence in probability law 
of $(Y_i(t))_{t \geq 0}, 1 \leq i \leq N$ 
is used.
This is RHS of (\ref{eqn:reducibility}) and 
the proof is completed. 
\qed
\vskip 0.3cm

\begin{prop}
\label{thm:Fredholm}
Let $\u=(u_1, \dots, u_N) \in \W_N(S)$ and
$\xi=\sum_{j=1}^N \delta_{u_j} \in \Conf_0(S)$.
Assume that there exists a pair $(Y, \cM_{\xi})$
satisfying conditions (M1)--(M3) and $\cD_{\xi}$
is defined by (\ref{eqn:D1}).
Then for any $M \in \N$, $0 \leq t_1 < \cdots < t_M \leq T < \infty$,
$f_{t_m} \in \cC_{\rm c}(S), 1 \leq m \leq M$,
the equality
\begin{align}
& \rE^{\u} \left[
\prod_{m=1}^M \prod_{j=1}^N
\Big\{1-(1-e^{f_{t_m}(Y_j(t_m))}) \Big\}
\cD_{\xi}(T, \Y(T)) \right]
=\mathop{{\rm Det}}_
{\bigotimes_{m=1}^M L^2(S, \lambda_{t_m})}
\Big[I-(1-e^{f_{\cdot}}) \sK_{\xi} \Big]
\label{eqn:Fredholm2}
\end{align}
holds, 
where RHS
is the multitime Fredholm determinant of
$\sK_{\xi}$ specified by the
spatio-temporal correlation kernel 
(\ref{eqn:KK1}).
\end{prop}

Let $\chi_{t}(\cdot) := 1-e^{f_t(\cdot)}$, $t \geq 0$.
LHS of (\ref{eqn:Fredholm2})
is an expectation of a usual determinant multiplied by
$\prod_{m=1}^M \prod_{j=1}^N (1-\chi_{t_m}(Y_j(t_m)) )$, 
while RHS is a multitime Fredholm determinant.
First we expand LHS with respect to
$\{\chi_{t_m}(Y_j(t_m)) : 1 \leq m \leq M, 1 \leq j \leq N\}$
and apply Lemma \ref{thm:reduce}. 
The expectation of each term in LHS will be calculated
by performing integrals using 
the transition density $p$ of the process $Y$
as an integral kernel,
while $p$ is involved in the integral representation
(\ref{eqn:KK1}) of the spatio-temporal correlation kernel $\bK_{\xi}$
for the multitime Fredholm determinant in RHS.
Therefore, simply to say, this equality is just obtained
by applying Fubini's theorem.
Since the quantities in (\ref{eqn:Fredholm2}) are multivariate
and multitime joint distribution is considered,
however, we also need combinatorics arguments
to prove Proposition \ref{thm:Fredholm}.
The proof was given in \cite{KT13,Kat14}.
Here we omit the proof of this lemma.

\vskip 0.3cm
\noindent{\it Proof of Theorem \ref{thm:DM_det}} \,
By (\ref{eqn:Xi1}),
the multitime Laplace transform of $\bP_{\xi}$ (\ref{eqn:GF1}) is 
written in the form
\[
\Psi_{\xi}^{\t}[\f]
=\bE_{\xi} \left[ \prod_{m=1}^N \prod_{j=1}^N
\Big\{1-(1-e^{f_{t_m}(X_j(t_m))}) \Big\} \right].
\]
By assumption of the theorem, it has DMR associated with $(Y,\cM_{\xi})$,
\[
\Psi_{\xi}^{\t}[\f]
=\rE^{\u} \left[
\prod_{m=1}^N \prod_{j=1}^N
\Big\{1-(1-e^{f_{t_m}(Y_j(t_m))}) \Big\}
\cD_{\xi}(T, \Y(T)) \right].
\]
Then Proposition \ref{thm:Fredholm} gives
a multitime Fredholm determinant expression
to this as
\[
\Psi_{\xi}^{\t}[\f]
=\mathop{{\rm Det}}_
{\bigotimes_{m=1}^M L^2(S, \lambda_{t_m})}
\Big[I-(1-e^{f_{\cdot}}) \sK_{\xi} \Big]
\]
with (\ref{eqn:KK1}). 
By Definition \ref{thm:determinantal},
the proof is completed. \qed

\subsection{Three applications}
\label{sec:three_applications}

In order to show applications of Theorem \ref{thm:DM_det},
we consider the following three kinds of interacting particle systems,
$((\Xi(t))_{t \geq 0}, \bP_{\xi}), \xi=\sum_{j=1}^N \delta_{u_j} \in \Conf_0(S)$ with
$\Xi(t, \cdot)=\sum_{j=1}^N \delta_{X_j(t)}(\cdot), t \geq 0$.
For each system of SDEs, $(B_j(t))_{t \geq 0}, 1 \leq j \leq N$
denote a set of independent one-dimensional standard Brownian motions 
(BMs) started at 0.
(From now on, BM means a one-dimensional standard
Brownian motion unless specially mentioned.)  
\begin{description}
\item[Process 1 :] \,
Noncolliding Brownian motions (the Dyson model with $\beta=2$);
\begin{align}
S &= \R,
\nonumber\\
X_j(t) &= u_j+ B_j(t)+ 
\sum_{\substack{ 1 \leq k \leq N,\\ k \not= j}} \int_0^t
\frac{ds}{X_j(s)-X_k(s)},\quad 1 \leq j \leq N, \quad t \geq 0.
\label{eqn:noncollBM}
\end{align}

\item[Process 2 :] \,
Noncolliding squared Bessel processes 
(BESQ$^{(\nu)}$) with $\nu \in (-1, \infty)$
(the Bru--Wishart process with $\beta=2$); 
\begin{align}
S &= \R_{\geq 0} := \{x \in \R: x \geq 0\},
\nonumber\\
X_j(t) &= u_j+ \int_0^t 2 \sqrt{X_j(s)} 
d B_j(s) + 2 (\nu+1) t
\nonumber\\
& 
+ \sum_{\substack{1 \leq k \leq N, \\ k \not= j}}
\int_0^t
\frac{4 X_j(s)ds}{X_j(s)-X_k(s)},
\quad 1 \leq j \leq N, \quad t \geq 0,
\label{eqn:noncollBESQ}
\end{align}
where if $-1 < \nu < 0$
the reflection boundary condition is assumed at
the origin.

\item[Process 3 :] \,
Noncolliding BM on a circle with a radius $r>0$;
(a trigonometric extension of the Dyson motion model
with $\beta=2$); \\
We solve the SDEs
\begin{equation}
\check{X}_j(t) = u_j+ B_j(t)
+\frac{1}{2 r} \sum_{\substack{1 \leq k \leq N, \cr k \not= j}}
\int_0^t
\cot \left( \frac{\check{X}_j(s)-\check{X}_k(s)}{2 r} \right) ds,
\quad 1 \leq j \leq N, \quad t \geq 0,
\label{eqn:noncollBM_S1}
\end{equation}
on $\R$ and then define the process on
$S=[0, 2 \pi r)$
by 
\begin{equation}
X_j(t)=\check{X}_j(t) \quad
\mbox{mod $2 \pi r$},
\quad 1 \leq j \leq N, \quad t \geq 0.
\label{eqn:noncollBM_S1b}
\end{equation}

Note that $\cot (x/2r)$ is a periodic function of $x$ with
period $2 \pi r$.
By the definition (\ref{eqn:noncollBM_S1b}),
measurable functions for Process 3 should be
periodic with period $2 \pi r$ in the following sense.
For $0 \leq t < \infty$, if an $\cF_{\Xi}(t)$-measurable
function $F$ is given in the form (\ref{eqn:F_gm}),
then for any $n \in \Z$,
\[
g_m((x_j+2 \pi r n)_{j=1}^N)
=g_{m}(\x), 
\quad 1 \leq m \leq M.
\]
\end{description}

The system (\ref{eqn:noncollBM_S1}) 
with the identification (\ref{eqn:noncollBM_S1b}) can be regarded as
a dynamical extension of the {\it circular unitary ensemble} (CUE)
of random matrix theory (see Section 11.8 in \cite{Meh04}
and Chapter 11 in \cite{For10}) \cite{HW96,CL01}.
The dynamics was studied in \cite{NF02}
and papers cited therein. 
Process 3 is a trigonometric extension of Process 1
and in the limit $r \to \infty$
Process 3 should be reduced to Process 1.
Since functions used to represent Process 3 in the present 
manuscript 
are all analytic with respect to $r$,
the hyperbolic extension will be similarly discussed \cite{CL01}.

Corresponding to the three interacting
$N$-particle systems, 
we consider the following three kinds of
one-dimensional processes.
The first one is BM on $S=\R$, whose transition density
started at $x \in \R$ is given by
\begin{equation}
p(t, y|x)=
p_{\rm BM}(t, y|x) =  \left\{ \begin{array}{ll}
\displaystyle{
\frac{e^{-(y-x)^2/(2t)}}{\sqrt{2 \pi t}}},
& t>0, y \in \R \cr
& \cr
\delta_x(\{y\}),
& t=0, y \in \R.
\end{array} \right.
\label{eqn:p}
\end{equation}
The second one is BESQ$^{(\nu)}$ with $\nu > \in (-1, \infty)$
on $S=\R_{\geq 0}$, which is 
given by the solution of the SDE, 
\[
Y(t)=u+\int_0^t 2 \sqrt{Y(s)} d B(s)
+2(\nu+1) t, \quad t \geq 0, \quad u > 0, 
\]
where $B$ is BM, and 
if $-1 < \nu < 0$ a reflecting wall
is put at the origin.
The transition density is given by
\begin{align}
p(t,y|x) &= p^{(\nu)}(t, y|x)
\nonumber\\
&= \left\{ \begin{array}{ll}
\displaystyle{
\left( \frac{y}{x} \right)^{\nu/2}
\frac{e^{-(x+y)/(2t)}}{2t}
I_{\nu} \left( \frac{\sqrt{xy}}{t} \right)},
& t>0, x>0, y \in \R_{\geq 0}, \cr
\displaystyle{
\frac{y^{\nu}e^{-y/(2t)}}{(2t)^{\nu+1} \Gamma(\nu+1)}},
& t >0, x=0, y \in \R_{\geq 0}, \cr
& \cr
\delta_x(\{y\}),
& t=0, x, y \in \R_{\geq 0},
\end{array} \right.
\label{eqn:p_nu}
\end{align}
where $I_{\nu}(x)$ is the modified Bessel function
of the first kind \cite{NIST10}
\[
I_{\nu}(x) = \sum_{n=0}^{\infty} 
\frac{(x/2)^{2n+\nu}}{n! \Gamma(n+1+\nu)}, 
\quad x \in \R_{\geq 0}.
\]
The third one is a Markov process on $S=[0, 2\pi r)$,
$r >0$, whose the transition density is given by
\begin{equation}
p(t, y|x)= p^r(t, y|x; N)= \left\{
\begin{array}{ll}
\displaystyle{\sum_{\ell \in \Z} p_{\rm BM}(t, y+2 \pi r \ell|x)}, \quad & \mbox{if $N$ is odd}, \cr
\displaystyle{\sum_{\ell \in \Z} (-1)^{\ell} 
p_{\rm BM}(t, y+2 \pi r \ell|x)}, \quad & \mbox{if $N$ is even},
\end{array} \right.
\label{eqn:prB1}
\end{equation}
$x, y \in [0, 2 \pi r)$, $t \geq 0$,
where $N$ is the number of particles of Process 3.
If we introduce the notation
\begin{equation}
\sigma_N(m)= \left\{
\begin{array}{ll}
m, \quad & \mbox{when $N$ is odd}, \cr
m-1/2, \quad & \mbox{when $N$ is even},
\end{array} \right.
\label{eqn:spin}
\end{equation}
for $m \in \Z$, then (\ref{eqn:prB1}) is written as
\begin{equation}
p^r(t, y|x; N)
= \frac{1}{2 \pi r}
\sum_{\ell \in \Z} e^{-\sigma_N(\ell)^2 t/2r^2
+ \sqrt{-1} \sigma_N(\ell) (y-x)/r}.
\label{eqn:pR2}
\end{equation}
It should be noted that this expression is
found in Nagao and Forrester \cite{NF02}.
This process shows a position on the
circumference $S=[0, 2 \pi r)$ of a Brownian motion
moving around on the circle $\S^1(r)$
(with alternating signed densities if $N$ is even).
In the following, we call this one-dimensional Markov process
$Y(t) \in [0, 2 \pi r)$, $t \geq 0$
simply `BM on $[0, 2 \pi r)$'.

We introduce integral transformations of function $f=f(W)$,  
\begin{equation}
\cI[f(W)|(t,x)]
= \left\{ \begin{array}{ll}
\displaystyle{
\int_{\R} dw \, f(\sqrt{-1} w) q(t,w|x),
}
& \quad \mbox{for Processes 1 and 3}, \cr
\displaystyle{
\int_{\R_{\geq 0}} dw \, f(-w) q^{(\nu)}(t,w|x),
}
& \quad \mbox{for Process 2},
\end{array} \right.
\label{eqn:I}
\end{equation}
with 
\begin{align}
\label{eqn:q}
q(t, w|x) &
\frac{e^{-(ix+w)^2/2t}}{\sqrt{2 \pi t}}, 
\\
\label{eqn:q_nu}
q^{(\nu)}(t,w|x) &= 
\left( \frac{w}{x} \right)^{\nu/2}
\frac{e^{(x-w)/2t}}{2t}
J_{\nu} \left( \frac{\sqrt{xw}}{t} \right),
\quad \nu > \in (-1, \infty), 
\end{align}
where 
$J_{\nu}$ is the Bessel function 
defined by (\ref{eqn:Bessel_function}). 

Note that $W$ in LHS of (\ref{eqn:I}) is a dummy variable,
but it will be useful in order to specify a function $f$.
For example, we can verify the following \cite{Kat14};
\begin{equation}
m_n(t, x) := 
\cI[W^n|(t,x)] = 
\left(\frac{t}{2}\right)^{n/2} H_n \left( \frac{x}{\sqrt{2t}} \right),
\quad t \geq 0, \quad n \in \N_0,
\quad \mbox{for Process 1},
\label{eqn:FMP1}
\end{equation}
where $\{H_n\}_{n \in \N_0}$ are the Hermite polynomials 
(\ref{eqn:Hermite1}),
and for $\nu \in (-1, \infty)$, 
\begin{equation}
m_n^{(\nu)}(t, x) := 
\cI[W^n|(t,x)] = (-1)^n n! (2t)^n
L^{(\nu)}_n \left( \frac{x}{2t} \right),
\quad t \geq 0, \quad n \in \N_0,
\quad \mbox{for Process 2},
\label{eqn:FMP2}
\end{equation}
where $\{L^{(\nu)}_n\}_{n \in \N_0}$ are the Laguerre polynomials 
(\ref{eqn:Laguerre1}).

We introduce sets of entire functions of $z \in \C$ \cite{Lev96}, 
$1 \leq k \leq N$, for
$\u \in \W_N(S)$, $\xi=\sum_{j=1}^N \delta_{u_j}$, $v \in S$,
$r>0$, 
\begin{equation}
\Phi_{\xi}^{v}(z)
=\left\{ \begin{array}{ll}
\displaystyle{
\prod_{\substack{1 \leq \ell \leq N, \cr u_{\ell} \not=v}}
\frac{z-u_{\ell}}{v-u_{\ell}},
}
& \quad \mbox{for Processes 1 and 2}, \cr
\displaystyle{
\prod_{\substack{1 \leq \ell \leq N, \cr u_{\ell} \not=v}}
\frac{\sin((z-u_{\ell})/2r)}{\sin((v-u_{\ell})/2r)},
}
& \quad \mbox{\rm for Process 3}. \cr
\end{array} \right.
\label{eqn:Phi}
\end{equation}
For $v \in S$, set 
\begin{equation}
\cM^{v}_{\xi}(t, x)
=\cI[\Phi_{\xi}^{v}(W)|(t,x)],
\quad (t, x) \in [0, \infty) \times S.
\label{eqn:cMs}
\end{equation}
Then the following is proved \cite{Kat14}.

\begin{prop}
\label{thm:DMR_3examples}
For Processes 1, 2, and 3, 
if $\xi=\sum_{j=1}^N \delta_{u_j} \in \Conf_0(S)$, then
$\cM_{\xi} =\{\cM_{\xi}^v(\cdot, \cdot) :
v \in \{u_1, \dots, u_N\} \}$ satisfies
the conditions {\bf (M1)}--{\bf (M3)}.
\end{prop}

We can prove the following \cite{Kat14}.
\begin{thm}
\label{thm:3_process}
For any $\xi \in \Conf_0(S)$, 
the three processes have DMRs 
associated with $(Y, \cM_{\xi})$
such that
\[
\mbox{$Y$ is given by} \, \left\{ \begin{array}{ll}
\mbox{BM on $\R$}, & \mbox{for Process 1}, \cr
\mbox{BESQ$^{(\nu)}$ on $\R_{\geq 0}$}, & \mbox{for Process 2}, \cr
\mbox{BM on $[0, 2 \pi r)$}, & \mbox{for Process 3}, \cr
\end{array} \right. 
\]
and $\cM_{\xi}=\{\cM_{\xi}^v(\cdot, \cdot)\}, v \in S$ is given by
(\ref{eqn:cMs}).
Then they are all DSPs
with the spatio-temporal correlation kernels 
\[
\bK_{\xi}(s,x;t,y)
=\int_{S} \xi(dv) p(s, x|v) \cM_{\xi}^{v}(t,y)
- \1(s>t) p(s-t,x|y),
\quad (s,x), (t,y) \in [0, \infty) \times S, 
\]
with 
\[
\cM^{v}_{\xi}(t, x)
=\cI[\Phi_{\xi}^{v}(W)|(t,x)],
\quad (t, x) \in [0, \infty) \times S,
\]
where $p$ is the transition density of the process $Y$.

\end{thm}

\subsection{Martingales for configurations with multiple points 
\label{sec:multi}}
For general $\xi \in \Conf(S)$
with $\xi(S)=N < \infty$, define
$\supp \xi=\{x \in S : \xi(x) > 0\}$ and let
$
\xi_{*}(\cdot)=\sum_{v \in \supp \xi} \delta_{v}(\cdot).
$
For $s \in [0, \infty)$, $v, x \in S$, $z, \zeta \in \C$, let
\begin{equation}
\phi_{\xi}^{v}((s,x);z,\zeta)
=\left\{ \begin{array}{ll}
\displaystyle{
\frac{p(s,x|\zeta)}{p(s,x|v)} \frac{1}{z-\zeta}
\prod_{\ell=1}^N
\frac{z-u_{\ell}}{\zeta-u_{\ell}},
}
& \, \mbox{for Processes 1 and 2}, \cr
\displaystyle{
\frac{p(s,x|\zeta)}{p(s,x|v)} \frac{1}{z-\zeta}
\prod_{\ell=1}^N
\frac{\sin((z-u_{\ell})/2r)}{\sin((\zeta-u_{\ell})/2r),}
}
& \, \mbox{\rm for Process 3}, \cr
\end{array} \right.
\label{eqn:phi}
\end{equation}
and
\begin{align}
\Phi_{\xi}^{v}((s,x); z) &= \frac{1}{2 \pi \sqrt{-1}}
\oint_{C(\delta_{v})} d \zeta \, 
\phi^v_{\xi}((s,x); z, \zeta)
\nonumber\\
&= {\rm Res} \,\Big[\phi^v_{\xi}((s,x); z, \zeta); \zeta=v \Big],
\label{eqn:PhiB}
\end{align}
where
$C(\delta_{v})$ is a closed contour on the complex plane $\C$
encircling a point $v$ on $S$
once in the positive direction. Define
\begin{equation}
\cM_{\xi}^{u}((s, x)|(t,y))
=\cI \left[\Phi_{\xi}^{u}((s,x); W) \Big| (t,y) \right],
\quad (s,x), (t,y) \in [0, \infty) \times S.
\label{eqn:cMB}
\end{equation}
Then it is easy to see that (\ref{eqn:KK1}) is rewritten as
\begin{equation}
\bK_{\xi}(s,x;t,y)
=\int_{S} \xi_{*}(dv) p(s, x|v) \cM_{\xi}^{v}((s,x)|(t,y))
- \1(s>t) p(s-t,x|y),
\label{eqn:K2}
\end{equation}
$(s,x), (t,y) \in [0, \infty) \times S$.

We note that, even though the systems of SDEs (\ref{eqn:noncollBM})-(\ref{eqn:noncollBM_S1})
cannot be solved for any initial configuration with
multiple points, $\xi \in \Conf(S) \setminus \Conf_0(S)$,
the kernel (\ref{eqn:K2}) with (\ref{eqn:cMB})
is bounded and integrable also for $\xi \in \Conf(S) \setminus \Conf_0(S)$.
Therefore, spatio-temporal correlations are given
by (\ref{eqn:rho1}) for any $0 \leq t_1 < \dots < t_M < \infty, M \in \N$
and finite-dimensional distributions 
are determined.
\begin{prop}
\label{thm:entrance_law}
Also for $\xi \in \Conf \setminus \Conf_0$,
the DSPs
with the spatio-temporal correlation kernels (\ref{eqn:K2})
are well-defined.
They provide the entrance laws for
the processes $((\Xi(t))_{t \geq 0}, \bP_{\xi})$.
\end{prop}

In order to give examples
of Proposition \ref{thm:entrance_law}, 
here we study the extreme case such that
all $N$ points are concentrated on an origin,
\begin{equation}
\xi=N \delta_0 \quad
\Longleftrightarrow \quad
\xi_{*}=\delta_0 \quad
\mbox{with} \quad
\xi(\{0\})=N.
\label{eqn:xi01}
\end{equation}
We consider Processes 1 and 2.
For (\ref{eqn:xi01}), (\ref{eqn:phi}) and (\ref{eqn:PhiB})
become
\begin{align*}
\phi_{N \delta_0}^0((s,x); z, \zeta)
&= \frac{p(s,x|\zeta)}{p(s,x|0)}
\frac{1}{z-\zeta} \left( \frac{z}{\zeta} \right)^N
\nonumber\\
&= \frac{p(s,x|\zeta)}{p(s,x|0)}
\sum_{\ell=0}^{\infty} \frac{z^{N-\ell-1}}{\zeta^{N-\ell}},
\end{align*}
and
\begin{align}
\Phi_{N \delta_0}^0((s,x);z) &=
\frac{1}{p(s,x|0)} \sum_{\ell=0}^{\infty} z^{N-\ell-1}
\frac{1}{2 \pi \sqrt{-1}} \oint_{C(\delta_0)} d \zeta \,
\frac{p(s,x|\zeta)}{\zeta^{N-\ell}}
\nonumber\\
&=
\frac{1}{p(s,x|0)} \sum_{\ell=0}^{N-1} z^{N-\ell-1}
\frac{1}{2 \pi \sqrt{-1}} \oint_{C(\delta_0)} d \zeta \,
\frac{p(s,x|\zeta)}{\zeta^{N-\ell}},
\label{eqn:xi03}
\end{align}
since the integrands are holomorphic when $\ell \geq N$,
where we assume $\nu \in (-1, \infty)$
and $C(\delta_0)$ is interpreted as 
$\lim_{\varepsilon \downarrow 0} C(\delta_{\varepsilon})$
for BESQ$^{(\nu)}$. 

For BM with the transition density (\ref{eqn:p}),
(\ref{eqn:xi03}) gives
\begin{align*}
\Phi_{N \delta_0}^0((s,x);z)
&= \sum_{\ell=0}^{N-1} z^{N-\ell-1}
\frac{1}{2 \pi \sqrt{-1}}
\oint_{C(\delta_0)} d \zeta \,
\frac{e^{x\zeta/s-\zeta^2/2s}}{\zeta^{N-\ell}}
\nonumber\\
&= \sum_{\ell=0}^{N-1} 
\left( \frac{z}{\sqrt{2s}} \right)^{N-\ell-1}
\frac{1}{2 \pi \sqrt{-1}}
\oint_{C(\delta_0)} d \eta \,
\frac{e^{2(x/\sqrt{2s}) \eta-\eta^2}}{\eta^{N-\ell}}
\nonumber\\
&= \sum_{\ell=0}^{N-1} 
\left( \frac{z}{\sqrt{2s}} \right)^{N-\ell-1}
\frac{1}{(N-\ell-1)!} H_{N-\ell-1} \left(\frac{x}{\sqrt{2s}} \right),
\end{align*}
where we have used the contour integral representation
of the Hermite polynomials \cite{Sze91}
\[
H_n(x)=\frac{n!}{2 \pi \sqrt{-1}} \oint_{C(\delta_0)} d \eta
\frac{e^{2 x \eta-\eta^2}}{\eta^{n+1}},
\quad n \in \N_0, \quad x \in \R.
\]
Thus by (\ref{eqn:FMP1}) its integral transformation is calculated as
\begin{align*}
& \cI \left[ \left.
\Phi_{N \delta_0}^0((s,x);W) \right| (t, y) \right]
\nonumber\\
& \quad 
= \sum_{\ell=0}^{N-1} \frac{1}{(N-\ell-1)!} H_{N-\ell-1}
\left( \frac{x}{\sqrt{2s}} \right)
\frac{1}{(2s)^{(N-\ell-1)/2}}
\cI[W^{N-\ell-1}|(t,y)]
\nonumber\\
& \quad 
= \sum_{\ell=0}^{N-1} \frac{1}{(N-\ell-1)!} H_{N-\ell-1}
\left( \frac{x}{\sqrt{2s}} \right)
\frac{1}{(2s)^{(N-\ell-1)/2}} 
m_{N-\ell-1}(t,y)
\nonumber\\
& \quad 
= \sum_{\ell=0}^{N-1} \frac{1}{(N-\ell-1)! 2^{N-\ell-1}} 
\left( \frac{t}{s} \right)^{(N-\ell-1)/2}
H_{N-\ell-1}
\left( \frac{x}{\sqrt{2s}} \right)
H_{N-\ell-1}
\left( \frac{y}{\sqrt{2t}} \right).
\end{align*}
Then we obtain the following,
\begin{align}
& \cM_{N \delta_0}^0((s,x)|(t,Y(t)))
=\sum_{n=0}^{N-1} \frac{1}{n! s^n} m_n(s,x) m_n(t,Y(t))
\nonumber\\ 
& \qquad \qquad =
\sum_{n=0}^{N-1} 
 \left(\frac{t}{s}\right)^{n/2}
\varphi_n \left( \frac{x}{\sqrt{2s}} \right)
\varphi_n \left( \frac{Y(t)}{\sqrt{2t}} \right),
\label{eqn:Hermite3B}
\end{align}
where
\[
\varphi_n(x)=\frac{1}{\sqrt{2^n n!}} H_n(x), \quad
n \in \N, \quad x \in \R.
\]
Similarly, for BESQ$^{(\nu)}, \nu \in (-1, \infty)$
with the transition density (\ref{eqn:p_nu}),
we obtain
\begin{align*}
\Phi_{N \delta_0}^{0}((s,x); z)
&= \frac{(2s)^{\nu} \Gamma(\nu+1)}{x^{\nu/2}}
\sum_{\ell=0}^{N-1} z^{N-\ell-1}
\frac{1}{2 \pi \sqrt{-1}} \oint_{C(\delta_0)} d \zeta \,
\frac{e^{-\zeta/2s}}{\zeta^{N-\ell+\nu/2}}
I_{\nu} \left( \frac{\sqrt{x \zeta}}{s}\right)
\nonumber\\
&= \Gamma(\nu+1) \sum_{\ell=0}^{N-1}
\left( -\frac{z}{2s} \right)^{N-\ell-1}
\frac{1}{\Gamma(N-\ell+\nu)}
L^{(\nu)}_{N-\ell-1} \left( \frac{x}{2s} \right),
\end{align*}
where we have used the contour integral representation
of the Laguerre polynomials
\[
L_n^{(\nu)}(x)
= \frac{\Gamma(n+\nu+1)}{x^{\nu/2}}
\frac{1}{2 \pi \sqrt{-1}} \oint_{C(\delta_0)}
d \eta \,
\frac{e^{\eta}}{\eta^{n+1+\nu/2}} J_{\nu}(2 \sqrt{\eta x}). 
\]
By (\ref{eqn:FMP2}), we have
\begin{align}
& \cM_{N \delta_0}^{0}((s,x)|(t, Y(t))
= \cI^{(\nu)} \left[ \left.
\Phi_{N \delta_0}^{0}((s,x); W) \right| (t, Y(t)) \right]
\nonumber\\
& \quad = \Gamma(\nu+1)
\sum_{n=0}^{N-1} \frac{1}{n! \Gamma(n+\nu+1) (2s)^{2n}}
m^{(\nu)}_n(s,x) m^{(\nu)}_n \left(t, Y(t) \right)
\nonumber\\
& \quad = 
\sum_{n=0}^{N-1}
\left( \frac{t}{s} \right)^n
\phi_n^{(\nu)} \left( \frac{x}{2s} \right)
\phi_n^{(\nu)} \left( \frac{Y(t)}{2t} \right),
\label{eqn:Laguerre3B}
\end{align}
where
\[
\phi_n^{(\nu)}(x)=\sqrt{\frac{n! \Gamma(\nu+1)}{\Gamma(n+\nu+1)}}
L_n^{(\nu)}(x),
\quad n \in \N_0, \quad x \in \R_{\geq 0}.
\]
The processes (\ref{eqn:Hermite3B}) and (\ref{eqn:Laguerre3B})
are continuous martingales.
Then we see
\[
\rE \left[ \cM_{N \delta_0}^0((s,x)|(t, \Y(t))) \right]
= \rE \left[ \cM_{N \delta_0}^0((s,x)|(0, \Y(0))) \right]=1
\]
for $(s,x) \in [0, T] \times S$, $0 \leq t < \infty$.

By the formula (\ref{eqn:K2}), we obtain the
spatio-temporal correlation kernels for Process 1 as
\begin{align*}
&\bK_{N \delta_0}^{1}(s,x;t,y) 
= p_{\rm BM}(s,x|0) \cM_{N \delta_0}^0((s,x)| (t,y))
-\1(s>t) p_{\rm BM}(s-t, x|y)
\nonumber\\
& \qquad= \frac{e^{-x^2/4s}/s^{1/4}}{e^{-y^2/4t}/t^{1/4}}
\bK^{(N)}_{\rm Hermite}(s,x;t,y)
\left(\frac{1}{\sqrt{2s}} \frac{e^{-(x/\sqrt{2s})^2}}{\sqrt{\pi}} \right)^{1/2}
\left(\frac{1}{\sqrt{2t}} \frac{e^{-(y/\sqrt{2t})^2}}{\sqrt{\pi}} \right)^{1/2}
\end{align*}
with
\begin{align*}
\bK^{(N)}_{\rm Hermite}(s,x;t,y)
&= 
\sum_{n=0}^{N-1} \left( \frac{t}{s} \right)^{n/2}
\varphi_n \left( \frac{x}{\sqrt{2s}} \right)
\varphi_n \left( \frac{y}{\sqrt{2t}} \right)
\nonumber\\
&- \1(s>t) 
\sum_{n=0}^{\infty} \left( \frac{t}{s} \right)^{n/2}
\varphi_n \left( \frac{x}{\sqrt{2s}} \right)
\varphi_n \left( \frac{y}{\sqrt{2t}} \right),
\end{align*}
and for Process 2
\begin{align*}
&\bK^{2}_{N \delta_0}(s,x;t,y) 
= p^{(\nu)}(s,x|0) \cM_{N \delta_0}^{0} ((s,x)| (t,y)) 
-\1(s>t) p^{(\nu)}(s-t, x|y)
\nonumber\\
&\qquad =
\frac{(x/2s)^{\nu/2} e^{-x/4s}/s^{1/2}}{(y/2t)^{\nu/2} e^{-y/4t}/t^{1/2}}
\bK^{(\nu, N)}_{\rm Laguerre} (s,x;t,y)
\left( \frac{(x/2s)^{\nu} e^{-x/2s}}
{\Gamma(\nu+1)} \frac{1}{2s} \right)^{1/2}
\left( \frac{(y/2t)^{\nu} e^{-y/2t}}
{\Gamma(\nu+1)} \frac{1}{2t} \right)^{1/2}, 
\end{align*}
with
\begin{align*}
\bK^{(\nu, N)}_{\rm Laguerre}(s,x;t,y) &=
\sum_{n=0}^{N-1} \left( \frac{t}{s} \right)^{n}
\phi^{(\nu)}_n \left( \frac{x}{2s} \right)
\phi^{(\nu)}_n \left( \frac{y}{2t} \right)
\nonumber\\
&- \1(s>t) 
\sum_{n=0}^{\infty} \left( \frac{t}{s} \right)^{n}
\varphi^{(\nu)}_n \left( \frac{x}{2s} \right)
\varphi^{(\nu)}_n \left( \frac{y}{2t} \right).
\end{align*}

The above results are summarized as follows.
\begin{prop}
\label{thm:extended_kernels}
\begin{description}
\item{\rm (i)} 
The process 1 (the noncolliding Brownian motions
= Dyson model with $\beta=2$) on $\R$
starting from $N \delta_0$ is the DSP,
$((\Xi(t))_{t \geq 0}, \bK^{(N)}_{\rm Hermite}, 
(d \lambda_t(dx))_{t \geq 0})$, 
with the spatio-temporal correlation kernel 
\[
\bK^{(N)}_{\rm Hermite}
(s, x; t, y)
= \begin{cases}
\displaystyle{
\sum_{n=0}^{N-1} \left( \frac{t}{s} \right)^{n/2}
\varphi_n \left( \frac{x}{\sqrt{2s}} \right)
\varphi_n \left( \frac{y}{\sqrt{2t}} \right)},
& \mbox{if $s < t$},
\cr
\displaystyle{
\sum_{n=0}^{N-1} 
\varphi_n \left( \frac{x}{\sqrt{2t}} \right)
\varphi_n \left( \frac{y}{\sqrt{2t}} \right)}, 
& \mbox{if $s = t$},
\cr
\displaystyle{
-\sum_{n=N}^{\infty} \left( \frac{t}{s} \right)^{n/2}
\varphi_n \left( \frac{x}{\sqrt{2s}} \right)
\varphi_n \left( \frac{y}{\sqrt{2t}} \right)},
& \mbox{if $s > t$},
\end{cases}
\]
$(s,x), (t,y) \in [0, \infty) \times \R$, where
$\varphi(x)=H_n(x)/\sqrt{2^n n!}, n \in \N_0$,
and with the time-dependent background measure
\begin{align*}
d \lambda_t(dx)
&= \sqrt{2 t} \circ \lambda_{\rm N(0,1/2)}(dx)
= \frac{e^{-(x/\sqrt{2t})^2}}{\sqrt{\pi}} \frac{dx}{\sqrt{2t}}
\nonumber\\
&=\frac{e^{-x^2/(2t)}}{\sqrt{2 \pi t}} dx 
=p_{\rm BM}(t, x|0) dx,
\quad (t, x) \in [0, \infty) \times \R.
\end{align*}
\item{\rm (ii)} 
The process 2 (the noncolliding BESQ$^{(\nu)}$
= Bru--Wishart process with $\beta=2$), $\nu \in (-1, \infty)$, 
on $\R_{\geq 0}$
starting from $N \delta_0$ is the DSP,
$((\widehat{\Xi}(t))_{t \geq 0}, \bK^{(N)}_{\rm Laguerre}, 
(d \widehat{\lambda}_t(dx))_{t \geq 0})$, 
with the spatio-temporal correlation kernel 
\[
\bK^{(N)}_{\rm Laguerre}
(s, x; t, y)
= \begin{cases}
\displaystyle{
\sum_{n=0}^{N-1} \left( \frac{t}{s} \right)^{n}
\phi_n \left( \frac{x}{2s} \right)
\phi_n \left( \frac{y}{2t} \right)},
& \mbox{if $s < t$},
\cr
\displaystyle{
\sum_{n=0}^{N-1} 
\phi_n \left( \frac{x}{2t} \right)
\phi_n \left( \frac{y}{2t} \right)}, 
& \mbox{if $s = t$},
\cr
\displaystyle{
-\sum_{n=N}^{\infty} \left( \frac{t}{s} \right)^{n}
\phi_n \left( \frac{x}{2s} \right)
\phi_n \left( \frac{y}{2t} \right)},
& \mbox{if $s > t$},
\end{cases}
\]
$(s,x), (t,y) \in [0, \infty) \times \R_{\geq 0}$, where
$\phi(x)=\sqrt{n! \Gamma(\nu+1)/\Gamma(n+\nu+1)}
L^{(\nu)}_n, n \in \N_0$,
and with the time-dependent background measure
\begin{align*}
d \widehat{\lambda}_t(dx)
&= (2t) \circ \lambda_{\Gamma(\nu+1,1)}(dx)
= \frac{1}{\Gamma(\nu+1)} \left( \frac{x}{2t} \right)^{\nu} 
e^{-x/(2t)} \frac{dx}{2t}
\nonumber\\
&=\frac{x^{\nu} e^{-x/(2t)}}{(2t)^{\nu+1} \Gamma(\nu+1)} dx
=p^{(\nu)}(t, x|0) dx,
\quad (t, x) \in [0, \infty) \times \R_{\geq 0}.
\end{align*}
\end{description}
\end{prop}

Note that,
for each time $t \in (0, \infty)$,
\begin{align*}
\bK^{(N)}_{\rm Hermite}(t, x; t, y)
&= K^{(N)}_{\rm Hermite} 
\left( \frac{x}{\sqrt{2t}}, \frac{y}{\sqrt{2t}} \right)
=\sqrt{2t} \circ K^{(N)}_{\rm Hermite},
\quad x, y \in \R,
\nonumber\\
\bK^{(N)}_{\rm Laguerre}(t, x; t, y)
&= K^{(N)}_{\rm Laguerre} 
\left( \frac{x}{2t}, \frac{y}{2t} \right)
=(2t) \circ K^{(N)}_{\rm Laguerre},
\quad x, y \in \R_{\geq 0}. 
\end{align*}
The spatio-temporal correlation kernels
$\bK^{(N)}_{\rm Hermite}$
and $\bK^{(N)}_{\rm Laguerre}$
are called the {\it extended Hermite and Laguerre kernels},
respectively, in random matrix theory (see, for instance, \cite{For10}).
Hence at each time $t \in (0, \infty)$, 
\begin{align*}
(\Xi(t), \bK^{(N)}_{\rm Hermite}, \lambda_t(dx))
=
(\Xi, \sqrt{2t} \circ K^{(N)}_{\rm Hermite}, 
\sqrt{2t} \circ \lambda_{\rN(0, 1/2)}(dx)),
\nonumber\\
(\widehat{\Xi}(t), \bK^{(N)}_{\rm Laguerre}, \widehat{\lambda}_t(dx))
=
(\Xi, (2t) \circ K^{(N)}_{\rm Laguerre}, 
(2t) \circ \lambda_{\Gamma(\nu+1, \nu)}(dx)),
\end{align*}
where $(\Xi, K^{(N)}_{\rm Hermite}, \lambda_{\rN(0, 1/2)}(dx))$
and $(\Xi, K^{(N)}_{\rm Laguerre}, \lambda_{\Gamma(\nu+1, \nu)}
(dx))$
and the DPPs studied in Section \ref{sec:classical}.

Here we would like to emphasize the fact that
these spatio-temporal kernels have been derived here by not following
the `bi(multiple)-orthogonal-function method' \cite{KT10}, 
but by only using proper martingales
determined by the chosen initial configuration (\ref{eqn:xi01}).
The above kernels determine the finite-dimensional distributions
and specify the entrance laws
for the systems of SDEs (\ref{eqn:noncollBM})
and (\ref{eqn:noncollBESQ}) from the state $N \delta_0$ \cite{KT04}.

\subsection{Relaxation phenomenon in Process 3}
\label{sec:relax}

As an application of Theorem \ref{thm:3_process},
here we study a typical non-equilibrium dynamics of Process 3,
that is, a relaxation phenomenon to the equilibrium.
For Processes 1 and 2, see \cite{KT09,KT10,KT11}.

For Process 3 we consider the following special initial configuration
\begin{equation}
\eta(\cdot)=\sum_{j=1}^N \delta_{w_j}(\cdot)
\quad \mbox{with} \quad 
w_j=\frac{2 \pi r}{N}(j-1), \quad
1 \leq j \leq N.
\label{eqn:ed1}
\end{equation}
It is an unlabeled configuration with equidistant spacing
on $[0, 2 \pi r)$.
In this case the entire function for Process 3 
given by (\ref{eqn:Phi}) becomes
\begin{align*}
\Phi_{\eta}^{w_k}(z)
&= \prod_{\substack{1 \leq \ell \leq N, \cr \ell \not=k}}
\frac{\sin(z/2r-(\ell-1)\pi/N)}{\sin((k-\ell)\pi/N)}
\nonumber\\
&= \prod_{n=1}^{N-1}
\frac{\sin[\{z/2r-(k-2)\pi/N\}+(n-1)\pi/N]}{\sin(n \pi/N)}.
\end{align*}
We use the product formulas
\[
\prod_{n=1}^{N-1} \sin \left( \frac{n \pi}{N} \right) 
=\frac{N}{2^{N-1}},
\qquad
\prod_{n=1}^{N} \sin \left[ x+\frac{(n-1)\pi}{N} \right]
=\frac{\sin(Nx)}{2^{N-1}},
\]
and obtain
\[
\Phi_{\eta}^{w_k}(z)
=\frac{1}{N} \frac{\sin[N\{z-2 \pi r(k-1)/N\}/2r]}
{\sin[\{z-2\pi r(k-1)/N\}/2r]}.
\]

It is easy to confirm the equality
\[
\frac{\sin(Nx)}{\sin x}
=\sum_{\substack{m \in \Z, \cr |\sigma_N(m)| \leq (N-1)/2}}
e^{2 \sqrt{-1} \sigma_N(m) x},
\quad N \in \N,
\]
where $\sigma_N$ is defined by (\ref{eqn:spin}). 
Then, the martingale function is given by
\begin{align*}
\cM_{\eta}^{w_k}(t,y) &=
\int_{\R} d \widetilde{y} \,
\frac{1}{N} 
\sum_{\substack{m \in \Z, \cr |\sigma_N(m)| \leq (N-1)/2}}
e^{2 \sqrt{-1} \sigma_N(m)\{y+\sqrt{-1} \widetilde{y}-2 \pi r (k-1)/N\}/2r}
p_{\rm BM}(t, \widetilde{y}|0)
\nonumber\\
&= \frac{1}{N} \sum_{\substack{m \in \Z, \cr |\sigma_N(m)|\leq (N-1)/2}}
e^{\sigma_N(m)^2 t/2r^2-\sqrt{-1} \sigma_N(m) y/r+2\sqrt{-1} (k-1) \sigma_N(m) \pi/N},
\end{align*}
where $p_{\rm BM}$ is given by (\ref{eqn:p}).

Now the spatio-temporal 
correlation kernel 
with respect to the Lebesgue measure $dx$ on
$[0, 2 \pi r)$
is obtained by the formula (\ref{eqn:KK1}),
\[
\widetilde{\bK}_{\eta}(s,x;t,y)
= \cG_{\eta}(s,x;t,y)-\1(s>t) p^r(s-t, x|y; N)
\]
with
\begin{align*}
\cG_{\eta}(s,x;t,y)
&= \sum_{k=1}^N p^r(s, x| w_k; N)
\cM_{\eta}^{w_k}(t,y)
\nonumber\\
&= \frac{1}{N} \frac{1}{2 \pi r}
\sum_{k=1}^N 
\sum_{\ell \in \Z} 
\sum_{\substack{m \in \Z, \cr |\sigma_N(m)|\leq (N-1)/2}}
h_{k,\ell,m}(s,x;t,y),
\end{align*}
where we use the expression (\ref{eqn:pR2})
for $p^r$ and obtain
\begin{align}
&h_{k,\ell,m}(s,x;t,y)
\nonumber\\
& \quad 
= e^{\{\sigma_N(m)^2 t - \sigma_N(\ell)^2 s\}/2r^2
+\sqrt{-1} \{\sigma_N(\ell) x- \sigma_N(m) y\}/r+2\sqrt{-1} (k-1)\{\sigma_N(m)-\sigma_N(\ell)\} \pi/N}.
\label{eqn:calK1}
\end{align}
By the equality
\[
\sum_{k=1}^N e^{2\sqrt{-1} (k-1)\{\sigma_N(m)- \sigma_N(\ell)\} \pi/N}
=N \sum_{k \in \Z} \1(\ell=m+k N),
\]
we obtain the following decomposition,
\[
\cG_{\eta}(s, x; t, y)=\sum_{k \in \Z} \cG_{\eta}^{(k)}(s, x; t, y)
\]
with
\begin{align*}
\cG_{\eta}^{(k)}(s, x; t, y)
&= \frac{1}{2 \pi r}
\sum_{\substack{m \in \Z, \cr |\sigma_N(m)|\leq (N-1)/2}}
e^{-\{\sigma_N(m+kN)^2-\sigma_N(m)^2\}s/2r^2}
\nonumber\\
& \qquad \qquad \qquad \times
e^{\sigma_N(m)^2(t-s)/2r^2-\sqrt{-1} \{\sigma_N(m)y-\sigma_N(m+kN)x\}/r}.
\end{align*}
Since $\sigma_N(m+kN)^2 > \sigma_N(m)^2$ if $m \in \Z, |\sigma_N(m)| \leq (N-1)/2$
and $k \not=0$, we see that
for $(s, x), (t, y) \in [0, \infty) \times [0, 2 \pi r)$
\[
\lim_{T \to \infty} \cG_{\eta}^{(k)}(s+T, x; t+T, y)
= \left\{ \begin{array}{ll}
\cG_{\rm eq}(t-s,y-x),
\quad &\mbox{if $k=0$}, \cr
0, \quad &\mbox{otherwise},
\end{array} \right.
\]
where 
\[
\cG_{\rm eq}(t, x)
= \frac{1}{2 \pi r} 
\sum_{\substack{m \in \Z, \cr |\sigma_N(m)|\leq (N-1)/2}}
e^{\sigma_N(m)^2 t/2r^2 -\sqrt{-1} \sigma_N(m) x/r}.
\]
In particular, when $s=t$ we have
\[
\cG_{\rm eq}(0, x)
= \frac{1}{2 \pi r} 
\sum_{\substack{m \in \Z, \cr |\sigma_N(m)|\leq (N-1)/2}}
e^{-\sqrt{-1} \sigma_N(m) x/r}
= \frac{1}{2 \pi r} 
\frac{\sin(N x/2r)}{\sin(x/2r)}.
\]
The results are summarized as follows.
Define
\[
\re_n(t,x) := e^{n^2 t-\sqrt{-1} nx}.
\]

\begin{prop}
\label{thm:relax}
Let $((\Xi(t))_{t \geq 0}, \bP_{\eta})$ be the Process 3
started at the configuration (\ref{eqn:ed1}).
It is a DSP, 
$((\Xi(t))_{t \geq 0}, \bK_{\eta}, \lambda_{[0, 2 \pi r)}(dx))$,
with the spatio-temporal kernel
\begin{align*}
\bK_{\eta}(s,x;t,y)
&=  \frac{1}{N} 
\sum_{k=1}^N 
\sum_{\ell \in \Z} 
\sum_{\substack{m \in \Z, \cr |\sigma_N(m)|\leq (N-1)/2}}
\frac{\re_{\sigma_N(m)}(t/2r^2, y/r-2\pi (k-1)/N)}
{\re_{\sigma_N(\ell)}(s/2r^2, x/r -2 \pi(k-1)/N)}
\nonumber\\
&-\1(s>t)
\sum_{\ell \in \Z}
\frac{\re_{\sigma_N(\ell)}(t/2r^2, y/r)}{\re_{\sigma_N(\ell)}(s/2r^2, x/r)},
\end{align*}
and with the uniform measure on the circle
$[0, 2 \pi r)$, 
$\lambda_{[0, 2 \pi r)}(dx)=dx/(2\pi r)= r \circ \lambda^{\AN}(dx)$.
This DSP shows a relaxation
phenomenon to an equilibrium DSP, 
$((\Xi(t))_{t \geq 0}, \bK^r_{\rm eq}, \lambda_{[0, 2 \pi r)}(dx))$,
where 
\[
\bK_{\rm eq}^r(s, x; t, y)
= \left\{ \begin{array}{ll}
\displaystyle{ 
\sum_{\substack{m \in \Z, \cr |\sigma_N(m)|\leq (N-1)/2}}
\frac{\re_{\sigma_N(m)}(t/2r^2, y/r)}{\re_{\sigma_N(m)}(s/2r^2, x/r)} }, 
& \mbox{if $s < t$}, \cr
& \cr
\displaystyle{
\frac{\sin[N(y-x)/2r]}{\sin[(y-x)/2r]}},
& \mbox{if $s=t$}, \cr
& \cr
\displaystyle{ -
\sum_{\substack{m \in \Z, \cr |\sigma_N(m)| > (N-1)/2}}
\frac{\re_{\sigma_N(m)}(t/2r^2, y/r)}{\re_{\sigma_N(m)}(s/2r^2, x/r)} 
},
& \mbox{if $s > t$}
\end{array} \right.
\]
for $(s,x), (t,y) \in [0, \infty) \times [0, 2 \pi r)$.
\end{prop}

Note that
\[
\bK^r_{\rm eq}(t, x; t, y)
=r \circ K^{\AN}(x, y),
\quad t \geq 0, \quad x, y \in [0, 2\pi r), 
\]
where $(\Xi, K^{\AN}, \lambda_{[0, 2 \pi)}(dx))$
was studied in Section \ref{sec:Lie_group}.
The above equilibrium DSP, \\
$((\Xi(t))_{t \geq 0}, \bK^r_{\rm eq}, \lambda_{[0, 2 \pi r)}(dx))$, 
is reversible with respect to 
the DPP, 
$(\Xi, r \circ K^{\AN}, r \circ \lambda^{\AN}(dx))$.
\vskip 0.3cm

\clearpage

\SSC
{Multiple Schramm--Loewner Evolutions (SLEs)
and Gaussian Free Fields (GFFs)}
\label{sec:SLE_GFF}
\subsection{Imaginary surface and SLE}
\label{sec:IS}

The present study has been motivated by the recent work
by Sheffield on the quantum gravity zipper and the AC geometry \cite{She16}
and a series of papers by Miller and Sheffield on the 
imaginary geometry \cite{MS16a,MS16b,MS16c,MS17}.
In both of them, a {\it Gaussian free field (GFF)} on 
a simply connected proper subdomain $D$ of the complex plane $\C$ 
(see, for instance, \cite{She07}) is coupled
with a {\it Schramm--Loewner evolution (SLE)} driven
by a Brownian motion moving 
on the boundary $\partial D$ \cite{Sch00,LSW04,Law05}.

Consider a simply connected domain $D \subsetneq \C$
and write $\cC_{\rm c}^{\infty}(D)$ for the space of real smooth functions
on $D$ with compact support.
Assume $h \in \cC_{\rm c}^{\infty}(D)$ and consider a smooth vector field
$e^{\sqrt{-1} (h/\chi + \theta)}$ with parameters $\chi, \theta \in \R$.
Then a {\it flow line} along this vector field, 
$\eta : (0, \infty) \ni t \mapsto \eta(t) \in D$
starting from 
$\lim_{t \to 0} \eta(t) =: \eta(0)=x \in \partial D$ 
is defined (if exists) as the 
solution of the ordinary differential equation (ODE)
\cite{She16,MS16a}
\begin{equation}
\frac{d \eta(t)}{dt} = e^{\sqrt{-1}\{h(\eta(t))/\chi + \theta\}},
\quad t \geq 0, \quad \eta(0) = x.
\label{eqn:flow1}
\end{equation}
Let $\widetilde{D} \subsetneq \C$ be another simply connected domain and
consider a conformal map $\varphi: \widetilde{D} \to D$.
Then we define the pull-back of the flow line $\eta$ by $\varphi$ as
$\widetilde{\eta}(t)=(\varphi^{-1} \circ \eta)(t)$.
That is, 
$\varphi(\widetilde{\eta}(t))=\eta(t)$, 
and the derivatives with respect to $t$ of the both sides of this equation 
gives 
$\varphi'(\widetilde{\eta}(t)) d \widetilde{\eta}(t)/dt=d \eta(t)/dt$
with $\varphi'(z) := d \varphi(z)/dz$.
We use the polar coordinate
$\varphi'(\cdot)=|\varphi'(\cdot)| e^{\sqrt{-1} \arg \varphi'(\cdot)}$,
where $\arg \zeta$ of $\zeta \in \C$ is a priori 
defined up to additive multiples of $2 \pi$, and hence
we have 
$d \widetilde{\eta}(t)/dt=
e^{\sqrt{-1}\{ (h \circ \varphi - \chi \arg \varphi')(\widetilde{\eta}(t))/\chi + \theta\}}
/|\varphi'(\widetilde{\eta}(t))|, t \geq 0$.
If we perform a time change $t \to \tau=\tau(t)$ by putting
$t=\int_0^{\tau} ds/|\varphi'(\widetilde{\eta}(s))|$ and 
$\widehat{\eta}(t):=\widetilde{\eta}(\tau(t))$, then 
the above equation becomes
\[
\frac{d \widehat{\eta}(t)}{dt}
= e^{\sqrt{-1}\{ (h \circ \varphi - \chi \arg \varphi')(\widehat{\eta}(t))/\chi + \theta\}},
\quad t \geq 0.
\]
Since a time change does not affect the geometry of a flow line,
we can identify $h$ on $D$ and 
$h \circ \varphi - \chi \arg \varphi'$ on $\widetilde{D}=\varphi^{-1}(D)$.
In \cite{She16,MS16a,MS16b,MS16c,MS17}, such a flow line
is considered also in the case that $h$ is given by 
an instance of a GFF defined as follows.

\begin{df}
\label{thm:def_GFF}
Let $D\subsetneq \mathbb{C}$ be a simply connected 
domain and $H$ be the Dirichlet boundary GFF 
following the probability law $\P$
(constructed in Section \ref{sec:GFF}). 
A GFF on $D$ is a random distribution $h$ 
of the form $h=H+u$, where $u$ is a 
deterministic harmonic function on $D$.
\end{df} 
\noindent
Since a GFF is not function-valued, 
but it is a {\it distribution-valued random field}
(see Remark \ref{thm:Remark3_1} in Section \ref{sec:GFF}),
the ODE in the form (\ref{eqn:flow1}) 
no longer makes sense mathematically in general. 
Using the theory of SLE, however, 
the notion of flow lines has been generalized as follows.

Consider the collection
\[
\sS:=\left\{(D,h) \Bigg| 
\substack{D\subsetneq \mathbb{C}:\ 
\mbox{\footnotesize simply connected} \\ 
h:\ \mbox{\footnotesize GFF on $D$}} \right\}.
\]
Fixing a parameter $\chi \in \R$, 
we define the following equivalence relation in $\sS$.
\begin{df}
\label{thm:def_IS}
Two pairs $(D,h)$ and
$(\widetilde{D}, \widetilde{h}) \in \sS$ are 
{\it equivalent} if there exists a conformal map 
$\varphi: \widetilde{D} \to D$ and
$\widetilde{h} 
\law= h \circ \varphi - \chi \arg \varphi^{\prime}$
in $\P$. 
In this case, we write 
$(D, h) \sim (\widetilde{D}, \widetilde{h})$.
\end{df}
\noindent
We call each orbit belonging to $\sS/\sim$ 
an {\it imaginary surface} \cite{MS16a} 
(or an {\it AC surface} \cite{She16}). 
That is, in this equivalence class, 
a conformal map $\varphi$ causes not only 
a coordinate change of a GFF as $h \mapsto h \circ \varphi$
associated with changing the domain of definition of the field as
$D \mapsto \varphi^{-1}(D)$,
but also an addition of a deterministic 
harmonic function $-\chi \arg \varphi'$ to the field.
Notice that this definition includes 
one parameter $\chi \in \R$.
Then the collection of its flow lines is named as
the {\it imaginary geometry} \cite{MS16a}
(or the {\it AC geometry} \cite{She16}). 

Consider the case in which $D$ is given by the
upper half-plane $\H :=\{z \in \C: \Im z > 0 \}$
with $\partial \H=\R \cup \{\infty \}$.
Let $(B(t))_{t \geq 0}$ be a one-dimensional
standard Brownian motion starting from the origin
following the probability law $\rP$.
We consider a chordal SLE$_{\kappa}$  
driven by $(\sqrt{\kappa} B(t))_{t \geq 0}$ on
$S :=\R$ with $\kappa \in (0, 4]$ \cite{Sch00,LSW04,Law05}.
We obtain a simple curve 
(called the {\it chordal SLE$_{\kappa}$ curve}) parameterized by time 
$\eta: (0, \infty) \ni t \mapsto \eta(t) \in \H$, 
such that $\lim_{t \to 0} \eta(t) =:\eta(0)=0$, 
$\lim_{t \to \infty} \eta(t)=\infty$, and 
at each time $t > 0$, 
the chordal SLE gives a conformal map
from $\H^{\eta}_t$ to $\H$,
where $\eta(0, t] :=\{\eta(s) : s \in (0, t]\}$ 
and $\H^{\eta}_t := \H \setminus \eta(0, t]$, $t > 0$
with $\H^{\eta}_0 := \H$.
In this manuscript, we will write the chordal SLE$_{\kappa}$ as 
$(g_{\H^{\eta}_t})_{t \geq 0}$.
Let $H(\cdot)$ be an instance of the GFF on $\H$
with the Dirichlet boundary condition on $\R$
following the probability law $\P$, which 
is independent of $(\sqrt{\kappa} B(t))_{t \geq 0}$ 
and hence of $(g_{\H^{\eta}_t})_{t \geq 0}$.
Instead of $H(\cdot)$ itself, we consider the following GFF on $\H$
by adding a deterministic harmonic function,
\begin{equation}
h(\cdot) := H(\cdot) - \frac{2}{\sqrt{\kappa}} \arg (\cdot).
\label{eqn:arg1}
\end{equation}
Given $\kappa \in (0, 4]$ for the SLE$_{\kappa}$,
fix the parameter $\chi$ as
$\chi = 2/\sqrt{\kappa}-\sqrt{\kappa}/2$. 
Note that the well-known relation between $\kappa$
and the {\it central charge} $c$ of 
conformal field theory 
(see, for instance, Eq.(6) in \cite{BB02})
is simply expressed using 
the present parameter $\chi$ as $c=1-6 \chi^2$.
Let
$f_{\H^{\eta}_t} := g_{\H^{\eta}_t} - \sqrt{\kappa} B(t)
=\sigma_{-\sqrt{\kappa} B(t)} \circ g_{\H^{\eta}_t}$,
where $\sigma_s$ denotes
a shift by $s \in \R$; $\sigma_s(z)=z-s, z \in \H$. 
Then we can prove that 
\cite{SS13,She16,MS16a}
\begin{equation}
\bra h, f \ket \law= 
\Big\bra h \circ f_{\H^{\eta}_t} - \chi \arg f_{\H^{\eta}_t}', f \Big\ket \quad
\mbox{in $\P \otimes \rP$},  
\label{eqn:chi_equivalence1}
\end{equation}
$\forall f \in \cC_{\rm c}^{\infty}(\H)$, at each $t \geq 0$, 
where the {\it pairing} $\bra \cdot, \cdot \ket$ is defined
by (\ref{eqn:inner_product}) below.
See also \cite{KM13}. 
We think that this equivalence in probability law (\ref{eqn:chi_equivalence1}) 
realizes the equivalence relation defined by Definition \ref{thm:def_IS},
where conformal maps $\varphi$ are
chosen from shifts of the chordal SLE$_{\kappa}$ 
$\{f_{\H^{\eta}_t} : t \geq 0\}$.
In other words, an imaginary surface whose representative is given by
$(\H, h)$ with (\ref{eqn:arg1}) is constructed as
a pair of time-evolutionary domains, 
$f^{-1}_{\H^{\eta}_t}(\H)=\H^{\eta}_t-\sqrt{\kappa} B(t), t \geq 0$, 
and a {\it stationary process of GFF}, 
$h \circ f_{\H^{\eta}_t}-\chi \arg f_{\H^{\eta}_t}', t \geq 0$,
defined on it.
It was proved \cite{She16,MS16a,MS16b,MS16c,MS17} that
the ray of this imaginary geometry 
starting from the origin 
is realized as the chordal SLE$_{\kappa}$ curve $\eta$
when $\kappa \in (0, 4]$.
Moreover, it was argued that, if $\chi=0$
(i.e., $\kappa=4$), the flow lines are identified with
the zero contour lines of the GFF $h$ \cite{SS13}.

Notice that $\arg z$ in (\ref{eqn:arg1}) is the imaginary part of 
the complex analytic function $\log z$.
Sheffield \cite{She16} studied another type of distribution-valued 
random field on $\H$ given by
$\widetilde{h}(\cdot) := \widetilde{H}(\cdot)+ (2/\sqrt{\kappa}) 
\Re \log (\cdot)
= \widetilde{H}(\cdot)+ (2/\sqrt{\kappa}) \log |\cdot|$,
where $\widetilde{H}(\cdot)$ is an instance of the
free boundary GFF on $\H$.
An equivalence class of pairs represented by $(D, \widetilde{h})$ 
is called a {\it quantum surface},
which gives a mathematical realization of the quantum gravity
\cite{DS11}. In \cite{She16}, this
quantum surface was shown to be stationary 
under a backward SLE.

In this last section, we generalize some of the above results
to the case in which the conformal maps
are generated by a multiple Loewner equation associated with
a multi-slit.
This section is based on the collaborations with
Shinji Koshida (Chuo University)
\cite{KK19,KK20}.
See also \cite{Kos19}. 

\subsection{Multiple SLEs} \label{sec:multiple_SLE}
\subsubsection{Loewner equation for single-slit and multi-slit}
\label{sec:chordal_LE}

Let $D$ be a simply connected domain in $\C$ which does not
complete the plane; $D \subsetneq \C$.
Its boundary is denoted by $\partial D$.
We consider a slit in $D$, which is defined as a trace
$\eta = \{\eta(t) : t \in (0, \infty) \}$
of a simple curve 
$\eta(t) \in D, 0 < t < \infty$; 
$\eta(s) \not= \eta(t)$ for $s \not=t$.
We assume that the initial point of the slit is located in $\partial D$,
$^{\exists}\eta(0) := \lim_{t \to 0} \eta(t) \in \partial D$.
Let $\eta(0,t]:=\{\eta(s) : s \in (0, t]\}$ and 
$D^{\eta}_t := D \setminus \eta(0, t], t \in (0, \infty)$.
The Loewner theory \cite{Low23} describes a slit $\eta$
by encoding the curve into a time-dependent
analytic function $g_{D^{\eta}_t} : t \in (0, \infty)$
such that 
\[
g_{D^{\eta}_t} : \mbox{conformal map} \, \, 
D^{\eta}_t \to D,
\quad t \in (0, \infty).
\]

By the {\it Riemann mapping theorem} 
(see, for instance, Section 6 in \cite{Ahl79}),
for $D \subsetneq \C$ and a point $z_0 \in D$,
there exists a unique analytic function
$\varphi(z)$ in $D$, normalized by
$\varphi(z_0)=0, \varphi'(z_0) > 0$, such that
\[
\varphi : \mbox{conformal map} \, \, 
D \to \D,
\]
where $\D$ denotes a unit disk; $\D:=\{z \in \C: |z| < 1\}$.
Loewner gave differential equation for $g_{D^{\eta}_t} $ 
in the case $D=\D$, which is called the Loewner equation \cite{Low23}.
Since a special case of the M\"obius transformation
\[
\sm(z) := \sqrt{-1} \frac{\alpha-z}{\alpha+z}, \quad |\alpha|=1,
\]
maps $\D$ to the upper half plane
$\H:=\{z \in \C: \Im z > 0\}$
with $\sm(0)=\sqrt{-1}, \sm(\infty) = -\sqrt{-1}$,
we can apply the Loewner theory to the case
with $D=\H$, in which
$^{\exists} \eta(0) := \lim_{t \to 0} \eta(t) \in \R$
and $\eta(0, t] \subset \H$ for $t \in (0, \infty)$ \cite{KSS68}.
For each time  $t \in (0, \infty)$,
$\H^{\eta}_t := \H \setminus \eta(0, t]$ is a simply connected
domain in $\C$ and there exists a unique analytic function
$g_{\H^{\eta}_t}$ such that
\[
g_{\H^{\eta}_t} : \mbox{conformal map} \, \, 
\H^{\eta}_t \to \H,
\]
which satisfies the condition
\[
g_{\H^{\eta}_t}(z)=z+ \frac{c_t}{z} + \rO(|z|^{-2})
\quad \mbox{as $z \to \infty$}
\]
for some $c_t > 0$, 
in which the coefficient of $z$ is unity and no constant term appears.
This is called the
{\it hydrodynamic normalization}.
The coefficient $c_t$ gives the 
{\it half-plane capacity} of 
$\eta(0, t]$ and denoted by 
$\hcap(\eta(0, t])$.
The following has been shown (see \cite{KSS68,Law05,dMG13}).

\begin{thm}
\label{thm:LE}
Let $\eta$ be a slit in $\H$ such that
$\hcap(\eta(0, t])=2t, t \in (0, \infty)$.
Then there exists a unique continuous 
driving function $V(t) \in \R, t \in (0, \infty)$ such that the 
solution $g_t$ of the differential equation
\begin{equation}
\frac{d g_t(z)}{dt} = \frac{2}{g_t(z)-V(t)}, \quad
t \geq 0, \quad g_0(z)=z,
\label{eqn:LE1}
\end{equation}
gives $g_t=g_{\H^{\eta}_t}, t \in (0, \infty)$.
\end{thm}

The equation (\ref{eqn:LE1}) is called
the {\it chordal Loewner equation}.
Note that at each time $t \in (0, \infty)$,
the tip of slit $\eta(t)$ and the value of $V(t)$
satisfy the following relations,
\begin{equation}
V(t)=\lim_{\substack{z \to 0, \cr \eta(t)+z \in \H^{\eta}_t}}
g_{\H^{\eta}_t}(\eta(t)+z)
\quad \Longleftrightarrow \quad
\eta(t)=\lim_{\substack{z \to 0, \cr z \in \H}} g_{\H^{\eta}_t}^{-1}(V(t)+z). 
\label{eqn:Ut1}
\end{equation}
Moreover, 
$V(t)=\lim_{s < t, s \to t} g_{\H^{\eta}_s}(\eta(t))$
and $t \mapsto V(t)$ is continuous
(see, for instance, Lemma 4.2 in \cite{Law05}).
We write 
\[
g_{\H^{\eta}_t}(\eta(t))=V(t) \in \R, \quad t \geq 0
\]
in the sense of (\ref{eqn:Ut1}). 

\vskip 0.3cm
\begin{example}
\label{thm:SLE0}
When the driving function is identically zero;
$V(t) \equiv 0, t \in (0, \infty)$,
the chordal Loewner equation 
$dg_{\H^{\eta}_t}(z)/dt=2/g_{\H^{\eta}_t}(z), t \geq 0$
is solved under the initial condition
$g_{\H^{\eta}_0}(z)=z \in \H$ as
$g_{\H^{\eta}_t}(z)^2=4t + z^2, t \geq 0$.
In this simple case, (\ref{eqn:Ut1}) gives
$\eta(t)=2 \sqrt{-1} t^{1/2}, t \geq 0$.
That is, the slit $\eta(0, t], t >0$
is a straight line along the imaginary axis 
starting from the origin, 
$\eta(0)=\lim_{t \to 0} \eta(t)=0$,
and growing upward as time $t$ is passing.
\end{example}
\vskip 0.3cm
\begin{example}
\label{thm:SLE0b}
The above example can be extended by introducing
one parameter $\alpha \in (0, 1)$ as follows.
Let
$\kappa=\kappa(\alpha)=4(1-2 \alpha)^2/\{\alpha(1-\alpha)\}$,
and consider the case such that
\[
V(t)=\begin{cases}
\sqrt{\kappa t}, & \mbox{if $\alpha \leq 1/2$}, \cr
-\sqrt{\kappa t}, & \mbox{if $\alpha > 1/2$}.
\end{cases}
\]
In this case, the inverse of $g_t$ is solved as
$g_{\H^{\eta}_t}^{-1}(z) = \left( z + 2 \sqrt{\frac{\alpha}{1-\alpha}} \sqrt{t} \right)^{1-\alpha}
\left(z - 2 \sqrt{\frac{1-\alpha}{\alpha}} \sqrt{t} \right)^{\alpha}$,
and the slit is obtained as
\[
\eta(t)=g_{\H^{\eta}_t}^{-1}(V(t))
=2 \left( \frac{1-\alpha}{\alpha} \right)^{1/2-\alpha}
e^{\sqrt{-1} \alpha \pi} t^{1/2},
\quad t \geq 0.
\]
The slit grows from the origin along a straight line in $\H$
which makes an angle $\alpha \pi$ with respect to the positive
direction of the real axis.
When $\alpha=1/2$, this is reduced to the result mentioned
in Example \ref{thm:SLE0}.
More detail for this example, see Example 4.12 in \cite{Law05}
and Section 2.2 in \cite{Kat15_Springer}.
\end{example}

Theorem \ref{thm:LE} can be extended 
to the situation such that $\eta$ in $\H$ is given by
a multi-slit \cite{RS17}.
Let $N \in \N :=\{1,2, \dots\}$ and assume that we have $N$ slits 
$\eta_i =\{\eta_i(t): t \in (0, \infty)\} \subset \H$, 
$1 \leq i \leq N$,
which are simple curves, disjoint with each other, 
$\eta_i \cap \eta_j = \emptyset, i \not= j$,
starting from $N$ distinct points 
$\lim_{t \to 0} \eta_i(t) =: \eta_i(0)$ on $\R$;
$\eta_1(0) < \cdots < \eta_N(0)$, 
and all going to infinity; $\lim_{t \to \infty} \eta_i(t)=\infty$, $1 \leq i \leq N$. 
A multi-slit is defined as a union of them, 
$\bigcup_{i=1}^N \eta_i$, and 
$\H^{\eta}_t := \H \setminus \bigcup_{i=1}^N \eta_i(0, t]$
for each $t > 0$ with $\H^{\eta}_0 := \H$.
For each time $t \in (0, \infty)$, 
$\H^{\eta}_t$
is a simply connected domain in $\C$ and then
there exists a unique analytic function $g_{\H^{\eta}_t}$ such that
\[
g_{\H^{\eta}_t} : \mbox{conformal map} \, \, 
\H^{\eta}_t \to \H,
\]
satisfying the hydrodynamic normalization condition
\[
g_{\H^{\eta}_t}(z)=z+ 
\frac{\hcap(\bigcup_{i=1}^N \eta_i(0, t])}{z} + \rO(|z|^{-2})
\quad \mbox{as $z \to \infty$}.
\]
\begin{thm}[\cite{RS17}]
\label{thm:mLE}
For $N \in \N$, 
let $\bigcup_{i=1}^N \eta_i$ be a multi-slit in $\H$ such that
$\hcap(\bigcup_{i=1}^N \eta(0, t])=2t, t \in (0, \infty)$.
Then there exists a set of weight functions 
$\lambda_i(t) \geq 0, t \geq 0, 1 \leq i \leq N$
satisfying $\sum_{i=1}^N \lambda_i(t)=1, t \geq 0$ and an
$N$-variate continuous driving function 
$\V(t)=(V_1(t), \dots, V_N(t)) \in \R^N, t \in (0, \infty)$ such that the 
solution $g_t$ of the differential equation
\begin{equation}
\frac{d g_t(z)}{dt} = \sum_{i=1}^N \frac{2 \lambda_i(t)}{g_t(z)-V_i(t)}, \quad
t \geq 0, \quad 
g_0(z)=z,
\label{eqn:mLE1}
\end{equation}
gives $g_t=g_{\H^{\eta}_t}, t \in (0, \infty)$.
\end{thm}

Roth and Schleissinger \cite{RS17} called (\ref{eqn:mLE1}) 
the {\it Loewner equation for the multi-slit} 
$\bigcup_{i=1}^N \eta_i$.
Similar to (\ref{eqn:Ut1}), the following relations hold,
\begin{align}
V_i(t) =
\lim_{\substack{z \to 0, \cr \eta_i(t)+z \in \H^{\eta}_t}}
g_{\H^{\eta}_t}(\eta_i(t)+z)
\quad &\Longleftrightarrow \quad
\eta_i(t)=\lim_{\substack{z \to 0, \cr z \in \H}} g_{\H^{\eta}_t}^{-1}(V_i(t)+z),
\nonumber\\
&
\quad 1 \leq i \leq N, \quad t \geq 0,
\label{eqn:Ut2}
\end{align}
and we write $g_{\H^{\eta}_t}(\eta_i(t))=V_i(t) \in \R$,
$1 \leq i \leq N, t \geq 0$ in this sense.

The Loewner equation for the multi-slit (\ref{eqn:mLE1}) 
given for $D=\H$ can be mapped to other simply 
connected domains by conformal maps.
Here we consider a conformal transformation
\begin{equation}
\varphi(z)=\sqrt{z} : \H \to \O,
\label{eqn:root_z}
\end{equation}
where $\O$ denotes the first orthant in $\C$;
$\O:=\{z \in \C : \Re z >0, \Im z >0 \}$.
We set 
\begin{equation}
\widehat{g}_t(z) =\sqrt{g_t(z^2)+c(t)}, \quad 
t \geq 0, \quad z \in \O
\label{eqn:ghat}
\end{equation}
with a function of time $c(t), t \geq 0$.
Then we can see that (\ref{eqn:mLE1}) is transformed to the 
following,
\begin{align}
\frac{d \widehat{g}_t(z)}{dt}
&=\sum_{i=1}^N \left( 
\frac{2 \widehat{\lambda}_i(t)}{\widehat{g}_t(z)-\widehat{V}_i(t)}
+\frac{2 \widehat{\lambda}_i(t)}{\widehat{g}_t(z)+\widehat{V}_i(t)} 
\right)
+\frac{2 \widehat{\lambda}_0(t)}{\widehat{g}_t(z)}, \quad t \geq 0,
\nonumber\\
\widehat{g}_0(z) &=z \in \O,
\label{eqn:LE_O1}
\end{align}
where 
$\widehat{V}_i(t)=\sqrt{V_i(t)+c(t)}, t \geq 0, 1 \leq i \leq N$ 
and 
$2 \sum_{i=1}^N \widehat{\lambda}_i(t)+\widehat{\lambda}_0(t)=(1/4) dc(t)/dt, t \geq 0$.
Here we can assume that $\widehat{V}_i(t) \in \R_+$
without loss of generality,
since, even if we allow 
$\widehat{V}_i(t) \in \R_+ \cup \sqrt{-1} \R_+ \cup \{0\}$,
we can transform the whole system by a (possibly random)
automorphism of $\O$ to the case that
$\widehat{V}_i(t) \in \R_+$.
The equation (\ref{eqn:LE_O1}) can be regarded as the
multi-slit version of the
{\it quadrant Loewner equation} considered in \cite{Tak14}.
The solution of (\ref{eqn:LE_O1}) gives the {\it uniformization map}
to $\O$;
\[
\widehat{g}_t=g_{\O^{\eta}_t} : \mbox{conformal map} \, \, 
\O^{\eta}_t \to \O,
\]
where $\O^{\eta}_t := \O \setminus \bigcup_{i=1}^N \eta_i(0, t]$, 
and $g_{\O^{\eta}_t}(\eta_i(t))=\widehat{V}_i(t) \in \R_{\geq 0}$,
$1 \leq i \leq N, t \geq 0$. 

\subsubsection{SLE }
\label{sec:SLE}

So far we have considered the problem where, 
given time-evolution of a single slit $\eta(0, t], t \geq 0$
or a multi-slit $\sum_{i=1}^N \eta(0, t], t \geq 0$ in $\H$,
time-evolution of the conformal map from $\H^{\eta}_t$ to $\H$,
$t \geq 0$ is asked.
The answers are given by the solution of the Loewner equation
(\ref{eqn:LE1}) in Theorem \ref{thm:LE} for a single slit
and by the solution of the multiple Loewner equation
(\ref{eqn:mLE1}) in Theorem \ref{thm:mLE} for a multi-slit,
which are driven by a single-value process
$(V(t))_{t \geq 0}$ and by a multi-variate process
$\V(t)=(V_1(t), \dots, V_N(t)) \in \R^N, t \geq 0$, 
respectively. The both processes are defined in $\R$
and deterministic.

For $\H$ with a single slit, Schramm considered
an inverse problem in a probabilistic setting \cite{Sch00}.
He first asked a suitable family of driving stochastic processes
$(X(t))_{t \geq 0}$ on $\R$.
Then he asked the probability law of a random slit in $\H$,
which will be determined by the relations (\ref{eqn:Ut1})
from $(X(t))_{t \geq 0}$ and the solution 
$g_t=g_{\H^{\eta}_t}, t \geq 0$
of the Loewner equation (\ref{eqn:LE1}).
Schramm argued that conformal invariance implies that the
driving process $(X(t))_{t \geq 0}$ should be 
a continuous Markov process which has in a particular parameterization
independent increments.
Hence $X(t)$ can be a constant time change of
a one-dimensional standard Brownian motion
$(B(t))_{t \geq 0}$, and it is expressed as
$(\sqrt{\kappa} B(t))_{t \geq 0} \law= (B(\kappa t))_{t \geq 0}$
with a parameter $\kappa >0$.
The solution of the Loewner equation 
driven by $X(t)=\sqrt{\kappa} B(t), t \geq 0$, 
\begin{equation}
\frac{d g_{\H^{\eta}_t}(z)}{dt} 
= \frac{2}{g_{\H^{\eta}_t}(z)-\sqrt{\kappa} B(t)},\quad t \geq 0, 
\quad g_{\H^{\eta}_0}(z)=z \in \H,
\label{eqn:SLE1}
\end{equation}
is called the {\it chordal Schramm--Loewner evolution} ({\it chordal SLE})
with parameter $\kappa >0$ and is written as
SLE$_{\kappa}$ for short. 

The following was proved by Lawler, Schramm, and Werner \cite{LSW04}
for $\kappa=8$ and by Rohde and Schramm \cite{RS05} for $\kappa \not=8$.
\begin{prop} 
\label{thm:SLE_curve}
By (\ref{eqn:Ut1}), a chordal SLE$_{\kappa}$
$g_{\H^{\eta}_t}, t \in (0, \infty)$ determines 
a continuous curve $\eta=\{\eta(t): t \in (0, \infty)\} \subset \H$
with probability one.
\end{prop}

The continuous curve $\eta$ determined by an SLE$_{\kappa}$ 
is called an {\it SLE$_{\kappa}$ curve}.
The probability law of an SLE$_{\kappa}$ curve depends on 
$\kappa$. As a matter of fact, SLE$_{\kappa}$ curve
becomes self-intersecting and can touch the real axis $\R$
when $\kappa > 4$,
so it is no more a slit, since a slit has been defined as
a trace of a continuous simple curve.

There are three phases of an SLE$_{\kappa}$ curve
as shown by the follows. 

\begin{prop}
\label{thm:SLE_3_phases}
For each case, 
the following statements hold with probability one.
\begin{description}
\item{\rm (i)} \quad
If $0 < \kappa \leq 4$, then the SLE$_{\kappa}$ curve
is simple, $\eta=\eta(0, \infty) \subset \H$, and
$\lim_{t \to \infty} |\eta(t)|=\infty$.

\item{\rm (ii)} \quad
If $4 < \kappa < 8$,
the SLE$_{\kappa}$ curve is self-intersecting, $\eta \cap \R \not= \emptyset$,
and hence at each time $t \in (0, \infty)$ 
the hull of the SLE$_{\kappa}$ curve can be defined
as the union of $\eta(0, t]$ and the finite domain in $\H$
enclosed by any segment of $\eta(0, t]$ and the real axis $\R$, which is denoted by $K_t$. 
Then
$\bigcup_{t >0} \overline{K_t} = \overline{\H} := \H \cup \R$
and hence $|\eta(t)| \to \infty$ as $t \to \infty$, 
but $\eta(0, \infty) \cap \H \not= \H$.

\item{\rm (iii)} \quad
If $\kappa \geq 8$, then $\eta$ is a space-filling curve.
That is, if we put $\eta[0, \infty) := \{0\} \cup \eta(0, \infty)$,
then $\eta[0, \infty) = \overline{\H}$.
\end{description}
\end{prop}

For proof of Proposition \ref{thm:SLE_3_phases}
and more detailed description of
the probability laws of an SLE$_{\kappa}$ curves
at special values of $\kappa$, see, for instance
\cite{Law05,Kat15_Springer,Kem17}.

\subsubsection{Multiple SLE}
\label{sec:mSLE}

For simplicity, we assume that 
$\lambda_i(t) \equiv 1/N, t \geq 0, 1 \leq i \leq N$
in (\ref{eqn:mLE1}) in Theorem \ref{thm:mLE}.
Then by a simple time change $t/N \to t$
associated with a change of notation, $g_{Nt} \to g_t =: g_{\H^{\eta}_t}$,
the Loewner equation for the multi-slit in $\H$ is written as
\begin{equation}
\frac{d g_{\H^{\eta}_t}(z)}{dt}
= \sum_{i=1}^N \frac{2}{g_{\H^{\eta}_t}(z)-X^{\R}_i(t)}, \quad t \geq 0, \quad 
g_{\H^{\eta}_0}(z)=z \in \H.
\label{eqn:mSLE1}
\end{equation}
Then we ask what is the suitable family of driving
stochastic processes of $N$ particles on $\R$,
$\X^{\R}(t)=(X^{\R}_1(t), \dots, X^{\R}_N(t)), t \geq 0$.

The same argument with Schramm \cite{Sch00} will give that
$\X^{\R}(t)$ should be a continuous Markov process.
Moreover, Bauer, Bernard, and Kyt\"ol\"a \cite{BBK05}, 
Graham \cite{Gra07}, and Dub\'edat \cite{Dub07} argued that 
$(X^{\R}_i(t))_{t \geq 0}, 1 \leq i \leq N$ are 
semi-martingales and the quadratic variations 
should be given by 
$\langle d X^{\R}_i, dX^{\R}_j \rangle_t= \kappa \delta_{ij} dt, t \geq 0$,
$1 \leq i, j \leq N$ with $\kappa > 0$.
Then we can assume that the system of SDEs for
$(\X^{\R}(t))_{t \geq 0}$ is in the form,
\begin{equation}
d X^{\R}_i(t)=\sqrt{\kappa} dB_i(t)+
F^{\R}_i(\X^{\R}(t)) dt, \quad t \geq 0,
\quad 1 \leq i \leq N,
\label{eqn:SDE_A1}
\end{equation}
where $(B_i(t))_{t \geq 0}, 1 \leq i \leq N$ are
independent one-dimensional standard Brownian motions,
$\kappa > 0$, and
$\{F^{\R}_i(\x)\}_{i=1}^N$ are suitable functions
of $\x=(x_1, \dots, x_N)$ which do not explicitly
depend on $t$. 

In the orthant system (\ref{eqn:LE_O1}),
we put $\widehat{\lambda}_i(t) \equiv r/(2N)$, 
$t \geq 0, r >0, 1 \leq i \leq N$
and $dc(t)/dt=4, t \geq 0$, 
and perform a time change
$r t/(2N) \to t$ associated with a change of
notation $\widehat{g}_{2N t/r} \to \widehat{g}_t =: g_{\O^{\eta}_t}$.
Then the Loewner equation in $\O$ is written as
\begin{align}
\frac{d g_{\O^{\eta}_t}(z)}{dt}
&= \sum_{i=1}^N \left( 
\frac{2}{g_{\O^{\eta}_t}(z)-X^{\R_{\geq 0}}_i(t)} 
+ \frac{2}{g_{\O^{\eta}_t}(z)+X^{\R_{\geq 0}}_i(t)} \right) 
+ \frac{4 \delta}{g_{\O^{\eta}_t}(z)}, \quad t \geq 0,
\nonumber\\
g_{\O^{\eta}_0}(z) &=z \in \O.
\label{eqn:mSLE2}
\end{align}
where $\delta := N(1-r)/r \in \R$.
We assume that the system of SDEs for
$\X^{\R_{\geq 0}}(t) \in (\R_{\geq 0})^N, t \geq 0$ 
is in the same form as (\ref{eqn:SDE_A1}),
\begin{equation}
d X^{\R_{\geq 0}}_i(t) =\sqrt{\kappa} d \widetilde{B}_i(t)+
F^{\R_{\geq 0}}_i( \X^{\R_{\geq 0}}(t)) dt, 
\quad t \geq 0,
\quad 1 \leq i \leq N,
\label{eqn:SDE_B1}
\end{equation}
with a set of independent one-dimensional
standard Brownian motions $(\widetilde{B}_i(t))_{t \geq 0}$, 
$1 \leq i \leq N$.

\subsection{GFF with Dirichlet boundary condition } \label{sec:GFF}
\subsubsection{Bochner--Minlos Theorem} 
\label{sec:GFF_Dirichlet}

Here we start with the classical {\it Bochner theorem},
which states that
a probability measure on a finite dimensional Euclidean space
is determined by a 
{\it characteristic function}
which is a Fourier transform of
the probability measure. 
Note that we have considered Laplace transforms
of probability measures in Section \ref{sec:DPP}
and multitime Laplace transforms
of probability measures in Section \ref{sec:SLG}.
First we define a functional of positive type.

\begin{df}
\label{thm:positive_type}
Let $\cV$ be a finite or infinite dimensional vector space.
A function $\psi: \cV \to \C$ is said to be 
a functional of positive type if for arbitrary $N \in \N$,
$\xi_1, \dots, \xi_N \in \cV$, and $z_1, \dots, z_N \in \C$, we have
\[
\sum_{n=1}^N \sum_{m=1}^N
\psi(\xi_n - \xi_m) z_n \overline{z_m} \geq 0.
\]
\end{df}
Then the following is proved.
\begin{lem}
\label{thm:positive_type2}
Let $\psi: \cV \to \C$ be a functional of positive type
on a vector space $\cV$. Then it follows that
{\rm (i)} 
$\psi(0) \geq 0$, 
{\rm (ii)} 
$\psi(\xi)=\overline{\psi(-\xi)}$ for all $\xi \in \cV$,
and
{\rm (iii)} 
$|\psi(\xi)| \leq \psi(0)$ for all $\xi \in \cV$.
\end{lem}
For $x, y \in \R^N$, the standard inner product is denoted by
$x \cdot y$ and we write $|x| := \sqrt{x \cdot x}$.
Let $\cB^N$ be the family of Borel sets in $\R^N$.
Then the following is known as the Bochner theorem.

\begin{thm}[Bochner theorem]
\label{thm:Bochner}
Let $\psi: \R^N \to \C$ be a continuous functional
of positive type such that $\psi(0)=1$. 
Then there exists a unique probability measure
$\rP$ on $(\R^N, \cB^N)$ such that
\[
\psi(\xi)=\int_{\R^N} e^{\sqrt{-1} x \cdot \xi} \, \rP(d x)
\quad \mbox{for \, $\xi \in \R^N$}.
\]
\end{thm}

If we consider the case that $\psi(\xi)$ 
is given by 
$\Psi(\xi):= e^{-|\xi|^2/2}, \xi \in \R^N$,
then 
the probability measure $\rP$ given by the Bochner theorem
is the finite-dimensional standard Gaussian measure,
\begin{align*}
\rP(d x) &=\frac{1}{(2\pi)^{N/2}} e^{-|x|^2/2} d x
\nonumber\\
&=\prod_{i=1}^N \lambda_{\rN(0,1)}(dx_i),
\quad x=(x_1, \dots, x_N) \in \R^N.
\end{align*}
Hence we can say that
the finite-dimensional standard Gaussian measure 
$\rP$ is determined by
the characteristic function $\Psi(\xi)$ as
\begin{align*}
\Psi(\xi) &=\int_{\R^N} e^{\sqrt{-1} x \cdot \xi} \rP(d x)
\nonumber\\
&= e^{-|\xi|^2/2} \quad \mbox{for \, $\xi \in \R^N$}.
\end{align*}

Now consider the case that $\cH$ is an infinite dimensional Hilbert space with inner product 
$\bra \cdot, \cdot \ket=\bra \cdot, \cdot \ket_{\cH}$
with $\| x \|=\| x \|_{\cH}=\sqrt{\bra x,x \ket_{\cH}}$,
$x \in \cH$.
The dual space of $\cH$ will be denoted by
$\cH^{\ast}$. 
Suppose that there were a probability measure $\rP$
on $\cH$ with a suitable $\sigma$-algebra such that
\[
\psi(\xi) 
= \int_{\cH} e^{\sqrt{-1} \bra x, \xi \ket} \rP(d x)
= e^{-\|\xi\|^2/2}
\quad \mbox{for $\xi \in \cH$}.
\]
Let $\{e_n\}_{n=1}^{\infty}$ be a complete orthonormal system (CONS)
of $\cH$.
If we set $\xi=t e_n, t \in \R$
for an arbitrary $n \in \N$, then
\[
\int_{\cH} e^{\sqrt{-1} t \bra x, e_n \ket} \rP(d x) = e^{-t^2/2},
\quad t \in \R.
\]
Since $x \in \cH$, we have $\bra x, e_n \ket \to 0$ as $n \to \infty$.
Therefore in the limit $n \to \infty$, 
the above equality gives $e^{-t^2/2}=1$, which is a contradiction.
This observation suggests that the application of the Bochner theorem
to an infinite dimensional space requires more consideration. 
The following arguments are base on \cite{Asa10}
and a note given by Koshida \cite{Koshida_note19}.

Let $D \subsetneq \C$ be a simply connected domain
that is bounded.
We consider the case $\cH=L^2(D, \mu(dz))$ with 
$\bra f, g \ket := \int_D f(z) g (z) d \mu(z)$, $f, g \in L^2(D, \mu(dz))$,
where $\mu(dz)$ is the Lebesgue measure on $D \subset \C$;
$\mu(dz)=dz d\overline{z}$.
Let $\Delta$ be the Dirichlet Laplacian acting on $L^2(D, \mu(dz))$.
Then $-\Delta$ has positive discrete eigenvalues so that
\begin{equation}
-\Delta e_n = \lambda_n e_n, \quad
e_n \in L^2(D, \mu(dz)), \quad n \in \N.
\label{eqn:lambda}
\end{equation}
We assume that the eigenvalues are labeled 
in a non-decreasing order;
$0 < \lambda_1 \leq \lambda_2 \leq \cdots$.
The system of eigenvalue functions $\{e_n\}_{n \in \N}$ forms
a CONS of $L^2(D)$. 
The following is known as the
{\it Weyl formula}
\begin{lem}
\label{thm:Weyl}
Let $D \subsetneq \C$ be a simply connected finite domain.
The eigenvalues $\{\lambda_n\}_{n \in \N}$ 
of the operator $-\Delta$ on $D$
exhibit the
following asymptotic behavior, 
\[
\lim_{n \to \infty} \frac{\lambda_n}{n} =\rO(1).
\]
\end{lem}

For two functions $f, g \in \cC_{\rm c}^{\infty}(D)$,
their {\it Dirichlet inner product} is defined as
\begin{equation}
\bra f, g \ket_{\nabla} 
:= \frac{1}{2 \pi} \int_D (\nabla f)(z) \cdot (\nabla g)(z)
\mu(d z).
\label{eqn:Dirichlet_IP}
\end{equation}
The Hilbert space completion of $\cC_{\rm c}^{\infty}(D)$ with respect
to this Dirichlet inner product
will be denoted by $W(D)$.
We write 
$\|f\|_{\nabla}=\sqrt{\bra f, f \ket_{\nabla}}, f \in W(D)$.
If we set $u_n=\sqrt{2\pi/\lambda_n} \, e_n, n \in \N$,
then by integration by parts, we have
\[
\bra u_n, u_n \ket_{\nabla}=
\frac{1}{2 \pi} \bra u_n, (-\Delta) u_m \ket=\delta_{nm},
\quad n, m \in \N.
\]
Therefore $\{u_n\}_{n \in \N}$ forms a CONS of $W(D)$.

Let $\widehat{\cH}(D)$ be the space of formal real infinite
series in $\{u_n\}_{n \in \N}$.
This is obviously isomorphic to $\R^{\N}$
by setting 
$\widehat{\cH}(D) \ni \sum_{n \in \N} f_n u_n 
\mapsto (f_n)_{n \in \N} \in \R^{\N}$.
As a subspace of $\widehat{\cH}(D)$, $W(D)$ is isomorphic to
$\ell^2(\N) \subset \R^{\N}$.
For two formal series 
$f=\sum_{n \in \N} f_n u_n$,
$g=\sum_{n \in \N} g_n u_n \in \widehat{\cH}(D)$
such that
$\sum_{n \in \N} |f_n g_n| < \infty$,
we define their {\it pairing} as
$\bra f, g \ket_{\nabla}:=\sum_{n \in \N} f_n g_n$.
In the case when $f, g \in W(D)$, their pairing of course
coincides with the Dirichlet inner product 
(\ref{eqn:Dirichlet_IP}).

Notice that, for any $a \in \R$, the operator
$(-\Delta)^a$ acts on $\widehat{\cH}(D)$ as
\[
(-\Delta)^a \sum_{n \in \N} f_n u_n
:= \sum_{n \in \N} \lambda_n^a f_n u_n,
\quad (f_n)_{n \in \N} \in \R^{\N}.
\]
Using this fact, we define
$\cH_a(D) := (-\Delta)^a W(D)$, $a \in \R$,
each of which is a Hilbert space with inner product
\[
\langle f, g \rangle_a
:=\bra(-\Delta)^{-a}f, (-\Delta)^{-a}g \ket_{\nabla},
\quad f, g \in \cH_a(D).
\]
We write $\|\cdot\|_a:=\sqrt{\langle \cdot, \cdot \rangle_a}, a \in \R$.

\vskip 0.3cm
\begin{example}
\label{thm:a=1/2}
When $a=1/2$, we have
\[
\langle f, g \rangle_{1/2}
=\Big\bra (-\Delta)^{-1/2} f, (-\Delta)^{-1/2}g \Big\ket_{\nabla}
=\frac{1}{2 \pi} \bra f, g \ket, \quad
f, g \in \cH_{1/2}(D).
\]
Therefore $\cH_{1/2}(D) = L^2(D, \mu(dz))$.
\end{example}


We can prove the following two lemmas.
\begin{lem}
\label{thm:Ha<Hb}
Assume $a< b$. Then
$\cH_a(D) \subset \cH_b(D)$.
\end{lem}
\noindent{\it Proof} \,
Let $f=\sum_{n \in \N} f_n u_n \in \widehat{\cH}(D)$
be a formal series.
Then we have
\[
\|f\|_b^2 
=\sum_{n \in \N} \lambda_n^{-2b} f_n^2
\leq \sum_{n=1}^{N-1}
(\lambda_n^{-2b}-\lambda_n^{-2a}) f_n^2 + \|f\|_a^2.
\]
Since the Weyl formula (Lemma \ref{thm:Weyl}) holds, 
we can take $N \in \N$ such that $\lambda_N > 1$.
Then the desired inclusion follows \qed.
\vskip 0.3cm
\begin{lem}
\label{thm:dual}
Let $a \in \R$ and fix $h \in \cH_a(D)$.
Then the assignment
\[
\bra h, \cdot \ket_{\nabla} : \cH_{-a}(D) \to \R
\quad \mbox{such that} \quad
\cH_{-a}(D) \ni f \mapsto \bra h, f \ket_{\nabla} \in \R
\]
is well-defined and continuous. 
In particular, $\cH_a(D)$ and $\cH_{-a}(D)$
makes a dual pair of Hilbert spaces with respect to 
the Dirichlet inner product $\bra \cdot, \cdot \ket_{\nabla}$. 
\end{lem}
\noindent{\it Proof} \,
For $h=\sum_{n \in \N} h_n u_n$
and $f=\sum_{n \in \N} f_n u_n$, 
Cauchy's inequality 
\[
\sum_{n \in \N} |h_n f_n|
=\sum_{n \in \N} | (\lambda_n^{-a} h_n) (\lambda^a f_n)|
\leq \left( \sum_{n \in \N} |\lambda_n^{-a} h_n|^2 \right)^{1/2}
\left( \sum_{n \in \N} |\lambda_n^{a} f_n|^2 \right)^{1/2}
\]
ensures that the pairing 
$\bra h,f \ket$ is well-defined.
Notice that
\[
\bra h, f \ket_{\nabla}
=\bra (-\Delta)^{-a} h, (-\Delta)^a f \ket_{\nabla}
=\bra (-\Delta)^{-2a}h, f \ket_{-a}.
\]
Since $(-\Delta)^{-2a} h \in \cH_{-a}(D)$ by
the assumption $h \in \cH_a(D)$ and Lemma \ref{thm:Ha<Hb}, 
then $\bra h, \cdot \ket_{\nabla}$ is continuous on $\cH_{-a}(D)$.
Therefore $\cH_{a}(D) \simeq \cH_{-a}(D)^{\ast}$.
\qed
\vskip 0.3cm

\vskip 0.3cm
\begin{rem}
\label{thm:Remark3_1}
Since $\cH_{1/2}(D)=L^2(D, \nu(dz))$ as mentioned
in Example \ref{thm:a=1/2}, 
the members of $\cH_a(D)$
with $a > 1/2$ cannot be functions, but are distributions.
\end{rem}
\vskip 0.3cm

Define 
\begin{equation}
\cE(D) := \bigcup_{a > 1/2} \cH_a(D).
\label{eqn:cE}
\end{equation}
Then its dual Hilbert space is identified with
$\cE(D)^{\ast} :=\bigcap_{a < -1/2} \cH_a(D)$ 
by Lemma \ref{thm:dual}, and 
\[
\cE(D)^{\ast} \subset W(D) \subset \cE(D)
\]
is established (by definition and Lemma \ref{thm:Ha<Hb}).
Here $(\cE(D)^{\ast}, W(D), \cE(D))$ is called a {\it Gel'fand triple}.
We set 
$\Sigma_{\cE(D)} 
=\sigma(\{\bra \cdot, f \ket_{\nabla} : f \in \cE(D)^{\ast} \})$.
On such a setting, the following is
obtained.
This theorem is the extension of the Bochner theorem 
(Theorem \ref{thm:Bochner})
and is called 
the {\it Bochner--Minlos theorem}
(see, for instance, \cite{Hid80,She07,Asa10}). 

\begin{thm}[Bochner--Minlos theorem]
\label{thm:BM}
Let $\psi$ be a continuous function of positive type
on $W(D)$ such that $\psi(0)=1$.
Then there exists a unique probability measure 
$\bP$ on $(\cE(D), \Sigma_{\cE(D)})$
such that
\begin{equation}
\psi(f)=\int_{\cE(D)} e^{\sqrt{-1} \bra h, f \ket_{\nabla}} \bP(d h)
\quad \mbox{for $f \in \cE(D)^{\ast}$}.
\label{eqn:BMeq1}
\end{equation}
\end{thm}
We will give the proof in Section \ref{sec:A} below.

Under certain conditions on $\psi$,
the domain of function $f$ for (\ref{eqn:BMeq1})
can be extended from $\cE(D)^{\ast}$ to $W(D)$
(see Proposition \ref{thm:extension1} in Section \ref{sec:B} below).
It is easy to verify that the functional
$\Psi(f):= e^{-\|f\|_{\nabla}^2/2}$ satisfies the conditions.
Then the following is established
with a probability measure $\P$ on $(\cE(D), \Sigma_{\cE(D)})$,
\begin{align}
\Psi(f) &=\int_{\cE(D)} e^{\sqrt{-1} \bra h, f \ket_{\nabla}} \P(d h)
\nonumber\\
&=e^{-\|f\|_{\nabla}^2/2}
\quad \mbox{for $f \in W(D)$}.
\label{eqn:GFF}
\end{align}

\begin{df}[Dirichlet boundary GFF]
\label{thm:GFF_Dirichlet}
A Gaussian free field (GFF) with Dirichlet boundary condition is defined 
as a pair $((\Omega, P),H)$ of a probability space $(\Omega, P)$ 
and an isotopy 
$H : W(D) \to L^{2}(\Omega, P)$ 
such that each $H(f)$, $f\in W(D)$ is a Gaussian random variable.
\end{df}

For each $f \in W(D)$, we write 
$\bra H, f \ket_{\nabla} \in L^{2}(\cE(D), \P)$ for the random variable defined by 
$h \mapsto \bra h, f \ket_{\nabla}$, $h \in \cE(D)$.
Then (\ref{eqn:GFF}) ensures that the pair of 
$((\cE(D), \P), H)$ gives a GFF with Dirichlet boundary condition. 
We often just call $H$ a Dirichlet boundary GFF without referring 
to the probability space $(\cE(D), \P)$.

\subsubsection{Conformal invariance of GFF} 
\label{sec:conformal_invariance}

Assume that $D, D' \subsetneq \C$ are simply connected domains
and let $\varphi: D' \to D$ be a conformal map.

\begin{lem}
\label{thm:conformal_inv}
The Dirichlet inner product (\ref{eqn:Dirichlet_IP}) is 
conformally invariant, 
that is,
\[
\int_{D} (\nabla f)(z) \cdot (\nabla g)(z) \mu(d z)
=\int_{D'} (\nabla (f \circ \varphi))(z) \cdot 
(\nabla (g \circ \varphi))(z) \mu(d z)
\quad \mbox{for $f, g \in \cC_{\rm c}^{\infty}(D)$}.
\]
\end{lem}
\noindent{\it Proof} \,
For $z \in D$ we write $z=x+\sqrt{-1} y, x, y \in \R$
and put $\varphi(z)=u(x, y)+\sqrt{-1} v(x,y)$
with real-valued functions $u$ and $v$.
Owing to the Cauchy-Riemann identities,
$\partial u/\partial x=\partial v/\partial y$,
$\partial u/\partial y=-\partial v/\partial x$,
the Jacobian for the transformation $\varphi$
is written as
\[
\frac{\partial (u,v)}{\partial (x,y)}
=\frac{\partial u}{\partial x} \frac{\partial v}{\partial y}
-\frac{\partial u}{\partial y} \frac{\partial v}{\partial x}
=\left(\frac{\partial u}{\partial x}\right)^2
+\left( \frac{\partial u}{\partial y}\right)^2.
\]
From the chain-rule and the Cauchy-Riemann identities again,
we have the equality
\[
(\nabla f \circ \varphi)(z) \cdot (\nabla g \circ \varphi)(z)
=\left( \frac{\partial f}{\partial u} \frac{\partial g}{\partial u}
+\frac{\partial f}{\partial v} \frac{\partial g}{\partial v} \right)
\left\{ \left(\frac{\partial u}{\partial x}\right)^2
+\left( \frac{\partial u}{\partial y}\right)^2 \right\}.
\]
Therefore, the statement is proved. \qed
\vskip 0.3cm

From the above lemma, we see that
$\varphi^{\ast} : W(D) \ni f \mapsto f \circ \varphi \in W(D')$
is an isomorphism. This allows one to consider
a GFF on an unbounded domain.
Namely, if $D'$ is bounded on which a Dirichlet GFF $H$ is defined, 
but $D$ is unbounded, 
we can define a family 
$\{\bra \varphi_{\ast} H, f \ket_{\nabla} : f \in W(D)\}$ by
$\bra \varphi_{\ast} H, f\ket_{\nabla} 
:= \bra H, \varphi^{\ast} f \ket_{\nabla}, f \in W(D)$
so as to have the following {\it covariance structure},
\begin{equation}
\E \Big[
\bra \varphi_{\ast} H, f \ket_{\nabla} 
\bra \varphi_{\ast} H, g \ket_{\nabla} \Big]
=\bra \varphi^{\ast} f, \varphi^{\ast} g \ket_{\nabla}
=\bra f, g \ket_{\nabla} \quad
\mbox{for $f, g \in W(D)$}.
\label{eqn:cov_str}
\end{equation}
Relying on the following formal computation
\begin{align*}
\bra \varphi_{\ast} H, f \ket_{\nabla}
=\bra H, \varphi^{\ast} f \ket_{\nabla}
&=\frac{1}{2 \pi} \int_{D'} (\nabla H)(z) 
\cdot (\nabla f \circ \varphi)(z) \mu(d z)
\nonumber\\
&=\frac{1}{2 \pi} \int_D (\nabla H \circ \varphi^{-1})(z) 
\cdot (\nabla f)(z) \mu(d z)
\end{align*}
we understand the equality $\varphi_{\ast} H = H \circ \phi^{-1}$.
By the fact (\ref{eqn:cov_str}) such that
the covariance structure does not change under 
a conformal map $\phi$, we say that
{\it the GFF is conformally invariant}. 

\subsubsection{The Green's function of GFF} 
\label{sec:Green}

Assume that $D \subsetneq \C$ is a simply connected domain.
In the previous subsections, we have constructed a family
$\{\bra H, f \ket_{\nabla} : f \in W(D)\}$ of random variables
whose covariance structure is given by
\[
\E \Big[
\bra H, f \ket_{\nabla} \bra H, g \ket_{\nabla} \Big]
=\bra f, g \ket_{\nabla} \quad
\mbox{for $f, g \in W(D)$}.
\]

By a formal integration by parts, we see that
\begin{align*}
\bra H, f\ket_{\nabla} 
&= \frac{1}{2\pi} \int_{D} (\nabla H)(z) \cdot (\nabla f)(z) \mu(d z)
=\frac{1}{2\pi} \int_{D} H(z) (-\Delta f)(z) \mu(d z)
\nonumber\\
&=\frac{1}{2 \pi} \bra H, (-\Delta) f \ket.
\end{align*}
Motivated by this observation, we define
\begin{equation}
\bra H, f \ket:=2 \pi \bra H, (-\Delta)^{-1} f \ket_{\nabla} 
\quad \mbox{for $ \in \sD((-\Delta)^{-1})$}, 
\label{eqn:inner_product}
\end{equation}
where
$\sD((-\Delta)^{-1})$ denotes the domain of
$(-\Delta)^{-1}$ in $W(D)$.
Note that if $D$ is bounded, then $(-\Delta)^{-1}$ is
a bounded operator, but if $D$ is unbounded,
then $(-\Delta)^{-1}$ is not defined on $W(D)$.
The action of $(-\Delta)^{-1}$ is expressed as
an integral operator and the integral kernel
is known as {\it the Green's function}. Namely,
\[
((-\Delta)^{-1} f)(z)
=\frac{1}{2 \pi} \int_{D} G_D(z,w) f(w) \mu(d w),
\quad \mbox{a.e.} \, z \in D, \quad
f \in \sD((-\Delta)^{-1}),
\]
where $G_D(z,w)$ denotes the Green's function of $D$
under the Dirichlet boundary condition.
Hence the covariance of $\bra H, f \ket$ and $\bra H, g \ket$ with
$f, g \in \sD((-\Delta)^{-1})$ is written as
\begin{equation}
\E[ \bra H, f \ket \bra H, g \ket] 
=\int_{D \times D} f(z) G_D(z, w) g(w) \mu(d z) \mu(d w).
\label{eqn:Green}
\end{equation}
When we symbolically write
\[
\bra H, f \ket=\int_{D} H(z) f(z) \mu(d z),
\quad f \in \sD((-\Delta)^{-1}),
\]
the covariance structure can be expressed as
\[
\E[H(z) H(w)]= G_D(z, w), \quad
z, w \in D, \quad n \not=w.
\]
The conformal invariance of GFF implies that
for a conformal map
$\varphi: D' \to D$, we have the equality,
\begin{equation}
G_{D'}(z, w)=G_D(\varphi(z), \varphi(w)), \quad z, w \in D'.
\label{eqn:G_conformal}
\end{equation}

\vskip 0.3cm
\begin{example}
\label{thm:G_H}
When $D$ is the upper half plane $\H$,
\begin{align*}
G_{\H}(z,w) &= \log \left|
\frac{z-\overline{w}}{z-w} \right|
=\log|z-\overline{w}| -\log|z-w|
\nonumber\\
&= \Re \log(z-\overline{w}) - \Re \log(z-w), 
\end{align*}
$z, w \in \H, z \not= w$.
\end{example}
\vskip 0.3cm
\begin{example}
\label{thm:G_O}
When $D$ is the first orthant $\O$,
\begin{align*}
G_{\O}(z,w) &= \log \left|
\frac{(z-\overline{w})(z+\overline{w})}{(z-w)(z+w)} \right|
\nonumber\\
&=\log|z-\overline{w}| + \log|z+\overline{w}|
-\log|z-w| - \log|z+w|, 
\nonumber\\
&=\Re \log (z-\overline{w}) + \Re \log(z+\overline{w})
- \Re \log (z-w) - \Re \log(z+w), 
\end{align*}
$z, w \in \O, z \not= w$.
\end{example}
\vskip 0.3cm

From the formula (\ref{eqn:Green}), we see that
$\cC_{\rm c}^{\infty}(D) \subset \sD((-\Delta)^{-1})$.
In the following, we will consider the
family of random variables $\{\bra H, f \ket : f \in \cC_{\rm c}^{\infty}(D) \}$
to characterize the GFF, $H$.

\subsection{GFFs coupled with stochastic log-gases} \label{sec:GFF_BPs}
\subsubsection{Dirichlet boundary GFF transformed by multiple SLE}
\label{sec:GFF_SLE}

Here we write the GFF with free boundary condition defined
on a simply connected domain $D \subsetneq \C$ as $H_D$.
Consider the transformation of $H_D$ by the multiple SLE,
\[
H_{D^{\eta}_t} := H_D \circ g_{D^{\eta}_t}, \quad t \geq 0.
\]
By (\ref{eqn:G_conformal}), the Green's function of $H_{D^{\eta}_t}, t \geq 0$
is given by
\[
G_{D^{\eta}_t}(z, w) = G_D(g_{D^{\eta}_t}(z), g_{D^{\eta}_t}(w)),
\quad z, w \in D^{\eta}_t := D \Big\backslash \bigcup_{i=1}^N \eta_i(0, t],
\quad t \geq 0.
\]
Using the explicit expressions of the Greens' functions
for $D=\H$ and $\O$, the following is obtained.

\begin{lem}
\label{thm:dG}
For $D=\H$ and $\O$, the increments of $G_{D^{\eta}_t}$ in time 
$t \geq 0$ are given as
\begin{align*}
d G_{\H^{\eta}_t}(z, w)
&=- \sum_{i=1}^N \Im \frac{2}{g_{\H^{\eta}_t}(z)-X^{\R}_i(t)}
\Im \frac{2}{g_{\H^{\eta}_t}(w)-X^{\R}_i(t)} dt,
\quad z, w \in \H^{\eta}_t, \quad t \geq 0,
\nonumber\\
d G_{\O^{\eta}_t}(z, w)
&=- \sum_{i=1}^N \Im 
\left( \frac{2}{g_{\O^{\eta}_t}(z)-X^{\R_{\geq 0}}_i(t)}
- \frac{2}{g_{\O^{\eta}_t}(z)+X^{\R_{\geq 0}}_i(t)} \right)
\nonumber\\
& \qquad \quad \times
\Im \left( \frac{2}{g_{\O^{\eta}_t}(w)-X^{\R_{\geq 0}}_i(t)}
- \frac{2}{g_{\O^{\eta}_t}(w)+X^{\R_{\geq 0}}_i(t)} \right) dt, 
\quad z, w \in \O^{\eta}_t, \, \, t \geq 0.
\end{align*}
\end{lem}
\noindent{\it Proof} \,
From the explicit expressions of
$G_{\H}$ and $G_{\O}$ given in 
Example \ref{thm:G_H} and \ref{thm:G_O}, 
we have
\begin{align*}
d G_{\H^{\eta}_t}(z, w)
&= \Re \frac{ d g_{\H^{\eta}_t}(z) - d \overline{g_{\H^{\eta}_t}(w)} }
{g_{\H^{\eta}_t}(z)- \overline{g_{\H^{\eta}_t}(w)} }
- \Re \frac{ d g_{\H^{\eta}_t}(z) - d g_{\H^{\eta}_t}(w) }
{g_{\H^{\eta}_t}(z)- g_{\H^{\eta}_t}(w) }, 
\nonumber\\
d G_{\O^{\eta}_t}(z, w)
&= \Re \frac{ d g_{\O^{\eta}_t}(z) - d \overline{g_{\O^{\eta}_t}(w)} }
{g_{\O^{\eta}_t}(z)- \overline{g_{\O^{\eta}_t}(w)} }
- \Re \frac{ d g_{\O^{\eta}_t}(z) - d g_{\O^{\eta}_t}(w) }
{g_{\O^{\eta}_t}(z)- g_{\O^{\eta}_t}(w) }
\nonumber\\
& + \Re \frac{ d g_{\O^{\eta}_t}(z) + d \overline{g_{\O^{\eta}_t}(w)} }
{g_{\O^{\eta}_t}(z)+ \overline{g_{\O^{\eta}_t}(w)} }
- \Re \frac{ d g_{\O^{\eta}_t}(z) + d g_{\O^{\eta}_t}(w) }
{g_{\O^{\eta}_t}(z) + g_{\O^{\eta}_t}(w) }. 
\end{align*}
By the multiple Loewner equation (\ref{eqn:mSLE1}) in $\H$, we see that
\begin{align*}
d g_{\H^{\eta}_t}(z)- d \overline{ g_{\H^{\eta}_t}(w)}
&= \sum_{i=1}^N \frac{2 dt}{g_{\H^{\eta}_t}(z)-X^{\R}_i(t)}
- \sum_{i=1}^N \frac{2 dt}{ \overline{g_{\H^{\eta}_t}(w)}-X^{\R}_i(t)}
\nonumber\\
&= - (g_{\H^{\eta}_t}(z)- \overline{ g_{\H^{\eta}_t}(w) }
\sum_{i=1}^N \frac{2 dt}{g_{\H^{\eta}_t}(z)-X^{\R}_i(t))
(\overline{g_{\H^{\eta}_t}(w)} -X^{\R}_i(t)) }. 
\end{align*}
Hence we have
\begin{align*}
d G_{\H^{\eta}_t}(z, w)
&= \Re \sum_{i=1}^N \frac{2 dt}{g_{\H^{\eta}_t}(z)-X^{\R}_i(t))
(\overline{g_{\H^{\eta}_t}(w)} -X^{\R}_i(t)) }
\nonumber\\
& - \Re \sum_{i=1}^N \frac{2 dt}{g_{\H^{\eta}_t}(z)-X^{\R}_i(t))
(g_{\H^{\eta}_t}(w)-X^{\R}_i(t)) }.
\end{align*}
For any two complex variables $\zeta$ and $\omega$, 
it is easy to verify the equality
$\Re \zeta \overline{\omega}-\Re \zeta \omega
=2 \Im \zeta \Im \omega$, 
and then the result is obtained.
Similarly by the multiple Loewner equation (\ref{eqn:mSLE2}) in $\O$, we have
\begin{align*}
\Re \frac{ d g_{\O^{\eta}_t}(z) - d \overline{g_{\O^{\eta}_t}(w)} }
{g_{\O^{\eta}_t}(z)- \overline{g_{\O^{\eta}_t}(w)} }
&= - \sum_{i=1}^N \Re 
\frac{2 dt}{(g_{\O^{\eta}_t}(z)-X^{\R_{\geq 0}}_i(t))
(\overline{g_{\O^{\eta}_t}(w)}-X^{\R_{\geq 0}}_i(t))}
\nonumber\\
& \quad - \sum_{i=1}^N \Re 
\frac{2 dt}{(g_{\O^{\eta}_t}(z)+X^{\R_{\geq 0}}_i(t))
(\overline{g_{\O^{\eta}_t}(w)}+X^{\R_{\geq 0}}_i(t))}
-\Re \frac{ 4 \delta dt}{g_{\O^{\eta}_t}(z) \overline{g_{\O^{\eta}_t}(w)}},
\nonumber\\
\Re \frac{ d g_{\O^{\eta}_t}(z) + d \overline{g_{\O^{\eta}_t}(w)} }
{g_{\O^{\eta}_t}(z)+ \overline{g_{\O^{\eta}_t}(w)} }
&= \sum_{i=1}^N \Re 
\frac{2 dt}{(g_{\O^{\eta}_t}(z)-X^{\R_{\geq 0}}_i(t))
(\overline{g_{\O^{\eta}_t}(w)}+X^{\R_{\geq 0}}_i(t))}
\nonumber\\
& +\sum_{i=1}^N \Re 
\frac{2 dt}{(g_{\O^{\eta}_t}(z)+X^{\R_{\geq 0}}_i(t))
(\overline{g_{\O^{\eta}_t}(w)}-X^{\R_{\geq 0}}_i(t))}
+\Re \frac{ 4 \delta dt}{g_{\O^{\eta}_t}(z) \overline{g_{\O^{\eta}_t}(w)}}.
\end{align*}
Again we use the equality $\Re \zeta \overline{\omega}-\Re \zeta \omega
=2 \Im \zeta \Im \omega$ for $\zeta, \omega \in \C$, and 
then we prove the lemma. \qed
\vskip 1cm

\subsubsection{Complex-valued logarithmic potentials and martingales}
\label{sec:martingales}

We have remarked in Section \ref{sec:log_gases2} that
the Dyson model and the Bru--Wishart processes
studied in random matrix theory can be regarded
as stochastic log-gasses defined on a line $S=\R$
and a half-line $S=\R_{\geq 0}$, respectively.
There the logarithmic potential are given by (\ref{eqn:potential}). 
Here we consider a {\it complex-valued logarithmic potential}
acting between a point $z$ in a two-dimensional domain
$D \subsetneq \C$ and $N$ points 
$\x=(x_1, \dots, x_N)$ on a part of the boundary $S \subset \partial D$.
For $D=\H$ and $\O$, we put
\begin{align}
\Phi_{\H}(z, \x) &= \sum_{i=1}^N \log(z-x_i),
\quad z \in \H, \quad \x \in \R^N, 
\nonumber\\
\Phi_{\O}(z, \x) &=
\Phi_{\O}(z, \x; q)
\nonumber\\
&=\sum_{i=1}^N \{ \log(z-x_i)+\log(z+x_i) \}
+q \log z,
\quad z \in \O, \quad \x \in (\R_{\geq 0})^N,
\label{eqn:c-potentials}
\end{align}
where the latter contains a real parameter $q \in \R$.

Now we consider time evolution of the complex-valued potential
$\Phi_{D}$ by putting $\x$ be the driving process
$(\X^S(t))_{t \geq 0}$ of the multiple SLE
$(g_{D^{\eta}_t})_{t \geq 0}$ and map the
function $\Phi_{D}(\cdot, \X^{S}(t))$ by 
$(g_{D^{\eta}_t})_{t \geq 0}$.
Note that by (\ref{eqn:Ut2}), 
\[
X^S_i(t)=\lim_{\substack{z \to 0, \cr \eta_i(t) + z \in D^{\eta}_t}}
g_{D^{\eta}_t}(\eta_i(t)+z) 
=: g_{D^{\eta}_t}(\eta_i(t)),
\quad 1 \leq i \leq N, \quad t \geq 0, 
\]
and we will write
\[
\Phi_D(g_{D^{\eta}_t}(z), \X^S(t))
=\Phi_D \Big(g_{D^{\eta}_t}(z), g_{D^{\eta}_t}(\bmeta(t)) \Big), 
\quad t \geq 0,
\]
where 
$\bmeta(t):=(\eta_1(t), \dots, \eta_N(t))$
and
$g_{D^{\eta}_t}(\bmeta(t)):=
(g_{D^{\eta}_t}(\eta_1(t)), \dots, g_{D^{\eta}_t}(\eta_N(t)))$,
$t \geq 0$.
That is, we consider the complex-valued potentials
representing interactions between images by the multiple SLE
of a point $z$ in the domail $D$ and $N$ tips of a multi-slit,
$\eta_i(t), 1 \leq i \leq N$ at each time $t \geq 0$.
We obtain the following.
\begin{lem}
\label{thm:c-potential_time}
For $D=\H$ and $\O$, 
the increments of the complex-valued potentials are
given as follows, 
\begin{align}
& d \Phi_{\H}(g_{\H^{\eta}_t}(z), \X^{\R}(t)) =
- \sum_{i=1}^N \frac{\sqrt{\kappa} dB_i(t)}
{g_{\H^{\eta}_t}(z)-X^{\R}_i(t)} 
\nonumber\\
& \quad - \sum_{i=1}^N \Bigg( F^{\R}_i(\X^{\R}(t))
-4 \sum_{\substack{1 \leq j \leq N, \cr j \not=i}}
\frac{1}{X^{\R}_i(t)-X^{\R}_j(t)} \Bigg)
\frac{dt}{g_{\H^{\eta}_t}(z)-X^{\R}_i(t)} 
-\left(1 - \frac{\kappa}{4} \right) d \log g_{\H^{\eta}_t}'(z), 
\nonumber\\
& \hskip 10cm z \in \H^{\eta}_t, \quad t \geq 0,
\label{eqn:dPhiH}
\\
& d \Phi_{\O}(g_{\O^{\eta}_t}(z), \X^{\R_{\geq 0}}(t), q) =
- \sum_{i=1}^N 
\Bigg( \frac{1}{g_{\O^{\eta}_t}(z)-X^{\R_{\geq 0}}_i(t)}
- \frac{1}{g_{\O^{\eta}_t}(z)+X^{\R_{\geq 0}}_i(t)}
\Bigg) \sqrt{\kappa} d \widetilde{B}_i(t)
\nonumber\\
& \quad 
- \sum_{i=1}^N \Bigg[
F^{\R_{\geq 0}}_i(\X^{\R_{\geq 0}}(t)) 
- \Bigg\{ 4 \sum_{\substack{1 \leq j \leq N, \cr j \not=i}}
\left( \frac{1}{X^{\R_{\geq 0}}_i(t)-X^{\R_{\geq 0}}_j(t)}
+ \frac{1}{X^{\R_{\geq 0}}_i(t)+X^{\R_{\geq 0}}_j(t)} \right)
\nonumber\\
& \hskip 3cm 
+2 (1+2 \delta + q) \frac{1}{X^{\R_{\geq 0}}_i(t)}
\Bigg\}
\Bigg] 
\Bigg( \frac{1}{g_{\O^{\eta}_t}(z)-X^{\R_{\geq 0}}_i(t)}
- \frac{1}{g_{\O^{\eta}_t}(z)+X^{\R_{\geq 0}}_i(t)}
\Bigg) dt
\nonumber\\
& \quad
- 4 \delta \left( 1 - \frac{\kappa}{4} -q \right)
\frac{dt}{(g_{\O^{\eta}_{t}}(z))^2}
- \left(1 - \frac{\kappa}{4} \right) d \log g_{\O^{\eta}_t}'(z), 
\quad z \in \O^{\eta}_t, \quad t \geq 0,
\label{eqn:dPhiO}
\end{align}
where $g'_{D^{\eta}_t}(z) := d g_{D^{\eta}_t}(z)/dz$.
\end{lem}
\noindent{\it Proof} \,
By It\^o's formula, 
\begin{align*}
d \log (g_{D^{\eta}_t}(z) \pm X^{S}_i(t))
&= \frac{d g_{D^{\eta}_t}(z) \pm d X^{S}_i(t) }{g_{D^{\eta}_t}(z) \pm X^{S}_i(t)}
-\frac{\kappa dt}{( 2 g_{D^{\eta}_t}(z) \pm X^{S}_i(t))^2},
\quad t \geq 0, \quad 1 \leq i \leq N, 
\end{align*}
for $(D, S)=(\H, \R)$ and $(\O, \R_{\geq 0})$, and
$d \log g_{\O^{\eta}_t}(z) = d g_{\O^{\eta}_t}(z)/g_{\O^{\eta}_t}(z)$.
Put the multiple Loewner equation (\ref{eqn:mSLE1}), (\ref{eqn:mSLE2})
and the SDEs of their driving processes (\ref{eqn:SDE_A1}), 
(\ref{eqn:SDE_B1}). 
For $(D, S)=(\H, \R)$, we obtain the equation,
\begin{align*}
& d \Phi_{\H}(g_{\H^{\eta}_t}(z), \X^{\R}(t)) 
= 2 \sum_{1 \leq i, j \leq N} 
\frac{dt}{(g_{\H^{\eta}_t}(z)-X^{\R}_i(t)) (g_{\H^{\eta}_t}(z)-X^{\R}_j(t))}
\nonumber\\
& \qquad \quad
- \sum_{i=1}^N \frac{\sqrt{\kappa} dB_i(t)}{g_{\H^{\eta}_t}(z)-X^{\R}_i(t)}
- \sum_{i=1}^N \frac{F^{\R}_i(\X^{\R}(t))}{g_{\H^{\eta}_t}(z)-X^{\R}_i(t)}
- \frac{1}{2} \sum_{i=1}^N \frac{\kappa dt}{(g_{\H^{\eta}_t}(z)-X^{\R}_i(t))^2}.
\end{align*}
Here we use the equalities,
\begin{align*}
& \sum_{1 \leq i, j \leq N}
\frac{1}{(g-x_i)(g-x_j)}
= \sum_{i=1}^N
\frac{1}{(g-x_i)^2}
+\sum_{1 \leq i \not= j \leq N}
\frac{1}{(g-x_i)(g-x_j)}
\nonumber\\
& \qquad
=\sum_{i=1}^N
\frac{1}{(g-x_i)^2}
+2 \sum_{1 \leq i \not= j \leq N}
\frac{1}{(g-x_i)(x_i-x_j)},
\quad 1 \leq i \leq N.
\end{align*}
Since we obtain from (\ref{eqn:mSLE1}) 
the equality,
\[
d \log g_{\H^{\eta}_t}'(z)
= - 2 \sum_{i=1}^N \frac{dt}{(g_{\H^{\eta}_t}(z)-X^{\R}_i(t))^2},
\]
the equality (\ref{eqn:dPhiH}) is verified. 
For $(D, S)=(\O, \R_{\geq 0})$, use the equalities
\begin{align*}
& \sum_{1 \leq i \not= j \leq N}
\left(\frac{1}{g-x_i} + \frac{1}{g+x_i} \right)
\left(\frac{1}{g-x_j} + \frac{1}{g+x_j} \right)
\nonumber\\
& \qquad 
= 2 \sum_{1 \leq i \not= j \leq N}
\left(\frac{1}{g-x_i} + \frac{1}{g+x_i} \right)
\left(\frac{1}{x_i-x_j} + \frac{1}{x_i+x_j} \right),
\nonumber\\
& \left(\frac{1}{g-x_i} + \frac{1}{g+x_i} \right) \frac{1}{g}
= \left(\frac{1}{g-x_i} - \frac{1}{g+x_i} \right) \frac{1}{x_i},
\quad 1 \leq i \leq N.
\end{align*}
and
\[
d \log g_{\O^{\eta}_t}'(z)
= - 2 \Bigg[ 
\sum_{i=1}^N 
\left\{ \frac{1}{(g_{\O^{\eta}_t}(z)-X^{\R_{\geq 0}}_i(t))^2}
+ \frac{1}{(g_{\O^{\eta}_t}(z)+X^{\R_{\geq 0}}_i(t))^2}
\right\}
+ \frac{2 \delta}{(g_{\O^{\eta}_t}(z))^2}
\Bigg] dt.
\]
Then the equality 
(\ref{eqn:dPhiO}) is proved. \qed
\vskip 0.3cm

If we assume that $\X^{\R}(t)$
is given by the $(8/\kappa)$-Dyson model $(\Y^{\R}(t))_{t \geq 0}$
satisfying the SDEs (\ref{eqn:DysonB}), the second term
in RHS of (\ref{eqn:dPhiH}) vanishes.
For (\ref{eqn:dPhiO}), first we put
$a=1-\kappa/4$ to make the third term in RHS become zero.
Then if we assume that $\delta=\nu$ and $\X^{\R_{\geq 0}}(t)$
is given by the $(8/\kappa, \nu)$-Bru--Wishart process 
$(\Y^{\R_{\geq 0}}(t))_{t \geq 0}$
satisfying the SDEs (\ref{eqn:Bru_WishartB}), the second term
in RHS of (\ref{eqn:dPhiO}) vanishes.

We note that a multiple SLE driven by the Dyson model
(or the Bru--Wishart process) is absolutely continuous
with respect to multiple of independent SLEs 
(see Section 3 in \cite{Gra07}).
Then the original SLE and multiple SLEs share many common properties.
For example, if we define 
$\tau^{\eta}_z := \sup\{t > 0: z \in D^{\eta}_t \}$, then
$\tau^{\eta}_z < \infty$ a.s. for any $z \in D$ \cite{RS05}.
Hence we obtain the following statements.

\begin{prop}
\label{thm:martingales}
Assume that 
\begin{equation}
q=1 - \frac{\kappa}{4}, \quad
\delta=\nu,
\label{eqn:paraA}
\end{equation}
and 
define
\begin{align}
\cM_{\H}(z, t)
&=- \Phi_{\H}(g_{\H^{\eta}_t}(z), \Y^{\R}(t))
- \left(1-\frac{\kappa}{4} \right)
\log g_{\H^{\eta}_t}'(z),
\quad z \in \H^{\eta}_t, \, t \geq 0, 
\nonumber\\
\cM_{\O}(z, t)
&=- \Phi_{\O}(g_{\O^{\eta}_t}(z), \Y^{\R_{\geq 0}}(t);  1-\kappa/4)
- \left(1-\frac{\kappa}{4} \right)
\log g_{\O^{\eta}_t}'(z),
\quad z \in \O^{\eta}_t, \, t \geq 0.
\label{eqn:martingales}
\end{align}
Then for each point $z \in \D$, 
$\cM_D(z, t \wedge \tau^{\eta}_z)$, 
$D=\H$ and $\O$, 
provide local martingales with increments,
\begin{align*}
d \cM_{\H}(z, t)
&= \sum_{i=1}^N \frac{\sqrt{\kappa} dB_i(t)}
{g_{\H^{\eta}_t}(z)-Y^{\R}_i(t)},
\quad z \in \H^{\eta}_t, \quad t \geq 0,
\nonumber\\
d \cM_{\O}(z, t)
&= \sum_{i=1}^N 
\Bigg( \frac{1}{g_{\O^{\eta}_t}(z)-Y^{\R_{\geq 0}}_i(t)}
- \frac{1}{g_{\O^{\eta}_t}(z)+Y^{\R_{\geq 0}}_i(t)}
\Bigg) \sqrt{\kappa} d \widetilde{B}_i(t),
\quad z \in \O^{\eta}_t, \quad t \geq 0.
\end{align*}
\end{prop}
\vskip 0.3cm

\subsubsection{Stationary evolution of GFFs 
coupled with stochastic log-gases}
\label{sec:stationaryGFF}

Now we consider a coupling of
$(H_{D^{\eta}_t}(z))_{t \geq 0}$ and
some functional of $(\cM_D(z, \Y^S(t)))_{t \geq 0}$;
\[
H_{D^{\eta}_t}(z)+ \alpha \cF[\cM_D(z, t) ],
\quad z \in D^{\eta}_t, \quad t \geq 0,
\]
where $\cF[\, \cdot \,]$ denotes a functional
and $\alpha$ is a coupling constant.

Comparing Lemma \ref{thm:dG} and Proposition \ref{thm:martingales}
we observe the fact that
\[
d \Big\langle \Im \cM_D(z, \cdot), 
\Im \cM_D(w, \cdot) \Big\rangle_t
=- \frac{\kappa}{4} d G_{D^{\eta}_t}(z, w),
\quad z, w \in D^{\eta}_t, \quad t \geq 0, 
\]
for $(D,S)=(\H, \R)$ and $(\O, \R_{\geq 0})$.
Hence we put
\[
\cF[\, \cdot \,] =\Im [\, \cdot \,] \quad \mbox{and} \quad
\alpha=\frac{2}{\sqrt{\kappa}},
\]
and define the time-dependent system of Gaussian field,
\begin{align}
H_D(z, t)
&:= H_{D^{\eta}_t}(z)+\frac{2}{\sqrt{\kappa}} 
\Im \cM_D(z, t)
\nonumber\\
&= H_{D^{\eta}_t}(z)
- \frac{2}{\sqrt{\kappa}} \Im \Phi_D(g_{D^{\eta}_t}(z), \Y^{S}(t))
- \chi \log g_{D^{\eta}_t}'(z),
\label{eqn:GFF_BP1}
\end{align}
$z \in D^{\eta}_t, \Y^S(t) \in S^N, t \geq 0$,
with 
\begin{equation}
\chi=\alpha \left( 1 - \frac{\kappa}{4} \right)
=\frac{2}{\sqrt{\kappa}} - \frac{\sqrt{\kappa}}{2}.
\label{eqn:chi}
\end{equation}
Since $\log \zeta = \log |\zeta| + \sqrt{-1} \arg \zeta$
for $\zeta \in \C$, where $\arg \zeta$ is a priori 
defined up to additive multiple of $2 \pi$, 
$H_D(z, t)$ defined by (\ref{eqn:GFF_BP1}) 
with (\ref{eqn:c-potentials}) and (\ref{eqn:paraA}) 
is written as follows, 
\begin{align}
H_{\H}(z, t) 
&= H_{\H^{\eta}_t}(z) 
-\frac{2}{\sqrt{\kappa}} \sum_{i=1}^N \arg (g_{\H^{\eta}_t}(z)-Y^{\R}_i(t))
-\chi \arg g_{\H^{\eta}_t}'(z), \quad z \in \H^{\eta}_t, \quad t \geq 0,
\nonumber\\
H_{\O}(z, t) 
&= H_{\O^{\eta}_t}(z) 
-\frac{2}{\sqrt{\kappa}} \sum_{i=1}^N 
\Big\{ \arg (g_{\O^{\eta}_t}(z)-Y^{\R}_i(t))
+ \arg (g_{\O^{\eta}_t}(z)+Y^{\R}_i(t)) \Big\}
\nonumber\\
& \quad
-\chi \arg g_{\O^{\eta}_t}(z)
-\chi \arg g_{\O^{\eta}_t}'(z), \quad z \in \O^{\eta}_t, \quad t \geq 0.
\label{eqn:GPP_BP2}
\end{align}
Note that if we put $t=0$ in (\ref{eqn:GFF_BP1}), we have
\begin{align}
H_D(z, 0)
&= H_D(z)+\frac{2}{\sqrt{\kappa}} \Im \cM_D(z, 0)
\nonumber\\
&
= H_D(z)-\frac{2}{\sqrt{\kappa}} \Im \Phi_D(z, \y^S)
\nonumber\\
&= \begin{cases}
\displaystyle{
H_{\H}(z)-\frac{2}{\sqrt{\kappa}} 
\sum_{i=1}^N \arg (z- y_i^{\R}), 
}
& \mbox{for $H=\H$},
\cr
\displaystyle{
H_{\O}(z)-\frac{2}{\sqrt{\kappa}} 
\left\{ \sum_{i=1}^N \arg (z- y_i^{\R_{\geq 0}})
+\sum_{i=1}^N \arg (z+ y_i^{\R_{\geq 0}}) \right\}
- \chi \arg z,
}
& \mbox{for $H=\O$},
\end{cases}
\label{eqn:time0}
\end{align}
where $\y^S=\Y^S(0) \in \W_N(S)$, since 
$D^{\eta}_0=D$, 
$g_{D^{\eta}_0}(z)=g_0(z)=z$ and $g_{D^{\eta}_0}'(z)=g_0'(z)=1, z \in D$.

It was argued in \cite{SS13} that
a GFF on a subdomain of $D$ can be regarded as 
a GFF on $D$. 
Following it, we regard $H_{D}(\cdot, t), t >0$ as a GFF on $\H$
so that the pairing $\bra H_{D}(\cdot, t), f \ket,  t > 0$ with
$f \in \cC^{\infty}_{\rm c}(D)$ makes sense.

\begin{thm}
\label{thm:mainGFF1}
Let $\kappa \in (0, 4]$. 
Assume that $(D, S)=(\H, \R)$ or $(\O, \R_{\geq 0})$
and $(\Y^S(t))_{t \geq 0}$ is the $(8/\kappa)$-Dyson model
if $S=\R$ and the $(8/\kappa, \nu)$-Bru--Wishart process
if $S=\R_{\geq 0}$.
Then 
$\{H_D(z, t)\}_{z \in D}, t \geq 0$
is stationary in the sense that
\begin{equation}
\bra H_D(\cdot, t), f \ket \law= \bra H_D(\cdot, 0), f \ket
\quad 
\mbox{in $\P \otimes \rP$, \quad 
$\forall f \in \cC_{\rm c}^{\infty}(D)$ \quad
at each time $t \geq 0$}.
\label{eqn:stationary0}
\end{equation}
\end{thm}
\noindent{\it Proof} \,
For any test function $f \in \cC_{\rm c}^{\infty}(D) \subset \sD((-\Delta^{-1}))$,
we have 
\[
d \left\langle \left\bra \frac{2}{\sqrt{\kappa}} 
\Im \cM_D(\cdot, \cdot), f \right\ket \right\rangle_t
=- d E_t(f),
\]
where
\[
E_t(f) := \int_{D^{\eta}_T \times D^{\eta}_t}
f(z) G_{D^{\eta}_t}(z, w) f(w) d \mu(z) d \mu(w),
\]
which is called the {\it Dirichlet energy}.
Since $D^{\eta}_t := D \setminus \bigcup_{i=1}^N \eta_i(0, t]$
is decreasing, $E_t(f)$ is non-increasing in time $t \geq 0$.
This implies that
$\bra (2/\sqrt{\kappa}) \Im \cM_D(\cdot, t), f \ket, t \geq 0$ is a Brownian motion
such that we can regard $-E_t(f)$ as time.
Let $T \in (0, \infty)$.
Then $\bra (2/\sqrt{\kappa}) \Im \cM_D(\cdot, T), f \ket$ is normally distributed with
mean $\bra (2/\sqrt{\kappa}) \Im \cM_D(\cdot, 0), f \ket$ and variance
$-E_T(f)-(-E_0(f))=-E_T(f)+E_0(f)$.
On the other hand, the random variable
$\bra H_{D^{\eta}_T}, f \ket :=\bra H_D \circ g_{D^{\eta}_T}, f \ket$
is also normally distributed with mean zero
and variance $E_t(f)$ by the conformal invariance of GFF.
Since the random variable $\bra H_{D^{\eta}_T}, f \ket$
is conditionally independent of 
$\bra (2/\sqrt{\kappa}) \Im \cM_D(\cdot, T), f \ket$, the sum
\[
\bra H_D(\cdot, T), f \ket
=\bra H_{D^{\eta}_T}(\cdot), f \ket +
\left\bra \frac{2}{\sqrt{\kappa}} \Im \cM_D(\cdot, T), f \right\ket
\]
is a normal random variable with mean
$\bra (2/\sqrt{\kappa}) \Im \cM_D(\cdot, 0), f \ket$ and
variance $(-E_T(f)+E_0(f))+E_T(f)=E_0(f)$.
These values of mean and variance
coincide with those of
$\bra H_D(\cdot,0), f \ket
=\bra H_D(\cdot)+(2/\sqrt{\kappa} \Im \cM_D(\cdot, 0), f \ket$.
Since $T \in (0, \infty)$ is arbitrary, the statement is proved. 
\qed
\vskip 0.3cm

Theorem \ref{thm:mainGFF1} implies that 
for $(D, S)=(\H, \R)$ and $(\O, \R_{\geq 0})$, 
by coupling the Dirichlet boundary GFF
$H$ on $D$ with the stochastic log-gas
$(\Y^{S}(t))_{t \geq 0}$ on $S$
via the multiple SLE driven by $(\Y^{S}(t))_{t \geq 0}$, 
we have a new kind of 
family of stationary processes of 
GFF following the probability law $\P \otimes \rP$ 
on $D \times S$. 
At the initial time, the process starts from
\begin{align*}
H_D(\cdot, 0)
&= \begin{cases}
\displaystyle{
H_{\H}(\cdot)-\frac{2}{\sqrt{\kappa}} 
\sum_{i=1}^N \arg (\cdot - y_i^{\R}), 
}
& \mbox{for $H=\H$},
\cr
\displaystyle{
H_{\O}(\cdot)-\frac{2}{\sqrt{\kappa}} 
\left\{ \sum_{i=1}^N \arg (\cdot - y_i^{\R_{\geq 0}})
+\sum_{i=1}^N \arg (\cdot + y_i^{\R_{\geq 0}}) \right\}
- \chi \arg \cdot,
}
& \mbox{for $H=\O$}.
\end{cases}
\end{align*}
Then we let the boundary points evolve according
to the stochastic log-gas and,
at each time $t > 0$, we consider the 
GFF $H_D(\cdot)+u_D(\cdot, t)$ on $D$, where
\begin{align*}
u_D(\cdot, t)
&= \begin{cases}
\displaystyle{
-\frac{2}{\sqrt{\kappa}} 
\sum_{i=1}^N \arg (\cdot - Y_i^{\R}(t)), 
}
& \mbox{for $H=\H$},
\cr
\displaystyle{
-\frac{2}{\sqrt{\kappa}} 
\left\{ \sum_{i=1}^N \arg (\cdot - Y_i^{\R_{\geq 0}}(t))
+\sum_{i=1}^N \arg (\cdot + Y_i^{\R_{\geq 0}}(t)) \right\}
- \chi \arg \cdot,
}
& \mbox{for $H=\O$}.
\end{cases}
\end{align*}
By definition, the pair $(D, H_D(\cdot)+u(\cdot, t))$
is equivalent to
$(D^{\eta}_t, H_D(\cdot, t))$.
Here the GFF $H_{D}(\cdot, t)$ can be extended to
a GFF on $D$ \cite{SS13}.
The stationary process $(H_D(\cdot, t))_{t \geq 0}$
can be regarded as a generalization of the one
considered by Miller and Sheffield \cite{She16,MS16a},
and, in particular, the equivalence class whose representative is given 
by $(D, H_D(\cdot, 0))$ is a generalization of
the imaginary surfaces (the AC surfaces)
studied by them. 
We note that $(\H, H_{\H}(\cdot, 0)) \sim (\O, H_{\O}(\cdot, 0))$ 
in the sense of Definition \ref{thm:def_IS}. 

\subsection{Proof of Theorem \ref{thm:BM}}
\label{sec:A}

We recall the {\it Riesz--Markov--Kakutani theorem} \cite{Asa10}.
Let $\cH$ be a compact Hilbert space and 
write $\cB_{\cH}$ for the family of Borel sets in $\cH$.
Then the space of real-valued continuous functions denoted by
$\cC(\cH)$ is a real Banach space with respect to
the supremum norm $\|\cdot\|_{\infty}$.

\begin{df}
\label{thm:positive_function}
A linear functional $\ell: \cC(\cH) \to \C$ is positive, 
if for an arbitrary non-negative function $f \in \cC(\cH)$,
we have $\ell(f) \geq 0$.
\end{df}
\begin{prop}[Riesz--Markov--Kakutani theorem]
\label{thm:RMK}
Let $\ell: \cC(\cH) \to \C$ be a positive linear functional.
Then there exists a unique finite measure
$\bP$ on $(\cH, \cB_{\cH})$ such that
\[
\ell(f)=\int_{X} f(x) \bP(d x), \quad f \in \cC(\cH).
\]
Moreover, $\bP(\cH)=\|\ell\|$ holds.
\end{prop}

Let 
$\Sigma_a:=\sigma(\{\bra \cdot, g \ket_{\nabla} : g \in \cH_{-a}(D)\})$
be a $\sigma$-algebra of $\cH_a(D)$.
Then the following proposition is proved.

\begin{prop}
\label{thm:Minlos}
Let $\psi : W(D) \to \C$ be a continuous functional
of positive type such that $\psi(0)=1$.
Then for each $a>1/2$, there exists a probability measure
$\bP$ on $(\cH_a(D), \Sigma_a)$ such that
\begin{equation}
\psi(f)=\int_{\cH_a(D)} e^{\sqrt{-1} \bra h, f \ket_{\nabla}} \bP(d h),
\quad f \in \cH_{-a}(D).
\label{eqn:Minlos}
\end{equation}
\end{prop}
\noindent{\it Proof} \,
The proof consists of two steps. \\
\underline{\it Step 1.} \,
Let $\widehat{\R}:=\R \cup \{\infty \}$ be a one-point
compactification of $\R$. Then
\[
Q :=\widehat{\R}^{\N}=\{ h=(h_n)_{n \in \N} : h_n \in \widehat{\R}, n \in \N\}
\]
is a compact Hausdorff space.
The family of Borel sets in $Q$ is denoted by $\cB_Q$.
Given $h=(h_n)_{n \in \N} \in Q$, 
we assign a real-valued function $q_n$ for each $n \in \N$ by
\[
q_n(h) = \begin{cases}
h_n, & h_n \not= \infty,
\cr
0, & h_n=\infty.
\end{cases}
\]
Then, it can be verified that $q_n, n \in \N$ are Borel measurable.
We write the space of real-valued continuous functions on $Q$
as $\cC(Q)$, which is a real Banach space with
respect to the supremum norm.
Let $\cC_{\rm fin}(Q)$ be the collection of continuous functions on $Q$
that depend on finitely many $f_n$'s, that is,
\[
\cC_{\rm fin}(Q) :=
\{f \in \cC(Q): ^{\exists \! \! \!}N \in \N, ^{\exists \! \!}\{i_1, \dots, i_N \} \subset \N,
f=f(f_{i_1}, \dots, f_{i_N}) \}.
\]
By a simple argument, it can be verified that
$\cC_{\rm fin}(Q)$ is dense in $\cC(Q)$. 

Note that the space $\widehat{\cH}(D)$ of formal series
is isomorphic to $\R^{\N}$, it can be identified with
an open set in $Q$.
Let $\cD$ be a subspace in $\widehat{\cH}(D)$ defined by
$\cD:= \bigoplus_{n \in \N} \R u_n$.
With $h \in \cC_{\rm fin}(Q)$, $f \in \cD$, 
we associate a Borel measurable function on $Q$,
\begin{equation}
F_h(f) := \bra h, f \ket_{\nabla}
=\sum_{n \in \N} q_n(h) \bra u_n, f \ket_{\nabla},
\label{eqn:functionF}
\end{equation}
which is a finite sum.

For $N \in \N$ and $\{i_1, \dots, i_N\} \subset \N$, we set
$\cD_{\{i_1, \dots, i_N\}} :=
\bigoplus_{n=1}^N \R u_{i_n} \simeq \R^N$.
Then we apply the Bochner theorem (Theorem \ref{thm:Bochner}) to 
$\psi_{\{i_1, \dots, i_N\}} := \psi |_{\cD_{\{i_1, \dots, i_N\}}}$
and obtain a probability measure $\bP_{\{i_1, \dots, i_N\}}$
on $(\cD_{\{i_1, \dots, i_N\}}, \cB^N)$ such that
\[
\psi_{\{1_i, \dots, i_N\}}(f)
=\int_{\cD_{\{i_1, \dots, i_N\}}} 
e^{\sqrt{-1} F_h(f)} \bP_{\{i_1, \dots, i_N\}}(d h),
\quad f \in \cD.
\]
Using this family $\{\bP_{\{i_1, \dots, i_N\}} : N \in \N, \{i_1, \dots, i_N\} \subset \N \}$
of probability measures, we define a linear functional $\ell : \cC_{\rm fin}(Q) \to \C$ by
\[
\ell(\varphi) :=
\int_{\cD_{\{i_1, \dots, i_N\}}} \varphi(h_{i_1}, \dots, h_{i_N})
\bP_{\{i_1, \dots, i_N\}}(d h),
\quad \varphi \in \cC_{\rm fin}(Q).
\]
Here we have chosen, for each $\varphi \in \cC_{\rm fin}(Q)$, 
a finite set $\{i_1, \dots, i_N\} \in \N$
such that $\varphi$ depends on $h_{i_1}, \dots, h_{i_N}$.
Then it can be verified that the functional $\ell$ is
well-defined independent of the choice of such finite sets.
Moreover, it is extended to a positive functionals
on $\cC(Q)$. 
Therefore, the Riesz--Markov--Kakutani theorem ensures that
there exists a unique probability measure $\bP$ on $(Q, \cB_Q)$ such that
\[
\ell(\varphi)=\int_{Q} \varphi(h) \bP(d h), \quad
\varphi \in \cC(Q).
\]
In particular, if we take $\varphi(h)=e^{\sqrt{-1} F_h (f)} \in \cC_{\rm fin}(Q)$ for
$f \in \cD$, 
where $F_h$ is defined by (\ref{eqn:functionF}), 
we have
\begin{equation}
\psi(f)=\int_{Q} e^{\sqrt{-1} F_h(f)} \bP(d h), \quad f \in \cD.
\label{eqn:chi_f1}
\end{equation}
\noindent
\underline{\it Step 2.} \,
By assumption, $\psi$ is continuous.
Therefore, for an arbitrary $\varepsilon >0$, there exists
$\delta >0$ such that, if $\|f\|_{\nabla} < \delta$, then
$|1-\psi(f)| < \varepsilon$. 
Let us fix such $\varepsilon$ and $\delta$.
Then, in particular, we have
\[
\Re (\psi(f)) > 1-\varepsilon, \quad \|f\|_{\nabla} < \delta.
\]
Since $\psi$ is of positive type, we have $|\psi(f)| \leq \psi(0) =1$,
and in particular $\Re(\psi(f)) \geq -1, f \in W(D)$.
If $\|f\|_{\nabla} \geq \delta$, then we have
$-1 > 1-\varepsilon - 2 \delta^{-2} \|f\|_{\nabla}^2$.
Thus
\[
\Re(\psi(f)) > 1- \varepsilon - 2 \delta^{-2} \|f\|_{\nabla}^2,
\quad \|f\|_{\nabla} \geq \delta.
\]
The same inequality also hold when $\|f\|_{\nabla} < \delta$.
Therefore
\begin{equation}
\Re(\psi(f)) > 1 - \varepsilon -2 \delta^{-2} \|f\|_{\nabla}^2,
\quad f \in W(D).
\label{eqn:ineq1}
\end{equation}
In particular, if we set $f=\sum_{n=1}^N f_n u_n \in 
\cD_{\{1, \dots, N\}}$, we have
\[
\Re(\psi(f)) > 1-\varepsilon - 2 \delta^{-2} \sum_{n=1}^N f_n^2.
\]
For $\alpha >0$, we introduce a probability measure
$\rP_{\alpha, N}$ on $(\cD_{\{1, \dots, N\}}, \cB^N)$ as
\[
\rP_{\alpha, N}(d f)
=\prod_{n=1}^N \sqrt{\frac{\lambda_n^{2a}}{2 \pi \alpha}}
e^{-\lambda_n^{2a} f_n^2/2 \alpha} df_n,
\quad f=(f_1, \dots, f_N) \in \R^N,
\]
where $\{\lambda_n\}_{n \in \N}$ 
are eigenvalues of $-\Delta$ 
as given by (\ref{eqn:lambda}).
When we put the integral expression
of $\psi(f)$ (\ref{eqn:chi_f1}) with (\ref{eqn:functionF}) 
for $f \in \cD_{\{1, \dots, N\}}$ 
into LHS of the
inequality (\ref{eqn:ineq1}) and then 
integrate the both sides of it
with respect to $\rP_{\alpha, N}(d f)$, we obtain 
\begin{equation}
\int_{Q} e^{-(\alpha/2) \sum_{n=1}^N \lambda_n^{-2a} q_n(h)^2} \bP(d h)
> 1-\varepsilon - 2 \alpha \delta^{-2} \sum_{n=1}^N \lambda_n^{-2a}.
\label{eqn:ineqB}
\end{equation}
Now we take the limit $N \to \infty$.
Note that the integrand of LHS 
of (\ref{eqn:Minlos}) is supported on
$\cH_{a}(D) = \{ h=\sum_{n \in \N} h_n u_n : \sum_{n \in \N}
(\lambda_n^{-a} h_n)^2 < \infty\}$.
The sum in RHS of (\ref{eqn:ineqB}) is shown to converge
\[
C:= \lim_{N \to \infty} \sum_{n=1}^N \lambda_n^{-2a}
\sim \sum_{n=1}^{\infty} n^{-2 a} < \infty
\]
relying on the Weyl formula (Lemma \ref{thm:Weyl}) 
and the assumption $a > 1/2$.
Therefore, we see that
\[
\int_{\cH_a(D)} e^{-(\alpha/2) \|h\|_a^2} \bP (d h)
> 1 - \varepsilon - 2 \alpha \delta^{-2} C.
\]
At the limit $\alpha \to 0$, this gives
$\bP(\cH_a(D)) > 1 -\varepsilon$. 
Since $\varepsilon >0$ is arbitrary, we have
$\bP(\cH_a(D))=1$, which allows us to
restrict the measure $\bP$ onto $(\cH_a(D), \Sigma_a)$
and have 
\[
\psi(f)=\int_{\cH_a(D)} e^{\sqrt{-1} \bra h, f \ket_{\nabla}} \bP (d h),
\quad f \in \cD.
\]
In this expression, it is obvious that the domain for $f$ can be extended
to $\cH_{-a}(D)$, Therefore, the proof is complete. \qed
\vskip 0.3cm

By the definition (\ref{eqn:cE}), 
Theorem \ref{thm:BM} is concluded from
Proposition \ref{thm:Minlos} proved above.

\subsection{On the domain
of functions for Theorem \ref{thm:BM}}
\label{sec:B}

We have constructed a family of random variables
$\{\bra H, f \ket_{\nabla} : f \in \cE(D)^{\ast} \}$ on $D \subsetneq \C$
so that the assignment $f \mapsto \bra H, f \ket_{\nabla}$
is almost surely continuous.
We show here that, under certain conditions,
the domain of test functions for the random field $H$
can be extended from $\cE(D)^{\ast}$ to $W(D)$ 
if we give up its continuity.

\begin{prop}
\label{thm:extension1}
Let $\psi: W(D) \to \C$ be a continuous functional
of positive type such that $\psi(0)=1$.
Suppose that $\psi$ further satisfies the following assumptions.
\begin{description}
\item{\rm \bf (A.1)} \,
For an arbitrary $N \in \N$, the function
$\psi(\sum_{n=1}^N t_n u_n), t_n \in \R, n=1, \dots, N$
is of $\cC^2$-class.

\item{\em \bf (A.2)} \,
For an arbitrary $f \in W(D)$, the infinite series
\[
\left\{ \left.
\sum_{n=1}^N \sum_{m=1}^N
\frac{\partial^2 \psi(t_n u_n+t_m(1-\delta_{nm}) u_m)}
{\partial t_n \partial t_m} 
\right|_{t_n=t_m=0}
\bra u_n, f \ket_{\nabla} \bra u_m, f \ket_{\nabla} : N \in \N \right\}
\]
converges.
\end{description}
Then there exists a family of random variables
$\{\bra H, f \ket_{\nabla} : f \in W(D)\}$ such that
\begin{description}
\item{\rm (i)} \quad
$\bra H, f \ket_{\nabla} \in L^2(\cE(D), \bP)$ for $f \in W(D)$.

\item{\rm (ii)} \quad
$\bra H, af+bg \ket_{\nabla}=a \bra H, f \ket_{\nabla}
+b \bra H, g \ket_{\nabla}$
for $a, b \in \R, f, g \in W(D)$.

\item{\rm (iii)} \quad
If $f \in \cE(D)^{\ast}$, then $\bra H, f \ket_{\nabla}$ coincides
with that given by Theorem \ref{thm:BM}.

\item{\rm (iv)} \quad
The following is established
\begin{equation}
\psi(f)=\int_{\cE(D)} e^{\sqrt{-1} \bra h, f \ket_{\nabla}} \bP(d h)
\quad \mbox{for $f \in W(D)$}.
\label{eqn:BM2}
\end{equation}
\end{description}
\end{prop}
\noindent{\it Proof} \,
From the assumption {\bf (A.1)}, (\ref{eqn:BMeq1}) gives
\[
\left. 
\frac{\partial^2 \psi(t_n u_n+t_m(1-\delta_{nm}) u_m)}
{\partial t_n \partial t_m}
\right|_{t_n=t_m=0}
=- \int_{\cE(D)} q_n(h) q_m(h) \bP(d h), 
\quad n, m \in \N.
\]
Set
\[
\bra h^{(N)}, f \ket_{\nabla}
:= \sum_{n=1}^N q_n(h) \bra u_n, f \ket_{\nabla}, \quad
f \in W(D), \quad N \in \N.
\]
Then for $N>M$,
\begin{align*}
&\int_{\cE(D)} 
\left|
\bra h^{(N)}, f \ket_{\nabla}- \bra h^{(M)}, f \ket_{\nabla} 
\right|^2 \bP(d h)
\nonumber\\
& \quad = \sum_{n=M+1}^N \sum_{m=M+1}^N \bra u_n, f \ket_{\nabla} 
\bra u_m, f \ket_{\nabla} \int_{\cE(D)} q_n(h) q_m(h) \bP(d h)
\nonumber\\
&\quad = - \sum_{n=M+1}^N \sum_{m=M+1}^N
\left.
\frac{\partial \psi(t_n u_n+(1-\delta_{nm}) t_m u_m)}{\partial t_n \partial t_m}
\right|_{t_n=t_m=0}
\bra u_n, f \ket_{\nabla} \bra u_m, f \ket_{\nabla}.
\end{align*}
By the assumption {\bf (A.2)}, this converges to 0 as $N, M \to \infty$.
Therefore the sequence $\{ \bra H^{(N)}, f \ket_{\nabla} : N \in \N\}$ is
a Cauchy sequence in $L^2(\cE(D), \bP)$ and the limit
\[
\bra H, f \ket_{\nabla} := \lim_{N \to \infty} \bra H^{(N)}, f \ket_{\nabla} 
\in L^2(\cE(D), \bP)
\]
exists and (i) is proved. 
The linearlity (ii) is obvious. 
By construction (iii) is concluded.
Since $\psi$ is continuous, we have
\[
\psi(f)=\lim_{N \to \infty} \psi
\left( \sum_{n=1}^N u_n \bra u_n, f \ket_{\nabla} \right)
=\lim_{N \to \infty} 
\int_{\cE(D)} e^{\sqrt{-1} \bra h^{(N)}, f \ket_{\nabla}} \bP(d h),
\quad f \in W(D).
\]
We see that
\begin{align*}
&\left| 
\int_{\cE(D)} \left(e^{\sqrt{-1} \bra h, f \ket_{\nabla}}
-e^{\sqrt{-1} \bra h^{(N)}, f \ket_{\nabla}} \right) 
\right| \bP(d h)
\leq 
\int_{\cE(D)} \left| e^{\sqrt{-1} \bra h, f \ket_{\nabla}}
-e^{\sqrt{-1} \bra h^{(N)}, f \ket_{\nabla}} \right| \bP(d h)
\nonumber\\
&\qquad \leq 
\int_{\cE(D)} \left|
\bra h, f \ket_{\nabla} - 
\bra h^{(N)}, f \ket_{\nabla} \right| \bP(d h)
\leq
\left( \int_{\cE(D)} \left|
\bra h, f \ket_{\nabla}
- \bra h^{(N)}, f \ket_{\nabla} \right|^2 \bP(d h) \right)^{1/2}
\nonumber\\
& \hskip 7cm
\to 0 \quad \mbox{as $N \to \infty$}.
\end{align*}
This implies (iv).
Then the proof is complete. \qed

\vskip 1cm
\noindent{\bf Acknowledgements} \,
This manuscript was prepared for the 
mini course given at 
`Workshop on Probability and Stochastic Processes'
held at Orange County, Coorg, India, 
from 23rd to 26th February, 2020,
which was organized by 
the Indian Academy of Sciences, Bangalore.
The present author would like to thank
Rahul Roy for kind invitation 
and stimulating discussion during the workshop.
The lectures are based on the joint work with
Tomoyuki Shirai (Kyushu University), Hideki Tanemura (Keio University),
and Shinji Koshida (Chuo University). 
The author is grateful to them for fruitful collaborations.
The present study has been supported by
the Grant-in-Aid for Scientific Research (C) (No.19K03674),
(B) (No.18H01124), and (S) (No.16H06338) 
of Japan Society for the Promotion of Science (JSPS).


\end{document}